\documentclass[12pt]{article}
\usepackage{amsmath}
\usepackage{amssymb}
\usepackage{amscd}
\usepackage{epsfig}
\usepackage{amsthm}

\def\mini{\scriptsize}
\def\t{\tilde}

\newcommand{\Isom}{\operatorname{Isom}}
\def\ola{\overrightarrow}
\def\half{\frac{1}{2}}
\def\no{\noindent}
\newtheorem{dfn}{Definition}[section]
\newtheorem{rem}[dfn]{Remark}
\newtheorem{thm}[dfn]{Theorem}
\newtheorem{mainthm}[dfn]{Main Theorem}
\newtheorem{defn}[dfn]{Definition}
\newtheorem{notation}[dfn]{Notation}

\newtheorem{lem}[dfn]{Lemma}

\newtheorem{prob}[dfn]{Problem}
\newtheorem{prop}[dfn]{Proposition}
\newtheorem{ass}[dfn]{Assumption}

\newtheorem{cor}[dfn]{Corollary}
\newtheorem{conj}[dfn]{Conjecture}

\newtheorem{condition}[dfn]{Condition}
\newtheorem{ex}[dfn]{Example}
\newtheorem{example}[dfn]{Example}

\newtheorem{lis}[dfn]{List}

\def\proof{\par\medskip\noindent{\it Proof: }}

\def\lra{\longrightarrow}
\def\Lra{\Longrightarrow}
\def\ra{\rightarrow}
\def\Ra{\Rightarrow}
\def\CR{\curvearrowright}
\def\acts{\CR}

\def\embed{\hookrightarrow}

\def\back{\backslash}

\def\C{{\mathbb C}}
\def\R{{\mathbb R}}
\def\H{{\mathbb H}}
\def\Z{{\mathbb Z}}
\def\K{{\mathbb K}}
\def\M{{\mathcal M}}

\def\O{{\mathcal O}}

\def\Q{{\mathbb Q}}

\def\N{{\mathbb N}}

\def\F{{\mathbb F}}
\def\g{{\mathfrak g}}
\def\h{{\mathfrak h}}
\def\k{{\mathfrak k}}
\def\p{{\mathfrak p}}
\def\a{{\mathfrak a}}

\def\al{\alpha}
\def\be{\beta}
\def\ga{\gamma}

\def\Ga{\Gamma}
\def\de{\delta}
\def\De{\Delta}

\def\Del{\Delta}
\def\Si{\Sigma}
\def\si{\sigma}

\def\la{\lambda}
\def\La{\Lambda}

\def\>{\rangle}
\def\<{\langle}

\def\3{\ss}
\def\8{\infty}
\def\geo{\partial_{\infty}}
\def\tits{\partial_{Tits}}
\def\tangle{\angle_{Tits}}
\newcommand{\oa}[1]{\stackrel{\ra}{#1}}
\def\ol{\overline}
\def\ul{\underline}
\def\ov{\overrightarrow}

\def\mf{\mathfrak }

\begin{document}

\title{The generalized triangle inequalities in symmetric spaces and buildings
with applications to algebra}
\author{Misha Kapovich, Bernhard Leeb and John J. Millson}
\date{June 8, 2005}

\maketitle

\begin{abstract}
In this paper we apply our results on the geometry of polygons in
infinitesimal symmetric spaces, symmetric spaces and buildings, \cite{KapovichLeebMillson1, KapovichLeebMillson2},
to four problems in algebraic group theory.
Two of these problems are generalizations of the problems of finding the constraints
on the eigenvalues (resp. singular values) of a sum (resp. product)
when the eigenvalues (singular values) of each summand (factor) are fixed.
The other two problems are related to the nonvanishing of the structure
constants of the (spherical)
Hecke and representation rings associated with a split reductive algebraic group
over $\Q$ and its complex Langlands' dual.
We give a new proof of the
``Saturation Conjecture'' for $GL(\ell)$ as a consequence of our solution
of the corresponding ``saturation problem'' for the Hecke structure constants
for all split reductive algebraic groups over $\Q$.
\end{abstract}

\section{Introduction}
\label{introduction}

In this paper we will examine and generalize the algebra problems
listed immediately below
from the point of view of spaces of non-positive curvature.

Let ${\mathbb F}$ be either the field $\R$ or $\C$,
and let $\K$ be a nonarchimedean valued field
with discrete valuation ring ${\mathcal O}$ and the value group $\Z$.
For simplicity, let us consider here and throughout this Introduction a split
reductive group $\ul{G}$ over $\Q$,
see chapter \ref{sect:organizing algebra} for a more general discussion.

\begin{itemize}
\item {\bf Q1. Eigenvalues of a sum.} Set $G:=\ul{G}({\mathbb F})$, let $K$ be a maximal
compact subgroup of $G$.
Let $\g$ be the Lie algebra of $G$,
and let $\g=\k+\p$ be its Cartan decomposition.
Give necessary and sufficient conditions on
$\alpha,\beta,\gamma\in\p/Ad(K)$
in order that there exist elements $A,B,C\in\p$
whose projections to $\p/Ad(K)$ are
$\al,\beta$ and $\ga$, respectively,
so that
$$A+ B +C = 0.$$

\item {\bf Q2. Singular values of a product.}
Let $G$ and $K$ be the same as above.
Give necessary and sufficient conditions on
$\alpha,\beta,\gamma\in K\backslash G/K$
in order that there exist elements
$A,B,C\in G$
whose projections to $K\backslash G/K$
are $\alpha,\beta$ and $\gamma$, respectively,
so that
$$ABC = 1.$$

\item  {\bf Q3. Invariant factors of a product.}
Set $G:=\ul{G}(\K)$ and $K:=\ul{G}(\O)$.
Give necessary and sufficient conditions on
$\alpha,\beta,\gamma\in
K\backslash G/K$
in order that there exist elements
$A,B,C\in G$
whose projections to
$K\backslash G/K$
are $\alpha,\beta$ and $\gamma$, respectively,
so that
$$
ABC = 1.$$

\item {\bf Q4. Decomposing tensor products.}
Let $\ul{G}^{\vee}$ be the Langlands' dual group of $\ul{G}$, see
Definition \ref{Langlandsdual}.
Give necessary and sufficient conditions on
highest weights $\al,\beta,\ga$
of irreducible representations $V_{\al}$, $V_{\beta}$, $V_{\ga}$
of $G^{\vee}:=\ul{G}^{\vee}(\C)$ so that
\begin{equation*}
(V_{\alpha} \otimes V_{\beta} \otimes V_{\gamma})^{G^{\vee}} \neq 0.
\end{equation*}
\end{itemize}

\begin{notation}
Throughout the paper we will denote by $Sol(\mathbf{Qi}, G)$ the sets of triples
$(\al,\be,\ga)$ which are solutions of
Problems {\bf Qi}, $i=1,...,4$ for the group $G$.
\end{notation}

We refer the reader to section \ref{linalg} for
{\em linear algebra} reformulations of the problems {\bf Q1, Q2, Q3}
in the case when $\ul{G}=GL(m)$
and to section \ref{geo=Hecke} for reformulation of {\bf Q3} in terms of the Hecke ring.

The above problems have a long history which is described in
detail in \cite{Fulton}, we briefly discuss it in section
\ref{linalg}.

To state our main results we need several notation. The quotient
spaces $\p/Ad(K)$ and $K\backslash G/K$ in Problems {\bf Q1} and
{\bf Q2} are naturally identified with the Euclidean Weyl chamber
$\De$ inside a Cartan subspace ${\mf a}\subset\p$.

The double coset space appearing in Problem {\bf Q3} is a more
subtle discrete object, namely the intersection of $\De$ with the
{\em cocharacter lattice} $L_{\ul{G}}\subset \a$ (see section
\ref{duality} for the precise definition).  The appearance of the
dual group $G^{\vee}$ in Problem {\bf Q4} is explained by the fact
that the cocharacter lattice of $\ul{G}$ is the same as the
character lattice of the dual group $\ul{G}^\vee$, see section
\ref{duality}. We let $\La$ denote the lattice
$$
\La:= \{(\al,\be,\ga)\in L_{\ul{G}}^3 : \al+\be+\ga\in
Q(R^\vee)\}.
$$
Here  $Q(R^\vee)$ is the coroot lattice of $\ul{G}$.
The importance of the lattice $\La$ comes from the inclusion
$$
Sol(\mathbf{Q3}, \ul{G}(\K))\subset \La.
$$
Next we need the number $k_R$, the {\em saturation factor} associated
with the root system $R$ of the group $\ul{G}$; it is computed in
the table \ref{tab} in section \ref{sec:satfactor}. For a split
rank $\ell$ simple algebraic group $\ul{G}$  over $\K$ these
numbers are defined as follows. Let $R$ be the associated root
system, let $\al_1,...,\al_\ell$ be the simple roots and $\theta$
be the highest root:
\begin{equation}\label{highestroot0}
\theta= \sum_{i=1}^\ell m_i \al_i.
\end{equation}
Then $k_R$ is the least common multiple of $m_1,...,m_\ell$.

\medskip
With these notation, the main result of this paper is the following
chain of relations between the problems {\bf Q1--Q4}:

\begin{mainthm}\label{mainthm}
1.
$$
Sol(\mathbf{Q1}, \ul{G}(\F))= Sol(\mathbf{Q2},
\ul{G}(\F))\supset Sol(\mathbf{Q3}, \ul{G}(\K))\supset
Sol(\mathbf{Q4}, \ul{G}^\vee(\C)).
$$

2.
$$
Sol(\mathbf{Q2}, \ul{G}(\F)) \cap k_R \cdot \La \subset
Sol(\mathbf{Q3}, \ul{G}(\K)) \subset Sol(\mathbf{Q2},
\ul{G}(\F)) \cap \La$$

3. For the root systems $B_2=C_2$ and $G_2$ we have
$$
Sol(\mathbf{Q2}, \ul{G}(\F))\cap \La\ne Sol(\mathbf{Q3},
\ul{G}(\K))
$$
and
$$
Sol(\mathbf{Q3}, \ul{G}(\K)) \ne Sol(\mathbf{Q4},
\ul{G}^\vee(\C)).
$$

4. Even if $\ul{G}$ is not necessarily split, for $G= \ul{G}(\F)$,
the solution set $Sol(\mathbf{Q1}, G)$ depends only on the
(finite) Weyl group of $G$. It is a convex homogeneous polyhedral
cone $D_3(G/K)$ contained in $\De^3$. A description of the
linear inequalities defining this cone can be found in section
\ref{polyhedron}.

5. Even if $\ul{G}$ is not necessarily split, for $G= \ul{G}(\K)$,
the solution set $Sol(\mathbf{Q3}, G)$ depends only on the affine
Weyl group of $G$.
\end{mainthm}

\begin{rem}
The equivalence of Problems {\bf Q1} and {\bf Q2} (known as
Thompson's Conjecture) for some classical groups (including
$GL(m)$) was proved by Klyachko in \cite{Klyachko2}, for all
complex semisimple groups by Alexeev, Meinrenken and Woodward in
\cite{AMW}, and in general case by Evens and Lu in the recent
paper \cite{EL}, using different methods.
\end{rem}

\begin{rem}
The above saturation constants $k_R$ are not necessarily the
smallest possible, for instance, in \S \ref{semigroup} we prove
that for the group $\ul{G}$ with the root system $G_2$,
$$
Sol(\mathbf{Q2}, \ul{G}(\F)) \cap 2 \cdot \La \subset
Sol(\mathbf{Q3}, \ul{G}(\K)) $$

\end{rem}

\begin{rem}
It is well-known,  that
$$
Sol(\mathbf{Q4}, \ul{G}(\C))\subset Sol(\mathbf{Q1},
\ul{G}(\C))\cap \La,
$$
see Theorem \ref{vectorsandngons} in the Appendix. In the
forthcoming paper \cite{KM} it is shown that
$$
Sol(\mathbf{Q3}, \ul{G}(\K))\cap k_R\cdot \La \subset
Sol(\mathbf{Q4}, \ul{G}(\C)).
$$
\end{rem}

\medskip
Despite of the algebraic appearance of the results of the Main
Theorem, its main source is {\em geometry}, the sole exception is
the relation between {\bf Q3} and {\bf Q4} which is proven via
{\em Satake correspondence}. In this paper, we reformulate the
algebra problems {\bf Q1}--{\bf Q3} as {\em geometric problems}
which are special cases of a geometric question raised and studied
in \cite{KapovichLeebMillson1, KapovichLeebMillson2}:

We fix a Euclidean Coxeter complex $(E,W_{aff})$ with the Euclidean Weyl chamber
$\De\cong E/W_{sph}$, and consider nonpositively curved metric spaces $X$
with geometric structures {\em modeled on} $(E,W_{aff})$.
The spaces of this kind which we are interested in are
symmetric spaces of nonpositive curvature,
their infinitesimal versions ({\em infinitesimal symmetric spaces},
see section \ref{geometries}),
and Euclidean buildings.
For such spaces $X$ there is a notion of $\De$-length
for oriented geodesic segments
which reflects the anisotropy of $X$. This leads to the following problem:

\begin{itemize}
\item
{\bf GTI: Generalized Triangle Inequalities.}
Give necessary and sufficient conditions on $\al, \beta, \ga\in\De$
in order that there exists a geodesic triangle in $X$
with $\De$-side lengths $\al,\beta$ and $\ga$.
\end{itemize}

\no We refer to the solution set of {\bf GTI} as $D_3(X)$.

As explained in section \ref{Geodesic polygons}, for a symmetric
space $X$, the problem {\bf GTI} is equivalent to the Singular
Value Problem {\bf Q2}. For an infinitesimal symmetric space $\p$,
it is equivalent to the Eigenvalues of a Sum Problem {\bf Q1}:

\begin{thm}
\label{equ} Suppose that $\ul{G}$ is a reductive algebraic group
over $\R$ and $\F=\C$ or $\F=\R$; let  $X=G/K$ be the symmetric space
corresponding to $G=\ul{G}(\F)$ and $\p$ be the infinitesimal
symmetric space. Then
$$
Sol(\mathbf{Q1}, G)=D_3(\p), \quad Sol(\mathbf{Q2}, G)=D_3(X).
$$
\end{thm}

The Invariant Factor Problem {\bf Q3} corresponds to the case when
$X$ is a Euclidean building and in this case the inequalities
defining $D_3(X)$ give only {\em necessary} conditions on
$(\al,\beta,\ga)$ to solve the algebraic problem {\bf Q3}. The key
difference with the problems {\bf Q1, Q2} is lack of homogeneity
in Euclidean buildings $X$: {\em Not all vertices of $X$ are
special} (except for the buildings of type $A_\ell$). This
eventually explains why for the Euclidean building $X$
corresponding to some groups $\ul{G}(\K)$ we have the inequality
$$
D_3(X)\cap \La \ne Sol(\mathbf{Q3},\ul{G}(\K)).
$$

We refer the reader to section \ref{Geodesic polygons}
for justification of the correspondence between algebraic and geometric problems.

Despite of the results of this paper, the structure of the
solution sets $Sol(\mathbf{Q3})$, $Sol(\mathbf{Q4})$ remains
unclear. It was known for a long time that $Sol(\mathbf{Q4})$ is
a semigroup. In contrast, we will see in section \ref{semigroup}
that the set $Sol({\bf Q3})$ is not a semigroup for the root
systems $B_2=C_2$ and $G_2$.

\medskip
The following result, a consequence of a recent theorem in logic
of M.~C.~Las\-kow\-ski \cite{Laskowski}, building on work of
S.~Kochen \cite{Kochen}, reveals the general structure of the
subset
$$
Sol({\bf Q3}, \ul{G}(\K))\subset L_{\ul{G}}^3.
$$
Define a subset of $L_{\ul{G}}^3$ to be {\em elementary} if it is
the set of solutions of a finite set of linear inequalities with
integer coefficients and a finite set of linear congruences. Then
Laskowski proves  the following

\begin{thm}[M.~C.~Laskowski, \cite{Laskowski}]
\label{las} There is an integer $N=N_{\ul{G}}$, depending only on
$\ul{G}$, such that for any nonarchimedean Henselian valued field
$\K$ with value group $\Z$ and residue characteristic greater than
$N$ we have
$$
Sol({\bf Q3}, \ul{G}(\K))\subset L_{\ul{G}}^3 \text{~~is a finite
union of elementary
 sets}.$$
\end{thm}

\begin{rem}
It follows from Theorem \ref{mainthm}, Part 5, that Theorem
\ref{las} holds for {\em all} complete nonarchimedean valued
fields $\K$. This is because all split groups $\ul{G}$, as $\K$
varies, have the same affine Coxeter group. Hence the set
$Sol({\bf Q3}, \ul{G}(\K))$ does not depend on $\K$. Therefore,
once the above statement is true for one of them it is true for
all of them.
\end{rem}

\begin{prob}
\label{lasproblem} 1. Find the corresponding inequalities and
congruences.

2. Is there an analogue of Laskowski's theorem that describes
the set Sol({\bf Q4})? If so find the inequalities and congruences
in this case as well.
\end{prob}

We can solve the Problem \ref{lasproblem} and thus the Problem {\bf Q3}
for the following groups:

\begin{thm}
\label{typeA} If $\ul{G}$ is covered by $SL(m)$ or $GL(m)$ then
$$
Sol({\bf Q3}, \ul{G}(\K))=Sol({\bf Q1}, \ul{G}(\F))\cap \La.
$$
In other words, say for the group $\ul{G}=SL(m)$, the solution set
$Sol({\bf Q3}, SL(m,\K))$ is the set of triples
$(\alpha,\beta,\gamma)$ in $D_3(Sym_m)\cap \Z^{3m}$ such that
$\alpha + \beta + \gamma\in Q(R^\vee)$. Here $Sym_\ell$ is the
infinitesimal symmetric space of zero trace symmetric $m\times m$
matrices and $Q(R^\vee)$ is the coroot lattice of $SL(m)$.
\end{thm}

A similar statement holds for a more general class of groups, see
Theorem \ref{atype}. For example, the conclusion of part (a) holds
for the groups $SL(\ell, \mathbb{D})$ where $\mathbb{D}$ is a
division algebra over $\K$, see Tableau des \'echelonages \cite[p.
29-30]{BT}.

We conclude by pointing out that for the group $GL(m)$ it was
known since the 1960-s that
$$
Sol(\mathbf{Q3}, GL(m))= Sol(\mathbf{Q4}, GL(m))
$$
and therefore we can use Theorem \ref{typeA} to give a new proof
of the Saturation Conjecture for $GL(m)$ (the theorem of Knutson
and Tao), see section \ref{saturationforgl}.

\medskip
This paper is organized as follows.

In chapter \ref{datum} we review root systems for algebraic
reductive groups and discuss Coxeter groups.


In chapter \ref{sect:organizing algebra} we set up the general algebra problem {\bf R}
which generalizes the setting of the Problems {\bf Q1--Q3}. We then discuss
in detail the parameter spaces for the Problems {\bf Q1--Q4} and their
mutual relation.

In chapter \ref{settingupgeometry} we first convert the problem
{\bf R} into an abstract geometry problem about existence of
polygons with the prescribed generalized side-lengths (Problem
\ref{linalgprob}). Next, we introduce a class of metric spaces
(metric spaces modeled on Euclidean Coxeter complexes) and restate
the abstract geometry problem for three classes of such spaces:
infinitesimal symmetric spaces, nonpositively curved symmetric
spaces and Euclidean buildings. We introduce the notion of refined
and coarse (the $\De$-length) generalized side-lengths for the
geodesic polygons in such metric spaces. For the infinitesimal
symmetric spaces and symmetric spaces, the problem of existence of
polygons with the prescribed $\De$-side lengths is adequate for
solving the corresponding algebra problems ({\bf Q1} and {\bf
Q2}), but in the case of buildings it is not. This establishes the
equalities and inclusions
$$
Sol({\bf Q3})\subset Sol({\bf Q1})= Sol({\bf Q2})=D_3(X)
$$
in Theorem \ref{mainthm} as well as Part 4 and 5 in this theorem.

In chapter \ref{gau} we describe the solution (given in
\cite{KapovichLeebMillson1, KapovichLeebMillson2}) to the problem
{\bf GTI} of existence of polygons with the prescribed $\De$-side
lengths in the above classes of metric spaces $X$. We also discuss
the relation of this solution to symplectic and Mumford quotients.
We describe the system of generalized triangle inequalities for
$X$ and give an explicit set of inequalities in the case of root
system of type $B_2$. {\em The reader who is not interested in the geometric ideas behind
the results of this paper can skip this chapter except for the definition of the mapping $\Phi_\psi$
for buildings which will be used in chapter \ref{k-computation}}.

In chapter \ref{building} we show that the geometry problem {\bf
GTI} solved in the previous chapter is not adequate (in the
building case) for solving the algebra problem {\bf Q3}. The key
tool here is the procedure of {\em folding} of geodesic triangles
in Euclidean buildings into apartments. The main application of
the examples constructed in this section is to establish the first
inequality in Part 3 of Theorem \ref{mainthm}.

In chapter \ref{k-computation} we show that in the case of
$A$-type root systems, the  solution of the geometry problem {\bf GTI}
given in chapter \ref{gau} solves the algebra problem {\bf Q3} as
well:
$$
Sol({\bf Q3}) =D_3(X)\cap \La.
$$
More generally, we establish the existence of the saturation factors
$k=k_{inv\; fact}$ and compute these numbers; modulo
multiplication by $k$ the geometry problem {\bf GTI} is equivalent
to the algebra problem {\bf Q3}. This establishes Part 2 of
Theorem \ref{mainthm}.

In chapter \ref{3&4} we reformulate problem {\bf Q3} in terms of
{\em spherical Hecke rings}, compare the algebraic problems {\bf
Q3} and {\bf Q4}, establish the inclusion
$$
Sol({\bf Q4}, \ul{G}^\vee(\C))\subset  Sol({\bf Q3}, \ul{G}(\K))
$$
(see Part 1 of Theorem \ref{mainthm}) and give our new proof of
the saturation conjecture for $GL(m)$.

In the Appendix, chapter \ref{append}, we relate the solution set of the problem {\bf Q4}
to Mumford quotients and show that the former forms a semigroup. Although these results
are known to experts, we include them for the sake of completeness.

\medskip
{\bf Acknowledgments.} During the work on this paper the first
author was visiting the Max Plank Institute (Bonn), he was also
supported by the NSF grants DMS-99-71404 and DMS-02-03045. The
first and the third authors were supported by the Mathematics
Department of the University of T\"ubingen during their stay there
in May of 2002. The third author was supported by the NSF grant
DMS-01-04006. The authors gratefully acknowledge support of these
institutions. We are grateful the referees of this paper for their
useful remarks. The authors are grateful to L.~ Ein, S.~Kumar,
T.~Haines, R.~Hemmecke, C.~Laskowski, J.~Stembridge, E.~Vinberg
and C.~Woodward for useful conversations. We would especially like
to thank Steve Kudla for suggesting we use the Satake transform to
connect Problems {\bf Q3} and {\bf Q4} and Jiu-Kang Yu for his
help with chapter \ref{3&4} and in particular for providing us
with the proof of Lemma \ref{basic}.

When we told G.~Lusztig of our Theorem \ref{q4->q3} establishing
the inclusion
$$
Sol({\bf Q4}, \ul{G}^\vee(\C))\subset Sol({\bf Q3}, \ul{G}(\K)),
$$
for local nonarchimedean valued fields $\K$, he informed us that
although he had not known the result before our message, it was an
easy consequence of his work in \cite{Lusztig2}. We should say
that our proof depends in an essential way on Lusztig's paper
\cite{Lusztig2}. We use his change of basis formulas, Lemma
\ref{expansionlemma}, and his realization  that the coefficients
in one of those formulas were Kazhdan-Lusztig polynomials for the
affine Weyl group, Lemma \ref {KazhdanLusztigestimate}.
Subsequently, alternative proofs of Theorem \ref{q4->q3} were
given by Tom Haines in \cite{Haines2} and by Kapovich and Millson
in \cite{KM}.

\tableofcontents

\newpage
\section{Roots and Coxeter groups}
\label{datum}

Let $\ul{G}$ be a reductive algebraic group over a field $\F$ and
$\ul{T}$ be a split torus in $\ul{G}$.
Our goal in this section
is to describe {\em the root datum}  associated to
the pair $(\ul{G}, \ul{T})$. The reader will find the definition of root
datum in \cite[\S 1]{Springer}.

\subsection{Split tori over $\F$}

We recall that the
algebraic group $\mathbb{G}_m$ is the affine algebraic group with
coordinate ring $\F[S,T](ST-1)$ and comultiplication $\Delta$ given by
$\Delta(T)=T\otimes T,\Delta(S)=S\otimes S$.

\begin{dfn}\label{splittorus}
An affine algebraic group $\ul{T}$
defined over $\F$ is a split torus of rank $l$ if it is isomorphic
to the product of $l$ copies of $\mathbb{G}_m$.
\end{dfn}

A character of an
algebraic group $\ul{T}$ defined over $\F$ is a morphism of algebraic
groups from $\ul{T}$ to $\mathbb{G}_m$. The product of two characters
and the inverse of a character are characters and accordingly the
set of characters of $\ul{T}$ is an abelian group denoted by $X^*(\ul{T})$.

\begin{lem}
Suppose that $\ul{T}$ is a split torus over $\F$ of rank $l$. Then the
character group of $\ul{T}$ is a lattice (i.e. free abelian group) of
rank $l$.
\end{lem}
\proof
We have $\F[\ul{T}]\cong \F[T_1,T_1^{-1},\cdots,T_l,T_l^{-1}]$. A character
of $\ul{T}$ corresponds to a Hopf algebra morphism from $\F[\mathbb{G}_m]$
to $\F[T_1,T_1^{-1},\cdots,T_l,T_l^{-1}]$. Such a morphism is determined
by its value on $T$. This value is necessarily a grouplike element
(this means $\Delta(f)=f\otimes f$). However the grouplike elements
of $\F[T_1,T_1^{-1},\cdots,T_l,T_l^{-1}]$ are the monomials in the
$T_i$'s and their inverses. The exponents of the monomial give the
point in the lattice.
\qed

\begin{cor}
$$Hom( \mathbb{G}_m, \mathbb{G}_m)\cong \Z.$$
\end{cor}

We note that the previous isomorphism is realized as follows. Any
Hopf-algebra homomorphism
$\psi$ of coordinate rings is of the form $T\to T^n$ for some integer $n$.
Then the above isomorphism sends $\psi$ to $n$.

\begin{dfn}
A cocharacter or a one-parameter (algebraic) subgroup of $\ul{T}$ is a morphism
$\phi:\mathbb{G}_m\to T$. The set of cocharacters of $\ul{T}$ will
be denoted $X_*(\ul{T})$.
\end{dfn}

\begin{lem}
Suppose that $\ul{T}$ is a split torus over $\F$ of rank $l$. Then the
cocharacter group of $\ul{T}$ is a free abelian group of
rank $l$.
\end{lem}
\proof
We have $\F[\ul{T}]\cong \F[T_1,T_1^{-1},\cdots,T_l,T_l^{-1}]$. A cocharacter
of $\ul{T}$ corresponds to a Hopf algebra morphism $\psi$ from
$\F[T_1,T_1^{-1},\cdots,T_l,T_l^{-1}]$
to $\F[\mathbb{G}_m]$. Such a morphism is determined
by its value on $T_1,\cdots, T_l$. Then $\psi$ corresponds to the lattice
vector $(m_1,\cdots,m_l)$ where $T^{m_i} = \psi (T_i)$.
\qed

We define an integer-valued pairing $\<~,~\>$ between characters and cocharacters
as follows. Suppose $\phi$ is a cocharacter of $\ul{T}$ and $\chi$ is a character.
Then $\chi\circ \phi \in Hom( \mathbb{G}_m, \mathbb{G}_m)\cong \Z$. We define
$\<\chi,\phi\>$ to be the integer corresponding to $\chi\circ \phi$.

We now describe two homomorphisms that will be useful in what follows.
Let $T_e(\ul{T})$ be the Zariski tangent space of $\ul{T}$ at the identity $e$.

\begin{dfn} \label{derivativeofacocharacter}
We define $\Phi:X_*(\ul{T})\lra T_e(\ul{T})$ by
$$\Phi(\lambda)= \lambda^{\prime}(1).$$
Here $1$ is the identity of $GL(1,\F)$ and $\lambda^{\prime}(1)$ denotes the derivative
at $1$.

We also define $\Phi^{\vee}:X^*(\ul{T})\lra T^*_e(\ul{T})$ by
$$\Phi^{\vee}(\lambda)= d\lambda|_e.$$
\end{dfn}

\begin{rem}
The character and cocharacter groups $X^*(\ul{T}), X_*(\ul{T})$
are multiplicative groups, the trivial (co)character
will be denoted by $1$. However, we will use the embeddings
$\Phi$ and $\Phi^{\vee}$ to identify them with additive groups.
This will be done for the most part in chapters
\ref{sect:organizing algebra} and \ref{3&4}.
\end{rem}

\subsection{Roots, coroots and the Langlands' dual}
\label{sect:Langlandsdual}

The reductive group $\ul{G}$ picks out a distinguished finite subset
of $X^*(\ul{T})$, the {\em relative root system} $R=R_{rel}(\ul{G}, \ul{T})$.
A character of $\ul{T}$ is a {\em root}
if it occurs in the restriction of the adjoint representation of $\ul{G}$
to $\ul{T}$. We let $Q(R)$ denote the subgroup of
$X^*(\ul{T})$ generated by $R$
and define $V:= Q(R)\otimes \R$.

We recall

\begin{dfn}\label{splitgroup}
The algebraic group $\ul{G}$ is {\em split over $\F$} if it has a maximal torus
$\ul{T}$ defined over $\F$, which is split.
\end{dfn}

From now on we assume $\ul{G}$ is split
over $\F$ and $\ul{T}$ is a maximal torus as in the above definition.

It is proved in \cite[\S 2]{Springer},
that $R\subset V$ satisfies the axioms of a root system.
Moreover in the same section it is proved
that to every root $\alpha \in R$ there is an associated coroot
$\alpha^{\vee} \in X_*(\ul{T})$ such that $\<\alpha, \alpha^{\vee}\> = 2$.
We let $R^{\vee}$ denote the resulting set of coroots, let $Q(R^{\vee})$ be
the subgroup of $X_*(\ul{T})$ they generate and
$V^{\vee}:=Q(R^{\vee})\otimes \R$. The root and coroot system $R$ and $R^{\vee}$ determine
(isomorphic) finite Weyl groups $W, W^{\vee}$
which acts on $V^{\vee}$ and $V$ respectively.
The action of the generators $s_{\alpha}, s_{\alpha^{\vee}}$ of the group $W, W^{\vee}$ on
$X^*(\ul{T})$ and $X_*(\ul{T})$ are determined by the formulae:
$$
s_{\alpha}(x):= x - \<x,\alpha^{\vee}\> \alpha \ \text{and}\
s_{\alpha^{\vee}}(u):= u - \<u,\alpha\> \alpha^{\vee}.
$$

We then have \cite[\S 2]{Springer}:

\begin{prop}
The quadruple $\Psi(\ul{G},\ul{T}):= (X^*(\ul{T}),R, X_*(\ul{T}), R^{\vee})$
is a {\em root datum}.
\end{prop}

\begin{dfn}
Let $\Psi = (X,R,X^{\vee},R^{\vee})$
and
$\Psi^{\prime}=(X^{\prime},R^{\prime},(X^{\prime})^{\vee},(R^{\prime})^{\vee})$
 be root data.
Then an {\em isogeny}  from $\Psi^{\prime}$ to $\Psi$ is a homomorphism $\phi$ from
$X^{\prime}$ to $X$ such that $\phi$ is injective with finite cokernel.
Moreover we require that $\phi$ induces a bijection from $R^{\prime}$ to $R$
and the transpose of $\phi$ induces a bijection of coroots.
\end{dfn}

Now suppose $f:\ul{G}\lra \ul{G}^{\prime}$ is a covering of algebraic groups.
If $\ul{T}$ is a maximal torus in $\ul{G}$ then its image
$\ul{T}^{\prime}$ is a maximal torus in
$\ul{G}^{\prime}$.
The induced map on characters gives rise to an isogeny of root data, denoted $\Psi(f)$.

Conversely, we have \cite[Theorem 2.9]{Springer}:

\begin{thm}
\begin{enumerate}
\item For any root datum $\Psi$ with reduced root system there exist a
connected split (over $\F$) reductive group $\underline{G}$ and a maximal
split torus $\ul{T}$ such that $\Psi = \Psi(\ul{G}, \ul{T})$. The pair
$(\ul{G},\ul{T})$ is unique up to isomorphism.

\item Let $\Psi = \Psi(\ul{G}, \ul{T})$ and $\Psi^{\prime} =
\Psi(\ul{G}^{\prime}, \ul{T}^{\prime})$ and $\phi$ be an isogeny from
$\Psi^{\prime}$ to $\Psi$. Then there is a covering $f:\ul{G}
\lra \ul{G}^{\prime}$ with the image of $\ul{T}$ equal to $\ul{T}^{\prime}$
such that $\phi = \Psi(f)$.
\end{enumerate}
\end{thm}

Before stating the next definition we need a lemma which we leave to
the reader.

\begin{lem}
If $(X,R,X^{\vee},R^{\vee})$ is a root datum so is $(X^{\vee},R^{\vee},X,R)$.
\end{lem}

We now have

\begin{dfn}
\label{Langlandsdual}
Let $\ul{G}$ be a (connected) split reductive group over $\Q$. Let
$\Psi=\Psi(\ul{G},\ul{T})=$ $(X,R,X^{\vee},R^{\vee})$ be the root datum
of $(\ul{G},\ul{T})$. Then the {\em Langlands dual}
$\ul{G}^{\vee}$ of $\ul{G}$ is the unique (up to isomorphism) reductive group
over $\Q$ which has the root datum $\Psi^{\prime}= (X^{\vee},R^{\vee},X,R)$.
\end{dfn}

In fact we will need only the complex points $G^{\vee}:= \ul{G}^{\vee}(\C)$
of $\ul{G}^{\vee}$ in what follows. We will accordingly abuse notation and
frequently refer to $G^{\vee}$ as the Langlands dual of $\ul{G}$. One has
$$
(\ul{G}^{\vee})^{\vee}\cong \ul{G}.
$$

\subsection{Coxeter groups}

In this section we review the properties of Coxeter groups, we
refer the reader for a more thorough discussion to \cite[Section
4.2]{Humphreys} and \cite{Bourbaki}. Let $R\subset V^*$ be a root
system of rank $n$ on a real Euclidean vector space $V$. We do not
assume that $n$ equals the dimension of $V$. Note that in the case
of semisimple Lie algebras, the space $V$ will be a Cartan
subalgebra $\a\subset \g$ with the Killing form.

We will sometimes identify $V$ and $V^*$ using the metric.
Let $Q(R)\subset V^*$ denote {\em the root lattice}, i.e. the integer span of $R$.
This subgroup is a lattice in $Span_{\R}(R)\subset V$, it
is a discrete free abelian subgroup of rank $n$.
Given a subgroup $\La\subset \R$ we define a collection ${\mathcal H}={\mathcal H}_{R,\La}$
of hyperplanes (called {\em walls}) in $V$ as the set
$$
{\mathcal H}:=\{ H_{\al,t}=\{v\in V: \al(v) =t\}, t\in \La, \al\in R \}.
$$
In this paper we will be mostly interested in the case when $\La$ is either $\Z$ or $\R$,
but much of our discussion is more general.
We define an {\em affine Coxeter group} $W_{aff}= W_{R,\La}$ as the group generated
by the reflections $w_H$ in the hyperplanes $H\in {\mathcal H}$.
The only reflection hyperplanes of the reflections in $W_{aff}$ are the elements of
${\mathcal H}$.
The pair $(E, W_{aff})$ is called a {\em Euclidean Coxeter complex}, where $E=V$
is the Euclidean space. The {\em vertices} of the Coxeter complex
are points which belong to the transversal intersections of $n$ walls in ${\mathcal H}$.
(This definition makes sense even if $n<\dim(V)$, only in this case there will be continuum of
vertices even if $\La=\Z$.) If $W_{aff}$ is trivial, we declare {\em each point} of $E$ a vertex.
We let $E^{(0)}$ denote the vertex set of the Coxeter complex.

\begin{defn}
An {\em embedding of Euclidean Coxeter complexes} is a map $(f,\phi): (E,W)\to (E',W')$,
where $\phi: W\to W'$ is a monomorphism of Coxeter groups and $f: E\to E'$
is a $\phi$-equivariant affine embedding.
\end{defn}

Let $L_{trans}$ denote the translational part of $W_{aff}$. If $\La=\Z$ then
$L_{trans}$ is the {\em coroot lattice} $Q(R^{\vee})$ of $R$.
In general, $L_{trans}= Q(R^{\vee})\otimes \La$.
The linear part $W_{sph}$ of $W_{aff}$ is a finite Coxeter group acting on $V$, it is
called a {\em spherical Coxeter group}.
The stabilizer of the origin $0\in E$ (which we will regard as a base-point $o\in E$)
in $W_{aff}$ maps isomorphically onto $W_{sph}$.
Thus $W_{aff}= W_{sph}\ltimes L_{trans}$.
A vertex of the Euclidean Coxeter complex
is called  {\em special} if its stabilizer in $W_{aff}$ is isomorphic to
$W_{sph}$. We let $E^{(0),sp}$ denote the set of special vertices of $E$.

\begin{rem}
The normalizer $N_{aff}$ of $W_{aff}$ (in the full group $V$ of translations on $E$) acts
transitively on the set of special vertices. The vertex set $E^{(0)}$ of
the complex $(E, W_{aff})$ contains $N_{aff}\cdot o$, but typically it is strictly
larger that. Moreover, in many cases $E^{(0)}$ does not form a group.
\end{rem}
We recall that the weight group $P(R)$ and the coweight group $P(R^{\vee})$
are defined by
\begin{align*}
 P(R)= & \{\lambda \in V^*: \la(v) \in \Z, \forall v\in R^{\vee}\},\\
 P(R^{\vee})= & \{v\in V: \al(v) \in \Z, \forall \al\in R\}.
\end{align*}

\begin{rem}
In the case when $n<\dim(V)$ our definition of
weights is different from the one in \cite{Springer}.
\end{rem}

Again, $P(R)$ and $P(R^{\vee})$ are lattices provided that $n=\dim(V)$,
otherwise they are  nondiscrete abelian subgroups of $V$.
We have the inclusions
$$
Q(R^{\vee})\subset P(R^{\vee}), \quad Q(R)\subset P(R).
$$
The normalizer $N_{aff}$ equals $P(R^{\vee}) \otimes \La$.

The spherical Coxeter groups $W_{sph}$ which appear in the above construction
act naturally on the sphere at infinity $S=\geo E$; the pair $(S,W_{sph})$ is called
a {\em spherical Coxeter complex}. The definitions of {\em walls}, {\em vertices}, etc., for
Euclidean Coxeter complexes generalize verbatum to the spherical complexes.
We will use the notation $\De_{sph}\subset S$ for the spherical Weyl chamber,
$\De_{sph}$ is the ideal boundary of the Euclidean Weyl chamber $\De\subset E$
(i.e. a fundamental domain for the action $W_{sph}\acts E$, which is bounded by walls).

\ From our viewpoint, the Euclidean Coxeter complex is a more fundamental object
than a root system. Thus, if the root system $R$
was not reduced, we replace it with a reduced root system $R'$ which has the same group
$W_{aff}$: If $\al, 2\al\in R$ we retain the root $2\al$ and eliminate the root $\al$.
We will assume henceforth that the root system $R$ is reduced.

\bigskip
{\bf Product decomposition of Euclidean Coxeter complexes}.
Suppose that $(E, W_{aff})$ is a Euclidean Coxeter complex associated with
the reduced root system $R$, let $R_1,...,R_s$ denote the decomposition of
$R$ into irreducible components. Accordingly, the Euclidean space $E$ splits as
the metric product
$$
E=E_0 \times \prod_{i=1}^s E_i,
$$
where $E_i$ is spanned by $R_i^{\vee}$, $1\le i\le s$.
This decomposition is invariant under the
group $W_{aff}$ which in turn splits as
$$
W_{aff}= \prod_{i=1}^s W_{aff}^i,
$$
where  $W_{aff}^i=W_{R_i,\La}$ for each $i=1,...,s$;
for $i=0$ we get the trivial Coxeter group $W^0_{aff}$.
The group $W_{aff}^i$ is the image of $W_{R_i,\La}$ under the natural
embedding of affine groups $Aff(E_i) \lra \prod_{i=1}^s Aff(E_i) = Aff(E)$.

Analogously, the spherical Coxeter group $W_{sph}$ splits as the direct product
$$
W_{sph}^0\times ...\times W_{sph}^s$$
(where $W_{sph}^0=\{1\}$). The Weyl chamber
$\De$ of $W_{sph}$ is the direct product of the Weyl chambers
$\De_0\times \De_1\times ...\times\De_s$, where $\De_i$ is a Weyl chamber of $W^i_{sph}$ and
$\De_0=E_0$. Similarly, the normalizer $N_{aff}$ of $W_{aff}$ splits as
$$
N_{aff}= V_0\times \prod_{i=1}^s N_{aff}^i,
$$
where $V_0$ is the vector space underlying $E_0$ and $N_{aff}^i$ is the
normalizer
of $W_{aff}^i$ in the group of translations of $E_i$. Note that for each $i=1,...,s$
the groups $W_{aff}^i$ and $N_{aff}^i$ act as lattices on $E_i$.
We observe that the vertex set of the complex $(E,W_{aff})$ equals
$$
E_0\times E_1^{(0)} \times... \times E_s^{(0)},
$$
where $E_i^{(0)}$ is the vertex set of the complex $(E_i,W_{aff}^i)$. Similarly,
the set of special vertices $E^{(0),sp}$ of $E$ equals
$$
E_0\times E_1^{(0),sp} \times... \times E_s^{(0),sp}.
$$

\begin{rem}
Our discussion of Coxeter groups was somewhat nongeometric; a more geometric
approach would be to start with an affine Coxeter
group and from this get root systems, etc.
\end{rem}

\section{The first three algebra problems and the parameter spaces $\Sigma$ for
$K\backslash \overline{G} /K$}
\label{sect:organizing algebra}

We will see in this chapter that the problems {\bf Q1}--{\bf Q3}
for reductive algebraic groups $G$ can
reformulated as special cases of a single algebraic problem as follows.
There is a group $\overline{G}$ (closely associated to $G$)
which contains $K$, a maximal bounded subgroup of $G$. The conditions
of fixing $\alpha$, $\beta$ and $\gamma$ in problems {\bf Q1}--{\bf Q3}
will amount to fixing three double cosets in $\Si=K\backslash \overline{G}/K$.
The problems {\bf Q1}--{\bf Q3} will be then reformulated as:

\begin{itemize}
\item {\bf Problem R($\ol{G}$)}: Find necessary and sufficient conditions on
$\al, \be, \ga \in \Si$ in order that there
exist $A, B, C\in \ol{G}$ in the double cosets represented by $\al, \be, \ga$ resp.,
such that $A\cdot B \cdot C = 1$.
\end{itemize}

\no We will now describe the groups $\overline{G}$ and $K$ for the problems {\bf Q1}--{\bf Q3}.
The main part of this chapter will then be occupied with describing the double coset spaces
$\Si=K\backslash \overline{G} /K$.
In section \ref{summ} we will also prove that the problem {\bf R($\ol{G}$)} agrees with the
Problem {\bf Q1} from the Introduction (the equivalence will be clear for two other problems).

\begin{enumerate}
\item For the Problem {\bf Q1}:
For $\F=\R$ or $\C$, let $\ul{G}$ be a connected reductive algebraic group over $\mathbb{R}$,
$G:= \ul{G}(\F)$ be a real or complex Lie group with Lie algebra $\mathfrak{g}$. Pick a
maximal compact subgroup $K$ of $G$. Let $\k$ denote the Lie algebra of $K$.
Then we have the orthogonal decomposition (with respect to the Killing form)
$$
\mathfrak{g} = \mathfrak{k} \oplus \mathfrak{p}.$$
We let  $\overline{G}$ be the {\em Cartan motion group}
$\overline{G}= K \ltimes \mathfrak{p}$.

\item For the problem {\bf Q2}: $\ul{G}$ and $K$ are the same as in 1, but now we take $\ol{G}= G=\ul{G}(\F)$.

\item For the problem {\bf Q3}:
Let $\K$ be a complete nonarchimedean valued field with a discrete valuation $v$ and the value group
$\Z=v(\K)\subset \R$.
Let $\mathcal{O}$ denote the ring of elements in $\mathbb{K}$ with nonnegative valuation.
Let $\ul{G}$ be a connected reductive algebraic group over $\K$, $\ol{G}:=G:= \ul{G}(\K)$
and $K:=G(\mathcal{O})$.
\end{enumerate}

\no In order to relate  Problems {\bf Q1}--{\bf Q3} with geometry
we will have to compute the double coset spaces $\Sigma = K\back \ol{G}/K$.
Moreover we will show that the parameter spaces $\Sigma$
admit canonical linear structures:
In the first two problems we use a Cartan decomposition $G=KAK$,
then  $\Sigma$ is described as the (positive) Weyl chamber $\De$ in $\a$,
the Lie algebra of $A$. In the case of the third problem, $\Si$ is identified with
the intersection $L_{\ul{G}} \cap \Delta\subset \a=X_*(\ul{T})\otimes \R$,
where $\ul{T}$ is a maximal $\K$-split torus of $\ul{G}$
and $L_{\ul{G}}$ is {\em extended cocharacter lattice} of $\ul{T}$ (see section
\ref{tough}).

\subsection{The generalized eigenvalues of a sum problem {\bf Q1} and the parameter space
$\Sigma$ of $K$-double cosets}\label{summ}

We continue with the notation of previous section, in particular,
$\ol{G}= K\ltimes\mathfrak{p}$. We now describe
the parameter space $\Sigma=K\backslash \ol{G}/K$.
We choose a Cartan subspace $\mathfrak{a}
\subset \mathfrak{p}$  (i.e. a maximal subalgebra of $\mathfrak{p}$,
necessarily abelian and reductive). We recall that the positive
Euclidean Weyl chamber $\Delta$ is the dual cone in
$\mathfrak{a}$ to the cone of positive restricted roots in
$\mathfrak{a}^*$.

\begin{lem}
The inclusion $\iota:\Delta \to K\ltimes
\mathfrak{p}$ given $\iota(\sigma)=(1,\sigma)$ induces a bijection onto
$\Sigma$.
\end{lem}
\proof
We first observe that since $(k_1,x)\cdot (k_2,0) = (k_1k_2,x)$, every double
coset in $K\backslash \ol{G}/K$ has a representative of the form $(1,x)$. Since
$(k,0)\cdot (1,x) \cdot (k,0)^{-1} = (1, Ad\, k (x))$ the lemma is clear.
\qed

We complete the proof of equivalence of {\bf R}($\ol{G}$) with Problem {\bf Q1} by observing
that for $g_i=(k_i,x_i)$, $i=1, 2, 3$,
\begin{equation}
\label{add}
(k_1,x_1)\cdot (k_2,x_2) \cdot (k_3,x_3) =
(k_1 k_2 k_3, x_1 +  Ad k_1 x_2 + Ad (k_1 k_2) x_3 ).
\end{equation}

We see from (\ref{add}) that if $g_1 g_2 g_3=1$ then putting
$a=x_1, b= Ad k_1 x_2, c= Ad (k_1 k_2) x_3$ we find
$$
a+b+c= 0. $$

Thus we find representatives $a, b$ and $c$ in $\p$ of the orbits
$Ad(K)(x_i), 1\leq i \leq 3$ whose sum equals zero
(i.e. we have a solution of the Eigenvalues of a Sum Problem).
Conversely, if $a, b, c\in \p$ solve the Eigenvalues of a Sum Problem
then
$$
g_1\cdot g_2\cdot g_3=(1,a)\cdot (1,b) \cdot (1,c) = (1,0)$$
and hence $g_1, g_2, g_3$ solve the double coset problem {\bf R}($\ol{G}$).

Thus we may find $A$, $B$ and $C$ in $\overline{G}$ in the required $K$-double cosets
with $A\cdot B \cdot C = 1$ if and only if there exist $a,b,c \in \p$ with
$a + b + c = 0$ in the required $Ad\, K$-orbits in $\p$.

\subsection{The generalized singular values of a product and the parameter space
$\Sigma$ of $K$-double cosets}

Let $\ol{G}:=G=\ul{G}(\F)$ and consider the Cartan decomposition $G=KAK$ where $A$ is the maximal
abelian subgroup of $\ol{G}$, whose Lie algebra is the Cartan subalgebra
${\mf a}\subset \p$ from the previous section.
We will identify $A$ and ${\mf a}$ via the exponential map, let
$A_\Delta:=\exp(\Delta)\subset A$.
Then we get the refinement of the Cartan decomposition:
$G=K A_\Delta K$.
It is now clear that the inclusion $\iota:\Delta\to G$ given by $\iota(\sigma)= \exp(\sigma)$
induces a bijection onto $\Si= K\backslash G/K$.
In what follows we will parameterize
double cosets in $\Si$ by their logarithms, which are vectors in $\Delta$.

\subsection{The generalized invariant factor problem and the parameter space $\Sigma$
of $K$-double cosets}
\label{tough}

The goal of this section is to give an explicit description of
the parameter space $\Sigma=\Si_{inv\; fact}$
for the Invariant Factor Problem {\bf Q3}. We start by noticing that the
discrete valuation $v: \K\to \Z$ admits a splitting $\kappa: \Z\to \K$:
it is given by sending $1\in \Z$ to the uniformizer $\pi\in \K$.
Throughout this section, $\ul{G}$ will denote a connected
reductive algebraic group over $\K$, which has  the (relative) rank $l$. Let
$\ul{T}\subset \ul{G}$ be a maximal $\K$-split torus, and let $X^*(\ul{T})$
and $X_*(\ul{T})$ be the lattices of characters and cocharacters. We denote by
$\ul{P}\subset \ul{G}$ be a minimal parabolic subgroup such that
$(\ul{P},\ul{T})$ is a parabolic pair, \cite[pg. 127]{Cartier}.

This data determines a  (relative) root system $R_{rel}$
and hence a (relative) finite Weyl group $W_{sph}$.
Let $\ul{Z}\subset \ul{G}$ and $\ul{N}\subset \ul{G}$
be respectively the centralizer and normalizer
(over $\K$) of the algebraic group $\ul{T}$.
As usual we let $G, N, Z$ and
$T$ denote the groups of $\K$-points of the corresponding algebraic
groups. We let $K$ denote $\ul{G}(\O)$ and $B:=\ul{T}(\mathcal{O})=K\cap T$.
We will frequently refer to $G$ as a {\em nonarchimedean reductive Lie group}.

As in the previous section,
we have to discuss the Cartan decomposition $G= K\cdot \Ga \cdot K$
 of the group $G$ (where $\Ga$ will be an appropriate subset of $Z$).
In the split case this decomposition will have the form
$G= K\cdot A_{\Delta} \cdot K$,
where $A_{\Delta}\subset A$ and $A$ will be a subgroup of $T$ so that
$A\times B=T$.

\medskip
We can now describe $\Sigma_{inv\; fact}(G)=K\back G/K$, the parameter space
for the generalized invariant factors problem for the nonarchimedean
reductive Lie group $G$.

{\bf Split Case.}
We first consider the case when $\ul{G}$ is split over $\K$
since the description
is more transparent in this case.
In this case there is an algebraic group  defined
and split over $\Q$ and in fact a group scheme over $\Z$ such that $\ul{G}$ is
obtained from it by
extension of scalars. We will abuse notation and also denote the above
group-scheme over $\Z$  by $\ul{G}$ as well.
Under this assumption the torus $\ul{T}$ is also defined and split over $\Q$,
hence the group of  real points of $\ul{T}$ is well-defined.
We thus consider the real torus $\ul{T}(\R)$ with the {\em real} Lie algebra $\a$.
The map $\Phi$ of Definition \ref{derivativeofacocharacter} (with $\F =\R$)
gives us an identification
$$
X_*(\ul{T}) \otimes \R \cong \a.
$$
Then the group $X_*(\ul{T})$ is identified with a lattice
$L:= L_{\ul{G}}$ in $\a$. Similarly we use the map $\Phi^{\vee}$
to identify $X^*(\ul{T})$ with a lattice in the dual space $\a^*$.
The  root system $R$ sits naturally in the dual space $\a^*$ and the
corresponding finite Weyl group
$W_{sph}$ acts on $\a$ and $\a^*$ in the usual way;
$W_{sph}$ leaves invariant the lattices $X_*(\ul{T})$ and $X^*(\ul{T})$.
Let $\De\subset E$ be the (positive) Weyl chamber of the group $W_{sph}$,
determined by the parabolic subgroup $\ul{P}$.
Our goal is to prove that the space $K\back G/K$ can be identified with the ``cone''
$$
\Delta_{L} =\Delta \cap L.
$$

By definition, a cocharacter $\phi$ is a homomorphism from the group
$\mathbb{G}_m(\K)$ to the group
$T$: note that the map $\Psi: X_*(\ul{T})\to T$ given by
$\Psi(\phi) = \phi(\pi)$ is injective (where $\pi$
is a fixed  uniformizer). We define the subgroup $A$ by
$$
A = \Psi(X_*(\ul{T})).$$

The following lemma is immediate (it suffices to prove it for $GL(1)$).

\begin{lem}
The group $T$ is a direct product $T= A\times B$.
\end{lem}

We let  $A_{\Delta_{L}}\subset A$ be the image of $\Delta_{L}= X_*(\ul{T}) \cap \Delta$ under
the map $\Psi$.

\begin{lem}\label{thecartandecomposition}
(See \cite[page 51]{Tits}.) We have the Cartan decomposition
$$
G= K \cdot A_{\Delta_{L}} \cdot K,$$
i.e., each element $g \in G$ has a
unique representation as a product $k_1 \cdot a \cdot k_2$,
where $a\in A_{\Delta_{L}}$ and $k_i\in K$.
\end{lem}

\begin{cor}
The inclusion $\iota: A_{\Delta_{L}} \subset A\subset T \subset G$ induces a bijection\newline
$A_{\Delta_{L}} \stackrel{\iota}{\lra} K \backslash G /K$.
\end{cor}

Now that we have identified the double coset space
$$
\Si_{inv\; fact}(G)=K \backslash G /K$$
with $\De_L$, it is clear that
the generalized invariant factor problem {\bf Q3} in the Introduction
agrees with the double coset problem {\bf R}($G$). However it may not
be clear at this point why these problems agree with
Problem {\bf P3} in case $G=GL(\ell)$.

To see this we note that we can give a second description of the
vector $\Phi(\lambda)$ if we choose coordinates in $\ul{T}$. Let
$\chi_1,\cdots, \chi_n$ be a basis for the group $X^*(\ul{T})$
and define $\chi: T\lra \K^n$ to be the map with components $\chi_j, 1\leq j
\leq n$.
We observe that we may use the derivative of $\chi_j$ at $e$ to obtain
coordinates on $T_e(\ul{T})$ and accordingly we obtain a map
$\dot{\chi}: \mathfrak{a} \lra \R^n$. We leave the following lemma
to the reader.

\begin{lem}
$\dot{\chi}(\Phi(\lambda)$ is the (integral) vector in $\R^n$ obtained by
applying the valuation $v$ to the coordinates of $\chi(\lambda(\pi))$.
\end{lem}

\medskip
{\bf General Case.} We now consider the general case when $\ul{G}$ is not necessarily split.
Most of the above discussion remains valid,
however the Cartan decomposition for $G$ has a slightly different form and for this reason
 we {\em cannot} use the cocharacter lattice as $L_{\ul{G}}$.
We no longer
can talk about {\em real points} of $\ul{T}$, so we do not have the interpretation of $\a$ as
a tangent space. We {\em define}
$\a$ as $X_*(\ul{T}) \otimes \R$. We set $V:= \a$.

\medskip
Let $E$ denote the affine space underlying $V$; as before, $\De\subset V$ is the positive
Weyl chamber for $W_{sph}$. Note that $X^*(\ul{Z})\subset X^*(\ul{T})$ is a subgroup
of finite index.

\begin{defn}
\label{nu}
Define the homomorphism
$\nu: Z \to V$ by the formula (see \cite{Tits})
\begin{equation}
\label{eq:nu}
\<\chi,\nu(z)\>=-v(\chi(z)), \forall z\in Z, \chi\in X^*(\ul{Z}).
\end{equation}
\end{defn}
\no Note that on the left hand side of the above formula, the pairing of
$\chi\in X^*(\ul{Z})\subset X^*(\ul{T})$ with $\nu(z)$
is coming from the pairing $\<\cdot , \cdot \>$ between $X_*(\ul{T})$ and $X^*(\ul{T})$.
On the right hand side of the equation (\ref{eq:nu}), the action of $\chi$ on $z$
is coming from the fact that each character $\chi\in X^*(\ul{Z})$ defines a character
of the group $Z=\ul{Z}(\K)$, i.e. we have $\chi(z)\in \K^{\times}$.

Let $Z_c$ denote the kernel of $\nu$ and let $L_{\ul{G}}$ denote the image
of $\nu$ in $V$.  Then $Z_c$ is a maximal compact subgroup of $Z$, see
\cite[pg. 135]{Cartier}.

\begin{defn}
\label{coL}
We will refer to the group $L:=L_{\ul{G}}:=\nu(Z)$ as the {\em extended cocharacter
lattice} of  $\ul{T}$.
\end{defn}
We have the inclusions
$$
X_*(\ul{T})=Hom(X^*(\ul{T}), \Z) \subset L_{\ul{G}}\subset Hom(X^*(\ul{Z}), \Z)
$$
with equality in the case when $\ul{G}$ is split over $\Q$. The quotient group $N/Z_c$
operates on the affine space $E$ faithfully, through a discrete group $\tilde{W}$.
The finite Weyl group $W_{sph}$ is the linear part of $\tilde{W}$ (the stabilizer of the origin)
and we have $\tilde{W}=W_{sph}\ltimes L_{\ul{G}}$.

We refer to \cite[section 1.2]{Tits} for more details.
Let $\Delta_L:= \Delta\cap L$ and $Z_{\Delta}:=\nu^{-1}(\Delta)$. Then one has:

\begin{thm}
\cite[section 3.3.3]{Tits}
\label{carta}
$G= K Z K$ and the map $K\back G/K \to V$ given by
$$
KzK \mapsto \nu(z), z\in Z_\Delta,
$$
descends to a bijection $K\back G /K \to \Delta_L$.
\end{thm}

This proves that in the nonsplit case we can also identify $\Delta_L$ with the parameter space
$\Si_{inv\; fact}(G)=K\back G/K$ for the Problem {\bf Q3}.

\subsection{Comparison of the parameter spaces for the four algebra problems}
\label{duality}

We can now compare the four algebra problems {\bf Q1--Q4} stated in the Introduction.
We have already compared the parameter
spaces for {\bf Q3} and {\bf Q4} in the Introduction.
To compare the first three problems we assume that we are given groups
$\ul{G}_1, \ul{G}_2, \ul{G}_3$ which are
connected reductive algebraic groups (possibly all the same), where the first
two are defined over $\R$, and $\ul{G}_3$ is over $\F_{(3)}=\K$,
where $\K$ is a nonarchimedean valued field with valuation $v$ and value group $\Z$.
Suppose that the fields $\F_{(1)}, \F_{(2)}$ are either $\R$ or $\C$.  We assume
that all three Lie groups $G_i=\ul{G}_i(\F_{(i)})$ have the same (relative)
rank $l$ and isomorphic (relative) Weyl groups
$W=W_{sph}$ acting on the appropriately chosen vector space $\a$. Let $\De\subset \a$
be a (positive) Weyl chamber. Let $L=L_{\ul{G}_3}$ be the extended cocharacter
lattice, see Definition \ref{coL}.
Then the parameter spaces $\Si_1, \Si_2$ for the Problems {\bf Q1}
and {\bf Q2} (for the groups $G_1, G_2$ resp.) are exactly the same: they are equal to $\De$.
The parameter space $\Si_{inv\; fact}(G_3)$ is the intersection $\De_L$ of the lattice $L$
with $\De$.

Note that in the special case when $\ul{G}=\ul{G}_i$, $i=2, 3$, is split over $\Q$,
$L=X_*(\ul{T})$ and we get
a canonical  isomorphism
$$
X_*(\ul{T})\otimes \R\to \a,
$$
described in the previous section (here $\a$ is the {\em real} Lie algebra of $\ul{T}(\R)$).
This isomorphism
identifies  $\Si_{inv\; fact}(G_3)$ with the set of ``integer points''
$\De_L$ in the Weyl chamber $\De$.

\subsection{Linear algebra problems}
\label{linalg}

In this section we reformulate problems {\bf Q1-- Q3} stated in the introduction as {\em linear algebra} problems
in the case when the group $\ul{G}$ is $GL(m)$.

Let $\alpha$, $\beta$ and $\gamma$ be $m$-tuples of real numbers arranged
in decreasing order. In Problem {\bf P3} we let $\mathbb{K}$ be a
complete, nonarchimedean valued field. We assume that the valuation $v$ is discrete and takes values in $\Z$.
We let $\mathcal{O}\subset \K$ be the subring of elements with nonnegative valuation.

In order to state Problems {\bf P2} and {\bf P3} below we recall some definitions from
algebra. The {\em singular values} of a matrix $M$ are the (positive) square-roots of
the eigenvalues of the matrix $MM^*$. To define the {\em invariant factors}
of a matrix $M$ with entries in $\K$ note first that it is easy to
see that the double coset
$$
GL(m,\O)\cdot M\cdot GL(m,\O)\subset GL(m,\K)$$
is represented by a
diagonal matrix $D:= D(M)$. The  invariant factors of $M$ are the  integers
obtained by applying the valuation $v$ to the diagonal entries of $D$.
If we arrange the invariant factors in decreasing order they are uniquely determined by $M$.

Now we can state the three algebra problems that interest us here,
following the presentation in \cite{Fulton}.

\begin{itemize}
\item {\bf P1. Eigenvalues of a sum.}
Give necessary and sufficient conditions on $\alpha$,
$\beta$ and $\gamma$ in order that there exist Hermitian matrices $A$, $B$
and $C$ such that the set of eigenvalues (arranged in decreasing order) of
$A$, resp. $B$, resp. $C$  is  $\alpha$, resp. $\beta$, resp. $\gamma$ and
$$A+ B +C = 0.$$

\item {\bf P2. Singular values of a product.}
Give necessary and sufficient conditions on $\alpha$,
$\beta$ and $\gamma$ in order that there exist
matrices $A$, $B$ and $C$ in $GL(m,\C)$
the logarithms of whose singular values are
$\al,\beta$ and $\ga$, respectively,
so that
$$ABC = 1.$$

\item  {\bf P3. Invariant factors of a product.}
Give necessary and sufficient conditions on the integer vectors $\alpha$,
$\beta$ and $\gamma$ in order
that there exist matrices $A$, $B$ and $C$
in $GL(m,\mathbb{K})$
with invariant factors $\alpha, \beta$ and $\ga$, respectively,
so that
$$ABC = 1.$$
\end{itemize}

We now see that the problems {\bf Pi} are equivalent to the corresponding problems {\bf Qi};
we consider the field $\F=\C$ since for $\F=\R$ the discussion is similar.

\medskip
1. We have $\mathfrak{g} = \mathfrak{gl}(m)$ and the Cartan subspace
$\p$ is the space of Hermitian $m$-by-$m$ matrices, $\mathfrak{a}$ is the
space of diagonal $m$-by-$m$ matrices with real entries and the cone
$\Delta\subset {\mf a}$ is the cone in which the diagonal entries are arranged in
decreasing order. The maximal compact subgroup $K=U(m)$.
Thus the parameter space $\Sigma$ for the
double coset $K\backslash GL(m,\C)/K$
is the cone of $m$-tuples of real numbers arranged in decreasing order.
This agrees with the parameter space of the Problem {\bf P1}.

2. Observe that the Cartan subgroup
$A\subset GL(m,\C)$ consists of diagonal matrices with positive diagonal entries,
$A_\Delta$ consists of matrices $D\in A$ with diagonal entries arranged in the decreasing order,
and the projection $f$ of $g$ to $KgK\in \Si=A_\Delta\cong \Delta$ is given by
first sending $g$ to  $h=\sqrt{gg^*}$, then diagonalizing $h$, arranging the
diagonal values in the decreasing order and then taking their logarithms.
Thus $f(g)$ is given by the vector whose components are the logarithms of
the singular values of the matrix $g$.
We conclude that Problem {\bf P2} is equivalent to the Problem {\bf Q2}.

3. The torus $T\subset G=GL(m,\K)$ consists of diagonal
matrices with entries in $\K$; we find that  $\Psi(\phi)=\phi(\pi)$ is the diagonal matrix
with diagonal entries $(\pi^{x_1},\cdots,\pi^{x_m})$, where $\phi$ corresponds
to the lattice vector $(x_1,\cdots, x_m)$. Therefore {\bf P3} $\iff$ {\bf Q3}.

\medskip
{\bf Historic sketch.}

\medskip
To unify the notation, in what follows we will refer to the problem {\bf Q4} for the group $GL(m)$ as {\bf P4}.

The complete solution of Problems {\bf P1--P4} and the relation
between them were established only recently due to the efforts of
several people. After much classical work, A.~Klyachko in
\cite{Klyachko1} proved that Sol({\bf P1}) is a polyhedral cone in
$\mathbb{R}^{3m}$, he also computed a finite system of linear
homogeneous inequalities describing Sol({\bf P1}) in terms of the
Schubert calculus in the Grassmannians $G(p,\C^m)$. This
computation was subsequently improved by P.~Belkale
\cite{Belkale}.

A.~Klyachko then proved in \cite{Klyachko2} that Sol({\bf
P1})=Sol({\bf P2}), which was known as {\em Thompson's
conjecture}.

The sets of solutions for Problems {\bf P3} and {\bf P4} are also the same,
namely they are integral points in the above polyhedral cone. This result is
due to P.~Hall,  J.\ Green  and T.\ Klein, see \cite{Klein1} and \cite{Klein2};
see also \cite[pg. 100]{Macdonald} for the history of this problem.

In this paper we show (Theorem \ref{q4->q3})  that
the inclusion
$$
Sol(\mathbf{Q4}) \subset Sol(\mathbf{Q3})
$$
is true for {\em all} split reductive groups over $\Q$. This was the harder of the
two implications for $GL(m)$, proved by T.~Klein
in \cite{Klein1} and \cite{Klein2}, see also \cite[pg. 94--100]{Macdonald}.

The exceptional (i.e. not true for  all split reductive groups) inclusion
$$
Sol({\mathbf P3}) \subset Sol({\mathbf P4})
$$
was first proved for $GL(m)$ by P.~Hall but not published. In fact
it follows from a beautiful and elementary observation of
J.~Green, which is set forth and proved in \cite[pg.
91-92]{Macdonald}, and which we will outline in \S
\ref{saturationforgl}.

The description of the solutions to Problem {\bf P4} as the set of integral points
in the polyhedral cone Sol({\bf P1}) (known as the {\em Saturation Conjecture})
is due to A.~Knutson and T.~Tao, \cite{KnutsonTao} (see also \cite{DerksenWeyman} for a proof using quivers and \cite{KM}
for a proof in the spirit of the present paper using Littelmann's path model).
The Knutson-Tao theorem
combined with the above equivalence of Problems {\bf P3} and {\bf P4}
establishes that the set of solutions to Problem {\bf P3} is also the
set of integral points in the above polyhedral cone.

In the present paper
we reverse this path: We first prove directly that the set of solutions to
Problem {\bf P3} is the set of integral points in the above polyhedral cone,
then using the above equivalence of Problems {\bf P3} and {\bf P4} we deduce
that this set is also the set of solutions to Problem {\bf P4}, thus providing yet another proof
of the saturation conjecture .

We refer to \cite{AgnihotriWoodward}, \cite{AMW}, \cite{Belkale},
\cite{KumarBelkale},  \cite{BerensteinSjamaar} \cite{Fulton},
\cite{EL}, \cite{KM}, \cite{Klyachko1}, \cite{Klyachko2} for more
details and further developments.

\section{The existence of polygonal linkages and solutions to the algebra problems}
\label{settingupgeometry}

In this chapter we will show that the algebra problems {\bf R}($\ol{G}$) from chapter
\ref{sect:organizing algebra} can be restated geometrically, in terms of the
existence of triangles with the prescribed ``side-lengths''
in three classes of spaces of nonpositive curvature, which are:

\begin{enumerate}
\item $X$ is an {\em infinitesimal symmetric space} (a Cartan motion space), $X=\mathfrak{p}$.
\item $X$ is a symmetric space of nonpositive curvature.
\item $X$ is an Euclidean building.
\end{enumerate}

\subsection{Setting up the general geometry problem}
\label{generalnonsense}

We start with a purely set-theoretic discussion. Let $X$ be a set with a base-point $o$,
and let $\ol{G}$ is a group acting on $X$. Let $K$ denote the stabilizer of $o$ in $G$.
In general, the $\ol{G}$-action is not transitive (such examples
will appear when $X$ is a Euclidean building), we let $Y$ denote the orbit $G\cdot o$.
For pairs of points $(x,x')\in X^2$ we define the invariant
$$
\si_{\bar{G}}(x,x')
$$
 as the projection of $(x,x')$ to $\ol{G}\backslash X^2$. We will regard $\si_{\bar{G}}$ as the
``generalized $\ol{G}$-invariant distance'' between points in $X$.
Given $\tau\in \ol{G}\backslash X^2$ we regard the set
$$
S(o, \tau)= \{x\in X: \si_{\bar{G}}(o,x)=\tau\},
$$
as the ``sphere of radius $\tau$'' centered at $o$.

For pairs of
points $(x,y)\in Y^2$ one can interpret the invariant $\si_{\bar{G}}(x,y)$ as follows:
Let $\Si$ denote the quotient space
$$
\Si=K\backslash Y= K\backslash \ol{G}/K.
$$
Then for pairs of points $(x,y)=(gK, hK)$ in $Y$ we can identify  $\si_{\bar{G}}(x,y)$
with the double coset $Kg^{-1}hK$. In other words, translate the pair
$(x,y)$ by $g^{-1}$ to $(o,z)=(o, g^{-1}(o))$ and then
project $g^{-1}(o)$ to the element $\si_{\bar{G}}(x,y)\in K\backslash Y$.
Then  two pairs $(x,y)$ and $(x',y')$ in $Y^2$
belong to the same $\ol{G}$-orbit if and only if $\si_{\bar{G}}(x,y)=\si_{\bar{G}}(x',y')$;
Hence, $\ol{G}\backslash Y\times Y$ can be bijectively identified with $\Si=K\backslash \ol{G}/K$.

We will use the notation $\Upsilon$ for the map $G\to \Si$.
For $\si\in \Sigma$ we let $\O_\si \subset \ol{G}$ denote $\Upsilon^{-1}(\si)$.

It is now clear that the Problem {\bf R}($\ol{G}$) stated in chapter
\ref{sect:organizing algebra} can be restated as a special case ($n=3$) of the following:

\begin{prob}
\begin{itemize}
\label{linalgprob}
Give conditions on the vector of generalized side-lengths $\oa{\sigma}=(\sigma_1,\cdots,\sigma_n)\in \Si^n$
that are necessary and sufficient in order that there exist
elements $g_1, g_2, \cdots, g_n$ in $\overline{G}$ such that
$$
\Upsilon(g_i)=\sigma_i , 1\leq i \leq n,\quad \mbox{and}
\quad \prod_{i=1}^n g_i=id.$$
\end{itemize}
\end{prob}

\no We set $\mathcal{O}=\prod_{i=1}^n \mathcal{O}_{\sigma_i}$.
Note that $(K\times K)^n$ acts on $\mathcal{O}$ by right and left multiplications.

Motivated by the connection between the moduli spaces of $n$-gons and symplectic
quotients when $Y$ is a complete simply-connected 3-dimensional Riemannian
manifold of constant curvature (see \cite{KapovichMillson96}, \cite{KapovichMillsonTreloar}),
\cite{Treloar}) we define the ``momentum map"
$\mu:\mathcal{O}\lra \overline{G}$ by
$$
\mu(g_1,...,g_n)= g_1 \cdot g_2 \cdots g_n.
$$
We define the analogue of the symplectic quotient of $\mathcal{O}$ by $K$, by
$$
\mathcal{O}/\!/K=\{g\in \mathcal{O} : \mu(g)=id \}/K.
$$
Here we divide out by the diagonal action of $K$, where $K$ acts on each
factor by {\em conjugation}. This action is the only action of $K$ on
$\mathcal{O}$ that we will use henceforth.

We let $\mathcal{Q}_{\si}(\ol{G})$ be the subset of $g\in \mathcal{O}$
such that $\mu(g)=id$. Thus $\mathcal{Q}_{\si}(\ol{G})$ is the set of
solutions of the problem \ref{linalgprob} and $\mathcal{O}/\!/K$ is its
quotient by $K$.

Our goal now is to reformulate the problem \ref{linalgprob}
in more geometric terms.

\medskip
An  {\em $n$-gon} (or {\em a closed $n$-gonal linkage})
in $X$ with the vertices $x_1,...,x_n\in X$
is an $n$-tuple $\oa{x}=(x_1,...,x_n)\in X^n$, regarded as a map $\Z/n\Z\to X$.
The {\em generalized side-lengths} of the polygon $\oa{x}$ are the elements
$$
\si_i= \si_{\bar{G}}(x_i, x_{i+1})\in \ol{G}\backslash X^2, \quad i\in \Z/n\Z.
$$
Thus for each $n$-gon we get a vector $\oa{\si}=(\si_1,...,\si_n)$ of its
generalized side-lengths. The group $\ol{G}$ acts naturally
on the space ${\mathcal P}ol_n(X)$ of all $n$-gons, preserving the generalized side-lengths of the polygons.

Fix $\oa{\si}=(\si_1,...,\si_n)$ and form the space
$$
{\mathcal P}ol_{n,\si}(X) :=\{ P\in {\mathcal P}ol_n(X): \si_{\bar{G}}(x_i, x_{i+1})=\si_i\},
$$
and its quotient, the {\em moduli space} of polygons with the fixed side-lengths in $X$:
\begin{equation}
\label{moduli}
\M_{n,\si}(X) :={\mathcal P}ol_{n,\si}(X)/\ol{G}.
\end{equation}
In the case when $X$ is a topological space and $\ol{G}$ acts homeomorphically on $X$ we give
the space $\M_{n,\si}(X)$ the quotient topology.

We note that  ${\mathcal P}ol_n(Y)$ sits naturally in  ${\mathcal P}ol_n(X)$ as the subset of
polygons with vertices in $Y$. The generalized side-lengths of the polygons
with vertices in $Y$ are regarded as elements of $\Si=K\backslash Y$.

We define a map
$\Phi:\mathcal{Q}_{\sigma}(\overline{G})\to  Y^{n}$
by
$$
\Phi(\oa{g})=(o,g_1\cdot o, g_1g_2 \cdot o, \cdots, g_1g_2\cdots g_{n-1} \cdot o).$$
We see that $\Phi$ is $K$-equivariant (recall that $K$ acts on $\mathcal{O}$
by conjugation on each factor).

\begin{lem}
\label{trick}
The map $\Phi$ induces a surjection $\ol{\Phi}$ from $\mathcal{O}/\!/K$ onto the moduli space
$\mathcal{M}_{n,\sigma} (Y)$.
\end{lem}
\proof We first verify that each polygon $\Phi(\oa{g})$ has the ``correct'' side-lengths.
By the left-invariance of $\si_{\bar{G}}$ we have
$$
\si_{\bar{G}}(\prod_{i=1}^k g_i\cdot o,\prod_{i=1}^{k+1} g_i\cdot o) =
\si_{\bar{G}}(o,g_k\cdot o)=\si_k.
$$
To prove surjectivity of the map $\ol{\Phi}$
let $\oa{y}=(y_1,...,y_n)\in \mathcal{P}ol_{n,\oa{\sigma}}(X)$. Choose $g$ such that
$g\cdot y_1=o$ and replace $\oa{y}$ by $g^{-1}\cdot \oa{y}$.
Since $\sigma_{\bar{G}}(y_1,y_2)=\sigma_1$ and (now) $y_1=o$ we find that $y_2\in S(o,\sigma_1)$.
Hence there exists $g_1\in \mathcal{O}_{\sigma_1}$ such that $y_2=g_1\cdot o$. Similarly, there exists
$g_2 \in \mathcal{O}_{\sigma_2}$ such that $g_1g_2\cdot o =y_3$.
Continuing in this way we get an $n$-tuple $(g_1,\ldots , g_{n-1}, g_n')\in \O$ such that
$$
(o,g_1\cdot o, g_1g_2 \cdot o, \cdots, g_1g_2\cdots g_{n-1} \cdot o, g_1g_2\cdots g_{n-1}g_{n}' \cdot o)
= (y_1,...,y_n, y_1),$$
where $y_1=o$. Thus $g_1g_2\cdots g_{n-1}g_{n}'=k\in K$ and
after replacing $g_n'$ with $g_n=g_n'k^{-1}$
we get an $n$-tuple $\oa{g}=(g_1,...,g_n)\in \O$  such that $\mu(\oa{g})=id$ and
$\Phi(\oa{g})=\oa{y}$. \qed

Thus we find that the problem \ref{linalgprob}
can be solved if and only the moduli space $\mathcal{M}_{n,\sigma}(Y)$
is nonempty.

To prove a version of Lemma \ref{trick} with a {\em bijective} map and to
get a better analogy with the symplectic quotients we have to make a further assumption about
the action $\ol{G}\acts Y$:

\begin{ass}
Through the rest of this section we assume that $\ol{G}$ contains
 a subgroup $B$ which acts simply-transitively on $Y$.
\end{ass}

Now, instead of the orbits $\O_{\si_i}$ and their product $\O$ we consider the intersections:
$$
\widehat\O_{\si_i} := \O_{\si_i} \cap B, \hbox{~~and~~} \widehat{O} := \O\cap B^n.
$$
The group $K$ no longer acts on $\widehat{O}$ by conjugations, instead one has the use a
{\em dressing action}, $dress(K)$, see \cite{KapovichMillsonTreloar} for the definition
in the case when $X=Y=\H^3$, $G=PSL(2,\C)$ and $B$ fixes a point at infinity.
In the case of infinitesimal symmetric spaces $X=\p$ (where $\p\oplus \k=\g$)
one takes $B=\p$ and the adjoint action of $K$ as the dressing action.
We now redefine
the symplectic quotient as follows:
$$
\widehat{O}//K:= (\widehat{O} \cap \mu^{-1}(id))/dress(K).
$$
Then one gets an analogue of Lemma \ref{trick} (which we do not need for the purposes of this paper),
we refer the reader to \cite{KapovichMillsonTreloar} for the discussion in the case of $X=\H^3$:

\begin{lem}
\label{trickk}
The map $\Phi$ induces a bijection $\widehat{\Phi}$ from $\widehat{\O}//K$ onto the moduli space
$\mathcal{M}_{n,\sigma}(Y)$. Moreover,
let $Y$ be a topological space and $G\subset Homeo(Y)$.
Then the quotient spaces $\widehat{\O}//K$ and ${\mathcal M}_{n,\si}(Y)$ have natural topology
and the bijection $\widehat{\Phi}$ is a homeomorphism.
\end{lem}

The above discussion has been completely formal, our next goal is to describe the spaces $X$
which can be used to analyze the problem \ref{linalgprob} for various groups $\ol{G}$.

\subsection{Geometries modeled on Coxeter complexes}
\label{geometries}

Fix a spherical or Euclidean Coxeter complex $(A,W)$, where $A$ is a Euclidean space $E$
or a sphere $S$ and $W=W_{aff}$ or $W=W_{sph}$ is a (possibly nondiscrete) Euclidean or a spherical
Coxeter group acting on $A$. If $W$ is discrete then the Coxeter complex $(A,W)$ is called
{\em discrete}. In the case of Euclidean Coxeter complexes we let $L_{trans}\subset W$
denote the translation subgroup of $W$.
Pick a special vertex $o\in E$ with stabilizer $W_{sph}$ (see
chapter \ref{datum}) and let $\De\subset E$ be a Weyl chamber of
$W_{sph}$. We will use the notation $\De_{sph}$ to denote the
ideal boundary of $\Delta$; $\De_{sph}$ is contained in the sphere
$S$ which will be regarded as the sphere at infinity of $E$; thus
$\De_{sph}$ is a fundamental domain for the action $W_{sph}\acts
S$.

Let $Z$ be a metric space. A {\em geometric structure} on $Z$ {\em modeled on} $(A,W)$
consists of an atlas of isometric embeddings $\varphi:A\embed Z$
satisfying the following compatibility condition:
For any two charts $\varphi_1$ and $\varphi_2$,
the transition map $\varphi_2^{-1}\circ\varphi_1$ is the restriction
of an isometry in $W$.
The charts and their images, $\varphi(A)=a\subset Z$, are called {\em apartments}.
We will sometimes refer to $A$ as the {\em model apartment}.
We will require that there are {\em plenty of apartments}
in the sense that any two points in $Z$ lie in a common apartment.
All $W$-invariant notions introduced for the Coxeter complex $(A, W)$,
such as walls, singular subspaces, chambers etc.,
carry over to geometries modeled on $(A, W)$.

One defines the group of automorphisms $Aut(A)$ of the model apartment
as the group of isometries of $A$ which normalize the subgroup
$W$. If $X$ is a space modeled on $(A, W)$ then an isometry $g: X\to X$
is an automorphism if it sends apartments to apartments and
for each pair of apartments $(A,\varphi_1), (A,\varphi_2)$
the composition $\varphi_2^{-1} \circ g\circ \varphi_1$
is the restriction of an automorphism of $A$.
The group of
automorphisms of $X$ is denoted $Aut(X)$.

Examples of the above geometries are provided by
symmetric spaces of noncompact type
and their infinitesimal analogues (infinitesimal symmetric spaces).
These are modeled on Euclidean Coxeter complexes with transitive
affine Weyl group.
In the case of a symmetric space $X$,
the apartments are the maximal flats.
The associated Coxeter complex has the form $(E,W_{aff})$
where $E$ is an apartment and $W_{aff}$ is the group generated by reflections
at singular hyperplanes.

Take a real or complex reductive Lie group $G$, the Lie algebra
$\g$ of $G$ has the decomposition $\g=\p \oplus \k$, where
$K\subset G$ is a maximal compact subgroup, let $\k$ denote the
Lie algebra of $K$. We will identify the Cartan subspace $\p$ with
the tangent space $T_pX$ to the symmetric space $X=G/K$ at the
point $p$ stabilized by $K$. The  subspaces $\p\subset \g$,
equipped with the affine action of $\ol{G}=\p \rtimes K$, are {\em
infinitesimal symmetric spaces} in the following sense. The
$(E,W_{aff})$-structure on $X$ induces a $(E,W_{aff})$-structure
on the {\em Cartan subspace} $\p=T_pX$, such that the apartments
in $\p$ are the translates of the Cartan subalgebras (i.e.\ the
maximal abelian subalgebras) in $\p$. The apartments through $0$
in $\p$ are the tangent spaces to the apartments through $p$ in
$X$. The term ``infinitesimal symmetric space" may be also
justified by noting that there is a one-parameter family of spaces
$X_{\epsilon}$ parameterized by $\epsilon \geq 0$ and all
isometric to $X$ for $\epsilon > 0$ such that $X_0 =
\mathfrak{p}$, see \cite[\S 5]{KapovichMillsonTreloar}.

The last kinds of geometries considered in this paper are
spherical and Euclidean buildings. We refer the reader to \cite{Ballmann} for
the definitions of $CAT(\kappa)$ metric spaces.

\begin{dfn}
A {\em spherical building} is a CAT(1)-space modeled
on a spherical Coxeter complex.
\end{dfn}

We will use the notation $\tangle$ for the metric in a spherical building.
Spherical buildings have a natural structure as
polysimplicial piecewise spherical complexes.
We prefer the geometric to the combinatorial view point
because it appears to be more flexible.

\begin{defn}
A discrete Euclidean building is a CAT(0)-space modeled on a discrete Euclidean Coxeter complex.
\end{defn}

In the non-discrete case the definition of Euclidean buildings is more subtle,
see \cite[section 4.1.2]{KleinerLeeb}.
We refer to \cite{KleinerLeeb} for a thorough discussion of buildings from the geometric viewpoint.

A building is called {\em thick}
if every wall is an intersection of apartments.
A non-thick building can always be equipped
with a natural structure of a thick building by reducing the Weyl group.

\begin{ex}
If $X$ is a symmetric space of noncompact type
or a thick Euclidean building modeled on the Coxeter complex $(E,W_{aff})$,
then its ideal boundary $\tits X$ is a thick spherical building
modeled on $(\tits E,W_{sph})$.
In the case that $X$ is a building,
the {\em spaces of directions} $\Si_xX$ are spherical buildings
modeled on $(\tits E,W_{sph})$.
The building $\Si_xX$ is thick if and only if $x$ is a special vertex of $X$.
We note that in the case when $X$ is a discrete building, $\Si_xX$ is just the
link of the point $x\in X$.
\end{ex}

Let $B$ be a spherical building modeled on a spherical Coxeter complex $(S, W_{sph})$.
We say that two points $x, y\in B$ are antipodal, if $\tangle(x,y)=\pi$; equivalently,
they are antipodal points in an apartment $s\subset B$ containing both $x$ and $y$.
The quotient map $S\ra S/W_{sph}\cong\De_{sph}$
induces a canonical projection
$\theta:B\ra \De_{sph}$
folding the building onto its model Weyl chamber.
The $\theta$-image of a point in $B$ is called its {\em type}.

\begin{rem}
To define $\theta(x)$ pick an apartment $s$ containing $x$ and a chart $\phi: S\to s$.
Then $\theta(x)$ is the projection of $\phi^{-1}(x)$ to  $S/W_{sph}\cong\De_{sph}$. We note that
this is clearly independent of $s$ and $\phi$.
\end{rem}

The same definition applies in the case of Euclidean buildings $B$. The difference however
is that the action $W_{aff}\acts E$ in general is no longer discrete, so we cannot
identify the image of $B\to E/W_{aff}$ with a simplex. If $W_{aff}$ acts as a lattice on $E$, then
$E/W_{aff}$ can be identified with a {\em fundamental alcove} for the action $W_{aff}\acts E$.

We now give two properties of the projection $\theta$:

1. If $h: a\to a'$ is an isomorphism of apartments in $B$
(i.e. ${\phi'}^{-1}\circ h\circ \phi\in W$)
then $\theta\circ h=\theta$.

2. If $B$ is a spherical building, $x, x'\in B$ which belong
to apartments $a, a'$ respectively
and $-x\in a, -x'\in a'$ are antipodal to $x, x'$,
then $\theta(x)=\theta(x')$ implies
$\theta(-x)=\theta(-x')$.
To prove this pick an isomorphism $h: a\to a'$.
Then (since $\theta(x)=\theta(x')$)
there exists $w\in W\acts a'$ such that  $w(h(x))=x'$.
Hence $w\circ h(-x)=-w\circ(x)=-x'$.
The claim now follows from 1.

\medskip
{\bf DeRham decomposition of Euclidean buildings.}
Suppose that $X$ is a thick Euclidean building modeled on a {\em reducible} discrete
Coxeter complex $(E,W_{aff})$. Consider the deRham decomposition
of the building $X$:
$$
X= X_0 \times X_1 \times ...\times X_s,
$$
where each $X_i$ is a thick Euclidean building modeled on the Euclidean Coxeter complex
$(E_i, W^i_{aff})$ and $X_0\cong E_0$ is the flat deRham factor.
For each $i=1,...,s$ the group $W^i_{aff}$ acts as an irreducible
affine Coxeter group on the Euclidean space $E_i$, for $i=0$ we get the trivial
Coxeter group $W^0_{aff}$. Accordingly,
the Euclidean space $E$ splits (metrically) as the product
$$
E_0 \times \prod_{i=1}^s E_i,
$$
this decomposition is invariant under $W_{aff}$ which in turn splits as
$$
W_{aff}= \prod_{i=0}^s W_{aff}^i,
$$
as explained in Section \ref{datum}.

\subsection{Bruhat-Tits buildings associated with nonarchimedean reductive
Lie groups}
\label{BTB}

In this section we describe properties of the Euclidean building
(the {\em Bruhat-Tits' building}) associated to a reductive
nonarchimedean Lie group by Bruhat and Tits. The reader familiar
with \cite{BT} can skip this section.

Let  $\K$ be a  valued field  with valuation $v$ and value group
$\Z$. Let $\O$ denote the subring in $\K$ which consists of
elements with nonnegative valuation. Let $\ul{G}$ be a connected
reductive algebraic group over $\K$ which has relative rank $l$.
The subgroup $K:=\ul{G}(\O)\subset G=\ul{G}(\K)$ is a maximal
bounded subgroup of $G$. In this section we review the properties
of a Euclidean building $X=X_{\ul{G}}$ attached to the group
$\ul{G}$. We refer the reader to \cite{Tits} and \cite[Chapter
7]{BT} for more details.

We begin by recalling the notation from section \ref{tough}.
Let $\ul{T}\subset \ul{G}$ be a maximal $\K$-split torus, $\ul{Z}\subset \ul{G}$
and $\ul{N}\subset \ul{G}$ be its centralizer and normalizer (over $\K$) respectively.
We also get the groups $N, Z$ and $T$ of $\K$-points of the corresponding algebraic
groups. Then $V:= X_*(\ul{T})\otimes \R$, where $X_*(\ul{T})$ is the group of
cocharacters of $\ul{T}$.
The Euclidean space $E$ is the affine space underlying $V$ with appropriately chosen
Euclidean metric.
Let $R_{rel}\subset V^*$ denote the relative root system of the pair $(\ul{G}, \ul{T})$; then
$W_{sph}$ is the finite Coxeter group corresponding to this system.

Definition \ref{nu} gives us a homomorphism $\nu: Z\to V$,
with the kernel $Z_c\subset Z$ and the image equal to the extended cocharacter lattice $L_{\ul{G}}$.
We also get the quotient group $N/Z_c$, which acts on the Euclidean space $E$
discretely and isometrically through a group
$\tilde{W}$. Unless $\ul{G}$ is semisimple, $\tilde{W}$
{\em is not} an affine Weyl group.

In the case when the group $\ul{G}$ is simply-connected and semisimple, one can take the
pair $(E, \tilde{W})$ as the Euclidean Coxeter complex $(E, W_{aff})$
of the Bruhat-Tits building $X_{\ul{G}}$ attached to the group $\ul{G}$.
In general however it is not the case and one has to do more work to define
$W_{aff}$. The construction of this group will be unimportant for us (we refer to \cite{Tits}
for the explicit construction), the important properties of $W_{aff}$
are the following:

\begin{enumerate}
\item $\dim(E)=l$, the relative rank of $\ul{G}$ or the dimension of $\ul{T}$.

\item $W_{aff}\subset \tilde{W}$ is a normal subgroup \cite[section 1.7, page 34]{Tits},
the index\newline
 $|\tilde{W}: W_{aff}|$ is finite in the semisimple case.

\item $\tilde{W}= W_{sph}\ltimes L_{\ul{G}}$.

\item $W_{aff}$ is an affine Coxeter group attached to a root system $R\subset V^*$.

\item The root system $R$ in general {\em is not the same}
 as $R_{rel}$ but they have the same rank
and the same finite Weyl group $W_{sph}$, \cite[section 1.7]{Tits}.

\item $W_{aff}=W_{sph}\ltimes L_{trans}$, where $W_{sph}$ is the finite Weyl group as above.

\item We have the inclusions
$$
Q(R^{\vee})=L_{trans}\subset L_{\ul{G}} \subset N_{aff}= P(R^{\vee}).
$$

\item If $\ul{G}$ is a simply-connected semisimple group then $W_{aff}=\tilde{W}$, \cite[section 1.13]{Tits}.
If $\ul{G}$ is split over an unramified extension of $\K$ then
$L_{\ul{G}}=X_*(\ul{T})$, \cite[section 1.3]{Tits}.

\item If $\ul{G}$ is split then the reduced (Bruhat-Tits) root system associated with $W_{aff}$
is the usual reduced root system $R_{rel}$ of the group $\ul{G}$
(note that $R_{rel}$ is also the absolute root system since we are in the split case).
\end{enumerate}

\no Here as usual, $L_{trans}$ is the translation subgroup of $W_{aff}$ and $N_{aff}$ is
the normalizer of $W_{aff}$ in $V$.
When we are dealing with {\em the root system associated with}
$\ul{G}$, as in 4 or 5 above, we will refer to $R$, and not to $R_{rel}$.

\begin{defn}
We will call $R$ the {\em Bruhat-Tits root system} associated with the algebraic group $\ul{G}$.
\end{defn}

We refer the reader to Tableau des \'echelonages in \cite[p.
29-30]{BT} for computation of the Bruhat-Tits root systems.

We note that the root systems $R$ and $R_{rel}$ are not very
different since they have the same finite Coxeter group. Thus, in
the irreducible case, both root systems are either isomorphic or
one of them is of type $B_l$ and the other of type $C_l$. In
particular, the {\em saturation factors} for $R$ and $R_{rel}$ are
exactly the same (see chapter \ref{sec:satfactor}).

Having described the properties of the Coxeter complex $(E, W_{aff})$, we will
describe the properties of the {\em Bruhat-Tits building} $X= X_{\ul{G}}$ associated
with the group $\ul{G}$:

\begin{lis}
\label{list}
\begin{enumerate}

\item $X$ is modeled on the Euclidean Coxeter complex $(E,W_{aff})$ described above.

\item The group $G$ acts on $X$ by (isometric) automorphisms.

\item The subgroup $K$ is the stabilizer of a special vertex $o$ in $X$.

\item  $G$ acts transitively on the set of apartments in $X$.

\item  For each apartment $a\subset X$ let $G_a$ be the stabilizer of
 $a$ in $G$. Then
the image of $G_a$ in $Aut(a)$ is the group $\tilde{W}$ (containing $W_{aff}$
as a normal subgroup).

\item  The building $X$ is thick (see \cite[Prop. 7.4.5]{BT}).

\item The group $W_{aff}$ acts discretely and (in the semisimple case) cocompactly on $E$.
\end{enumerate}
\end{lis}

\begin{rem}
The properties in the List \ref{list}, except the last one, hold in the case
of symmetric spaces and their infinitesimal analogues.
\end{rem}

We note that the affine space $E$ has the distinguished point $o$,
the origin (corresponding to the trivial cocharacter).
The corresponding vertex in $X$ is also denoted by $o$, it
is stabilized by $K=\ul{G}(\O)$.

\subsection{Geodesic polygons}
\label{Geodesic polygons}

In this paper we will be considering polygons in the
metric spaces modeled on the Coxeter complexes, which were discussed in section \ref{geometries}.
With the exception of the rank zero spherical buildings,
all such metric spaces $X$ are {\em geodesic} and thus we define a {\em geodesic polygon}
in $X$  as a polygon with the vertices $z_1\dots z_n$
together with a choice of {\em sides} $\ol{z_iz_{i+1}}$,
i.e. the geodesic segments connecting $z_i$ to $z_{i+1}$. We note that in the case of
metric spaces modeled on Euclidean Coxeter complexes the sides are uniquely determined by the
polygon  $z_1\ldots z_n$. We recall that {\em $n$-gons} in $X$ are regarded as maps $\Z/n\Z\ra X$.

Let $(A,W)$ be a spherical or Euclidean Coxeter complex.
The complete invariant of a pair of points $(x,y)\in A^2$ with respect to the action
$W\acts A$, is its image $\si_{W}(x,y)$ under the canonical projection to
$A\times A/W$. We define the {\em refined length} of a geodesic segment $\ol{xy}$
as $\si_{ref}(x,y):=\si_{W}(x,y)$. This notion carries over to geometries
modeled on the Coxeter complex $(A,W)$:
For a pair of points $(x,y)$
pick an apartment $a$ containing $x,y$
and, after identifying $a$ with the model apartment $A$,
let $\si_{ref}(x,y)$ be the projection to $A\times A/W$.

In the case of Euclidean Coxeter complexes there are extra structures associated
with the concept of refined length.
Given a Euclidean Coxeter complex $(E, W_{aff})$,
pick a special vertex $o\in E$.
Then we can regard $E$ as a vector space $V$, with the origin $0=o$. Let $\De\subset E$ denote a Weyl chamber
 of $W_{sph}$, the tip of $\De$ is at $o$.

Suppose that $L$ is a subgroup of the group $V$ of all translations of $E$ so that:
\begin{equation}
\label{lattice}
 L_{trans}\subset L \subset N_{aff},
\end{equation}
where $N_{aff}\subset V$ is the normalizer of $W_{aff}$. Since $o$ is the origin in $E$,
we will identify the orbits $L\cdot o$ and $L_{trans}\cdot o$ with $L$ and $L_{trans}$ respectively.
We define the set of {\em refined $L$-integral lengths}
as the subset
$$
(L\times L)/W_{aff}\subset (E\times E)/W_{aff}.$$
If $L=L_{trans}$ then we have a natural bijection
$$
(L_{trans}\times L_{trans})/W_{aff}\cong \De \cap L_{trans}.
$$
Since $(E, W_{aff})$ is a Euclidean Coxeter complex, there is also a coarser notion of {\em $\De$-length}
obtained from composing $\si_{ref}$ with the natural forgetful map
\[ E\times E/W_{aff}\ra E/W_{sph}\cong\De .\]
To compute the $\De$-length $\si(x,y)$
we regard the oriented geodesic segment $\ol{xy}$
as a vector in $E$ and project it to $\De$.

Again, the concepts of $\De$-length, $L$-integral $\De$-lengths, etc., carry over to the geometries
modeled on $(E,W_{aff})$.
Note that $\De$-length and refined length coincide for symmetric spaces
and their infinitesimal analogues
because the affine Weyl group acts transitively.
We define the set of $L$-{\em integral} $\De$-lengths
as the subset $\De\cap L\subset\De$.
A segment with $L$-integral $\De$-length has $L$-integral refined length
iff its vertices lie in the distinguished orbit $L\cdot o$ (identified with $L$).
For segments with endpoints of type $W_{aff}\cdot o$,
the notions of $\De$-length and refined length are equivalent.

Given a collection $\tau=(\tau_1,...,\tau_n)\in \De^n$ of $\De$-lengths we define
the {\em moduli space $\M_{n,\tau}(X)$}
as the quotient of the collection of geodesic polygons in $X$
with the $\De$-side-lengths $\tau$ by the action of the group $G$.
We give $\M_{n,\tau}(X)$ the quotient topology.

We are now ready to state the questions which will be (for $n=3$)
the geometric counterparts to the algebra questions in chapter \ref{sect:organizing algebra}:

\begin{prob}
\label{mainques0}
Let $X$ be a symmetric space of noncompact type or an infinitesimal symmetric space,
or a thick Euclidean building.
Describe the set $D_n(X)\subset \De^n$ of $\De$-side
lengths which occur for geodesic $n$-gons in $X$.
\end{prob}

\begin{rem}
The above problem was solved in the papers
\cite{KapovichLeebMillson1, KapovichLeebMillson2}.
In these papers  the notation ${\mathcal
P}_n$ was used in place of $D_n$. We refer the reader to section
\ref{polyhedron} for the description of the solution of Problem
\ref{mainques0}.
\end{rem}

\begin{prob}
\label{geometricquestion}
Let $X$ be a thick Euclidean building.
Describe the set
$$
D_n^{ref,L}(X)\subset (W\backslash L\times L)^n\subset (W\backslash A\times A)^n$$
of refined $L$-integral side-lengths which occur for $n$-gons in $X$.
\end{prob}

One of the key results concerning the above problem is the
following

\begin{thm}
[Transfer Theorem, see \cite{KapovichLeebMillson2}]
\label{transfer}  Suppose that $X, X'$ are thick Euclidean
buildings modeled on the Coxeter complexes $(A,W), (A',W')$
respectively, $L, L'$ are lattices as in (\ref{lattice}), $\iota:
(A,W)\to (A',W')$ is an embedding of Coxeter complexes such that
$\iota(L)\subset L'$. Then $\iota$ induces an embedding
$$
\iota_*: D_n^{ref,L}(X) \to D_n^{ref,L'}(X')
$$
\end{thm}

Below we explain how given an algebraic problem
{\bf Q1}--{\bf Q3} one finds a metric space $X$
modeled on a Euclidean Coxeter complex,
so that the geometric problem \ref{geometricquestion} (for $n=3$)
is equivalent to the corresponding algebraic problem.

\medskip
As we explained in the beginning of chapter \ref{sect:organizing algebra}, with each
problem {\bf Qi}, {\bf i=1, 2, 3}, we can associate a pair of groups $K\subset \ol{G}$.

1. For the Problem {\bf Q1} we take
$\overline{G}= K \ltimes \mathfrak{p}$
 and let $X:= \p$.
Then $X$ is a metric space modeled on the Euclidean Coxeter complex
$(A,W)=(E, W_{aff})$.

2. For the Problem {\bf Q2} we take $\overline{G}=G$
and let $X$ be the symmetric space $X=G/K$.

In both cases $\ol{G}$ acts transitively on $X$ and we apply
Lemma \ref{trick} to the transitive action $\ol{G}\acts X$,
to conclude that {\bf Q1} and {\bf Q2} are equivalent to
Problem \ref{mainques0} for the corresponding space $X$.

\begin{ex}
Let $G=GL_m(\C)$, $K=U(m)$. Then the symmetric space $X$ associated with $G$ is $G/K$, which
is the space ${\mathbb P}_m$ of positive-definite Hermitian $m\times m$ matrices. The
model apartment $A$ in ${\mathbb P}_m$ consists of diagonal matrices with positive diagonal entries.
The finite Weyl $W_{sph}={\mathcal S}_m$ group acts on $A$ by permutations  of the diagonal elements;
the affine Weyl group $W_{aff}$ is isomorphic to $W_{sph}\ltimes \R^{m-1}$.
Thus we get a geometric model for analyzing the Problem {\bf P2} from the Introduction.
\end{ex}

3. For the Problem {\bf Q3} we assign to the group $G=\ul{G}(\K)$
the Euclidean (Bruhat-Tits) building $X$ as it was
explained in the previous section.
We let $Y$ be the $G$-orbit of the special vertex $o\in X$.

Although, unlike in the previous two examples,
the group $G$ does not act transitively on the building $X$, this
group acts transitively on the subset $Y\subset X$
and we apply Lemma \ref{trick} to the action $G\acts Y$
to see that the problem {\bf Q3}
(or, equivalently, {\bf R}($G$)) is equivalent to finding necessary and
sufficient conditions for existence of geodesic triangles in $X$ whose
vertices are in $Y$ and whose $\Delta$-lengths are the prescribed elements of $\Delta_L$,
where $L=L_{\ul{G}}$. More generally, there is a surjective map
$$
pr: W\back A\times A\to G\back X\times X,
$$
and by applying the Transfer Theorem \ref{transfer} we get:

\begin{prop}
Suppose that there exists a polygon $P$ in $X$ whose
$\si_{G}$-side-lengths
are $(\si_1,...,\si_n)$. Then for any choice of
$\tilde\si_i\in pr^{-1}(\si_i)$,
there exists a polygon $\tilde{P}$ in $X$
whose $\si_{ref}$-side-lengths are
$\tilde\si_1,...,\tilde\si_n$.
\end{prop}

Therefore we will stick to the notion of the refined
side-length $\si_{ref}$ through the rest of the paper.

\section{Weighted configurations, stability and the relation to polygons}
\label{gau}

Let $X$ be a symmetric space of nonpositive curvature or a Euclidean building.
Recall that
the ideal boundary $B=\tits X$ has the structure of a spherical building, the metric
on $B$ is denoted by $\tangle$.
Given a Weyl chamber $\De$ in $X$, we get a spherical Weyl chamber
$\De_{sph}=\geo \De\subset \tits X$.
We will identify $\De_{sph}$ with the unit vectors in $\De$.
Recall that there is a canonical projection $\theta:\tits X\ra\De_{sph}$,
see section \ref{geometries}.

Take a collection of weights $m_1,...,m_n\geq0$ and define a finite measure space
$(\Z/n\Z, \nu)$ where the measure $\nu$ on $\Z/n\Z$ is given by $\nu(i)=m_i$.
An $n$-tuple of ideal points $(\xi_1,...,\xi_n)\in B^n$ together with
$(\Z/n\Z, \nu)$ determine a {\em weighted configuration at infinity}, which is a map
$$
\psi: (\Z/n\Z, \nu)\to \tits X.
$$
The {\em type} $\tau(\psi)=(\tau_1,\dots,\tau_n)\in\De^n$
of the weighted configuration $\psi$
is given by $\tau_i=m_i\cdot\theta(\xi_i)$.
Let $\mu=\psi_*(\nu)$ be the pushed forward measure on $B$.
We define the {\em slope} of a measure $\mu$ on $B$ with finite total mass $|\mu|$ as
\[ slope_{\mu}(\eta) = -\int_B \cos\tangle(\xi,\eta) \;d\mu(\xi) .\]
In this paper we will consider only measures with finite support.

\begin{dfn}
[Stability]
A measure $\mu$ on $B$ (with finite support)
is called {\em semi\-stable}
if $slope_{\mu}(\eta)\geq0$
and {\em stable}
if $slope_{\mu}(\eta)>0$
for all $\eta\in B$.
\end{dfn}

There is a refinement of the notion of semistability
motivated by the corresponding concept in geometric invariant theory.

\begin{dfn}[Nice semistability]
A measure $\mu$ on $B$ (with finite support) is called {\em nice semistable}
if $\mu$ is semistable and
$\{slope_{\mu}=0\}$ is a subbuilding or empty.
In particular, stable measures are nice semistable.
\end{dfn}

A weighted configuration
$\psi$ on $B$
is called {\em stable, semistable} or {\em nice semistable}, respectively,
if the corresponding measure $\psi_{\ast}\nu$ has this property.


For the purposes of this paper,
i.e.\ the study of polygons,
nice semistability plays a role in the case of symmetric spaces
and  infinitesimal symmetric spaces only. We note however that for these spaces, existence
of a semistable configuration $\psi$ on $\tits X$ implies existence of a nice semistable
configuration on $\tits X$, which has the same type as $\psi$, see \cite{KapovichLeebMillson1}.

\begin{ex}
(i)
Let $B$ be a spherical building of rank $0$.
Then a measure $\mu$ on $B$ is stable
iff it contains no atoms of mass $\geq\half|\mu|$,
semistable iff it contains no atoms of mass $>\half|\mu|$,
and nice semistable iff it is either stable or consists of two atoms of equal mass.

(ii)
Suppose that $B$ is a unit sphere
and regard it as the ideal boundary of a Euclidean space $E$,
$B=\tits E$. Each semistable measure $\mu$ has slope zero everywhere.
\end{ex}

Define the subset $\De^n_{ss}(B)\subset \De^n$ consisting of those $n$-tuples $\tau\in \De^n$
for which there exists a weighted semistable configuration on $B$ of type $\tau$.

Suppose now that $G$ is a reductive complex Lie group, $K\subset G$ is a maximal
compact subgroup, $X=G/K$ is the associated symmetric space. Then the spaces of weighted
configurations in $\tits X$ of the given type $\tau\in \De^n$ can be identified with products
$$
F=F_1\times ... \times F_n
$$
where $F_{i}$'s are smooth complex algebraic varieties (generalized flag varieties)
on which the group $G$ acts transitively. Hence $G$ acts on $F$ diagonally.

In case $X$ is the symmetric space associated to a complex Lie group, the
notions of stability (semistability, etc.) introduced above coincide
with corresponding notions from symplectic geometry, and, in the case where the
weights $\tau_i$'s
are $L$-integral (i.e., belong to $L=L_{\ul{G}}$)
they also coincide with the concepts of stability (semistability, etc.)
used in the Geometric Invariant Theory; for a proof of this see \cite{KapovichLeebMillson1}.

\subsection{Gauss maps and associated dynamical systems}
\label{Gaussmaps}

We now relate polygons in $X$ (where $X$ is a metric space modeled
on a Euclidean Coxeter complex) and weighted configurations on the
ideal boundary $B$ of $X$, which plays a key role in
\cite{KapovichLeebMillson1} and \cite{KapovichLeebMillson2}. If
$X$ is an infinitesimal symmetric space $\p$ we identify the
visual boundary of $X$ with the Tits boundary $\tits X'$ of the
corresponding symmetric space $X'$ via the exponential map $\p\to
X'$. Thus for all three geometries, the ideal boundary $B$ is a
spherical building.

Consider a (closed) polygon $P=x_1x_2\dots x_n$ in $X$,
i.e. a map $\Z/n\Z\ra X$.
The distances $m_i=d(x_i, x_{i+1})$
determine a finite measure $\nu$ on $\Z/n\Z$ by $\nu(i)=m_i$.
The polygon $P$ gives rise to a collection $Gauss(P)$
of {\em Gauss maps}
\begin{equation}
\label{mapfromindex}
\psi:\Z/n\Z \lra \tits X
\end{equation}
by assigning to $i$ an ideal point $\xi_i\in\tits X$
so that the geodesic ray $\ol{x_i\xi_i}$ (originating at $x_i$ and asymptotic to $\xi_i$)
passes through $x_{i+1}$. This construction, in the case of the hyperbolic
plane, already appears in the letter of Gauss to W.~Bolyai, \cite{Gauss}.
Taking into account the measure $\nu$,
we view the maps $\psi:(\Z/n\Z,\nu)\ra\tits X$
as {\em weighted configurations} of points on $\tits X$.
Note that if $X$ is a Riemannian symmetric space
and the $m_i$'s are all non-zero,
there is a unique Gauss map
due to the unique extendability of geodesics.
On the other hand,
if $X$ is a Euclidean building then, due to the branching of geodesics,
there are in general infinitely many Gauss maps.
However, the corresponding weighted configurations are of the same type,
i.e.\ they project to the same weighted configuration on $\De_{sph}$.

The following crucial observation explains why
the notion of semistability is important for studying closed polygons.

\begin{lem}
[\cite{KapovichLeebMillson1}, \cite{KapovichLeebMillson2}]
\label{polygonstosemistable} For each Gauss map $\psi$ the pushed
forward measure $\mu=\psi_{\ast}\nu$ is semistable. If $X$ is a
symmetric space or an infinitesimal symmetric space then the
measure $\mu$ is nice semistable.
\end{lem}

We are now interested in finding polygons
with prescribed Gauss map.
Such polygons will correspond to the fixed points
of a certain dynamical system.
For $\xi\in \tits X$ and $t\geq0$,
we define the map $\phi:= \phi_{\xi,t}: X\to X$
by sending $x$ to the point at distance $t$ from $x$
on the geodesic ray $\ol{x\xi}$. Since $X$ is nonpositively curved, the map $\phi$ is
1-Lipschitz.
Fix now a weighted configuration
$\psi:(\Z/n\Z,\nu)\ra \tits X$
with non-zero total mass.
We define the 1-Lipschitz self-map
\[
\Phi=\Phi_{\psi}:X\lra X
\]
as the composition $\Phi_n\circ\dots\circ\Phi_1$
of the maps
$\Phi_i=\phi_{\xi_i,m_i}$.
The fixed points of $\Phi$
are the first vertices of closed polygons $P=x_1\ldots x_n$
so that $\psi$ is a Gauss map for $P$.
Since the map $\Phi$ is 1-Lipschitz, and the space in question is complete and
has nonpositive curvature,
the existence of a fixed point for $\Phi$
reduces (see \cite{KapovichLeebMillson2}) to the existence of a bounded
orbit for the dynamical system $(\Phi^n)_{n\in \N}$ formed by the
iterations of $\Phi$. Of course, in general, there is no reason to
expect that $(\Phi^n)_{n\in \N}$ has a bounded orbit: for instance, if the
support of the measure $\mu=\psi_*(\nu)$
is a single point, all orbits are unbounded.

One of  our results is that under the appropriate
semi-stability assumption on $\psi$
the system $(\Phi^n)_{n\in \N}$ has a bounded orbit:

\begin{thm}[\cite{KapovichLeebMillson1}, \cite{KapovichLeebMillson2}]
\label{stabletofixed}
Suppose that $X$ is either a symmetric space or a Euclidean building with one vertex.
Suppose that $\psi$ is a nice semistable weighted configuration on
$\tits X$ (in the symmetric space case)
or a semistable configuration (in the building case).
Then $\Phi_\psi$ has a fixed point.
\end{thm}

This theorem also holds for arbitrary Euclidean buildings: It was proven in an early
version of \cite{KapovichLeebMillson2} under the local compactness assumption, this proof was
superseded by a proof by Andreas Balser \cite{Balser} who had removed the local compactness assumption.

Combining the above result with the Transfer Theorem \ref{transfer} we get:

\begin{thm}
[\cite{KapovichLeebMillson1, KapovichLeebMillson2}]
\label{poly-ss}
Suppose that $X$ is a symmetric space of nonpositive curvature or an infinitesimal symmetric space
or a Euclidean building with a model Weyl chamber $\De$. Then $D_n(X)=\De_{ss}^n(\tits X)$.
\end{thm}

The equivalence of the Problems {\bf Q1} and {\bf Q2} (and consequently
{\bf P1} and {\bf P2}) in the Introduction follows immediately
from the above theorem since for an infinitesimal symmetric space $\p$
and the corresponding symmetric space $X$ the Tits boundaries are the same.
As another corollary of the combination of Theorems \ref{transfer} and \ref{poly-ss}
we get:

\begin{thm}
[\cite{KapovichLeebMillson1, KapovichLeebMillson2}]
\label{mainth}
Let $X$ be either a thick Euclidean building, a symmetric space or an infinitesimal sym\-metric space.
Then $D_n(X)$ depends only on the associated spherical Coxeter complex
and not on the type of the geometry.
\end{thm}

For instance, suppose that $X$ is a nonpositively curved symmetric space and $X'$ is a Euclidean building
which have isomorphic finite Weyl groups and the same rank. Then
$D_n(Cone(\tits X))=\De_{ss}^n(\tits X)= D_n(X)$, $D_n(Cone(\tits X'))=\De_{ss}^n(\tits X')= D_n(X')$,
where $Cone(\tits X)$ and $Cone(\tits X')$ are Euclidean cones over the Tits boundary
of $X, X'$, i.e. 1-vertex buildings.
Since $Cone(\tits X)$ and $Cone(\tits X')$ have isomorphic spherical Weyl groups,
their affine Weyl groups are isomorphic
as well, so the transfer theorem implies that
$$D_n(Cone(\tits X))=D_n(Cone(\tits X'))$$
 as required.

In the case when the group $G$ is complex, for the ideal boundaries $B$ of symmetric spaces $X=G/K$
one constructs the moduli space of semistable weighted configurations on $B$ as follows.

Given a type $\tau$ (so that for each $i$ the vector $\tau_i$ is nonzero), the set of semistable
configurations $Conf_{\tau,sst}(B)$ (resp. nice semistable configurations $Conf_{\tau,nsst}(B)$) of
type $\tau$ on $B$ has a natural topological structure. Define a relation $\sim$ on
$Conf_{\tau,sst}(B)$ by $\psi\sim \psi'$ if
$$
\ol{G\cdot \psi}\cap \ol{G\cdot \psi'}\ne \emptyset.
$$
One then verifies that $\sim$ is an equivalence relation, see \cite{HeinzerLoose} and \cite{Sjamaar}.
Define the moduli space $\M_{\tau}(B)$ of semi-stable
configurations on $B$ of type $\tau$ as the quotient $Conf_{\tau,sst}(B)/\sim$.
We note that the moduli space (see (\ref{moduli})) $\M_{n,\tau}(X)$
of $n$-gons in $X$ with the given
side-lengths $\tau\in \De^n$, is also a compact topological space.

\begin{thm}
The moduli space $\M_{\tau,sst}(B)$ is Hausdorff.
The map $P\mapsto Gauss(P)$ defines a natural
homeomorphism  $h: \M_{n,\tau}(X)\to \M_{\tau,sst}(B)$.
\end{thm}
\proof It was proven in \cite{HeinzerLoose} and \cite{Sjamaar} that $\M_{\tau,sst}(B)=Conf_{\tau,nsst}(B)/G$,
and that $\M_{\tau,sst}(B)$ is Hausdorff. The surjection $h$ is clearly continuous, it is also easily seen to be
injective. Therefore the map
$$
h: \M_{n,\tau}(X)\to Conf_{\tau,nsst}(B)/G.
$$
is also a homeomorphism. \qed

\begin{rem}
For the purposes of this paper we only need to know that
$\M_{\tau,sst}(B)\ne \emptyset$ iff $\M_{n,\tau}(X)\ne \emptyset$.
\end{rem}

Below are few more details concerning the symplectic nature
of the moduli space $\M_{\tau,sst}(B)$ in the case when $G$ is
a complex Lie group. Recall that $\g= \p \oplus \k$
and $\p =\sqrt{-1}\k$.
For $\tau\in \De^n\subset \a^n \subset \p^n$
we identify each $\tau_i$ with the element $\al_i:=\sqrt{-1}\tau_i$
of the Lie algebra $\k$.
We consider the product $M:=\prod_{i=1}^n \O_i$ of the orbits $\O_i:= Ad(K)(\al_i)$, $i=1,...,n$.
The manifold $M$ carries a natural symplectic structure which is invariant under the diagonal adjoint action
of the group $K$. Let $f: M\to \k$ be the {\em momentum mapping} of this action.
Since the map $f$ is given by the formula
$$
f: (\be_1,..,\be_n)\mapsto \sum_{i=1}^n  \be_i,
$$
we obtain an identification of the moduli space
$\M_{n,\tau}(\p)$ of polygons in $\p$ with the {\em symplectic quotient}:
$$
\M_{\tau,sst}(B)=Conf_{\tau}(B)/\!/G:=M/\!/ G:= \{ \xi: f(\xi)=0\}/K.
$$
In the case when all $\tau_i$ belong to the cocharacter lattice $L$, this quotient is the same as the
{\em Mumford quotient} of the projective variety $Conf_{\tau}$ by the group $G$.

\subsection{The polyhedron $D_n(X)$}
\label{polyhedron}

 One of the main results of
\cite{KapovichLeebMillson1} is a description of $D_n(X)$ (where
$X$ is an infinitesimal symmetric space) in terms of the Schubert
calculus in the Grassmannians associated to complex and real Lie
groups $G$ (i.e. the quotients $G/P$ where $P$ is a maximal
parabolic subgroup of $G$). Below we describe $D_n(X)$ for the
three classes of metric spaces considered in the present paper. We
first do it in the context of symmetric spaces of {\em noncompact
type} (i.e. their deRham decomposition contains no Euclidean
factor\footnote{In the context of Lie groups it corresponds to the
case of semisimple algebraic groups.}) since the description in
this class is more natural.

Let $X$ be a symmetric space of noncompact type
and $G$ the identity component of its isometry group.
The ideal boundary $\tits X$
is a spherical building modeled on
a spherical Coxeter complex $(S,W)$
with model spherical Weyl chamber $\De_{sph}\subset S$.
We identify $S$ with an apartment in $\tits X$.
Let $\De$ denote the Euclidean Weyl chamber of $X$.
We identify $\De_{sph}$ with $\tits\De$.

Let $B$ be the stabilizer of $\De_{sph}$ in $G$.
For each vertex $\zeta$ of $\tits X$
one defines the generalized Grassmannian $Grass_{\zeta}=G\zeta=G/P$.
(Here $P$ is the maximal parabolic subgroup of $G$ stabilizing $\zeta$.)
It is a compact homogeneous space
stratified into $B$-orbits called {\em Schubert cells}.
Every Schubert cell is of the form $C_{\eta}=B\eta$ for a unique
vertex $\eta\in W\zeta\subset S^{(0)}$ of the spherical Coxeter complex.
The closures $\ol{C_{\eta}}$ are called {\em Schubert cycles}.
They are unions of Schubert cells and represent well defined elements
in the homology $H_{\ast}(Grass_{\zeta},\Z_2)$.

For each vertex $\zeta$ of $\De_{sph}$ and
each $n$-tuple $\oa\eta=(\eta_1,\dots,\eta_n)$
of vertices in $W\zeta$
consider the following homogeneous linear inequality for $\xi\in\De^n$:
\begin{equation*}
\label{stabineq}
\tag{$\ast_{\zeta;\oa\eta}$}
\sum_i \xi_i\cdot\eta_i \leq0
\end{equation*}
Here we identify the $\eta_i$'s with unit vectors in $\De$.

Let $I_{\Z_2}(G)$ be the subset consisting of all data
$(\zeta,\oa\eta)$
such that the intersection of the Schubert classes
$[\ol C_{\eta_1}],\dots,[\ol C_{\eta_n}]$ in $H_{\ast}(Grass_{\zeta},\Z_2)$
equals $[pt]$.

\begin{thm}
[\cite{KapovichLeebMillson1}]
\label{thm:stabineq}
$\De^n_{ss}(\tits X)\subset\De^n$ consists of all solutions $\xi$
to the system of inequalities {\em (\ref{stabineq})}
where $(\zeta,\oa\eta)$ runs through $I_{\Z_2}(G)$.
\end{thm}

\begin{rem}
This system of inequalities depends on the Schubert calculus
for the generalized Grassmannians $G/P$ associated to the group $G$. It is one
of the results of \cite{KapovichLeebMillson2}
that the set of solutions depends only on the spherical Coxeter complex.
\end{rem}

Typically, the system of inequalities in Theorem
\ref{thm:stabineq} is redundant. If $G$ is a {\em complex} Lie
group one can use the complex structure to obtain a smaller system
of inequalities. In this case, the homogeneous spaces
$Grass_{\zeta}$ are complex manifolds and the Schubert cycles are
complex subvarieties and hence represent classes in {\em integral}
homology. Let $I_{\Z}(G)\subset I_{\Z_2}(G)$ be the subset
consisting of all data $(\zeta,\oa\eta)$ such that the
intersection of the Schubert classes $[\ol C_{\eta_1}],\dots,[\ol
C_{\eta_n}]$ in $H_{\ast}(Grass_{\zeta},\Z)$ equals $[pt]$.

The following analogue of Theorem \ref{thm:stabineq} was proven
independently and by completely different methods in
\cite{BerensteinSjamaar} and in \cite{KapovichLeebMillson1}:

\begin{thm}
[Stability inequalities] \label{thm:stabineqC} $\De^n_{ss}(\tits
X)$ consists of all solutions $\xi$ to the system of inequalities
{\em (\ref{stabineq})} where $(\zeta,\oa\eta)$ runs through
$I_{\Z}(G)$.
\end{thm}

We now consider the general case when $X$ is a symmetric space or a Euclidean building
which splits as $X_0\times X_1$, where $X_0$ is the flat deRham factor of $X$.
In the case when $X=\p$ is an infinitesimal symmetric space we consider the decomposition
$X=X_0\times X_1$ corresponding to the orthogonal decomposition
$\p=\p_0\oplus \p_1$, where $\p_0$ is the Lie algebra of the split part of the
of the central torus of $G$ (where $\g$ is the Lie algebra of $G$).
In this case we again refer to $\p_0:=X_0$ as the Euclidean
deRham factor of $X$.
If $\p=T_p(X')$, where $X'$ is a nonpositively curved symmetric space,
then the above decomposition of $\p$ is the infinitesimal version of the
splitting off the deRham factor of $X'$.

Let $E=E_0\times E_1$ be the corresponding decomposition of the Euclidean Coxeter complex,
the Weyl chamber $\De\subset E$ splits as $\De_0\times \De_1$, where $\De_1\subset E_1$
is the Weyl chamber for the action $W_{sph}\acts E_1$ and $\De_0=E_0$.
Let $D_n(X_i)\subset (\De_i)^n$ denote the side-lengths
polyhedron for the space $X_1$. It is clear that
$$
D_n(X_0)=\{(\si_1,...,\si_n): \sum_{i=1}^n \si_i=0\}.$$
Then we get:
\begin{equation}
\label{split}
D_n(X)=D_n(X_0)\times D_n(X_1).
\end{equation}
We refer the reader to Proposition \ref{prop:reduction} for the explanation.

Combining our results one obtains the following recipe
for determining the polytope $D_n(X)$ for any of the spaces $X$
as in Theorem \ref{mainth}
(i.e., infinitesimal symmetric spaces, symmetric spaces and Euclidean buildings):
Given $X$, first of all, split $X$ as $X_0\times X_1$, where $X_0$ is the
Euclidean deRham factor. In view of the formula (\ref{split}) it suffices to describe
$D_n(X_1)$; so we let $X:= X_1$.
Then find a complex semisimple Lie group $G$ of noncompact type,
whose spherical Coxeter complex is isomorphic to the one of $X$.
Let $X'$ be the symmetric space $G/K$ associated with $G$.
Using Schubert calculus for the Grassmannians $G/P$ associated to
$G$ as in Theorem \ref{thm:stabineqC}
above, compute the system of stability inequalities describing $\De_{ss}^n(\tits X')=D_n(X')$.
The polytopes $D_n(X)$ and $D_n(X')$ are equal.
We note that although the polyhedron
$D_n(X)$ depends only on the spherical Coxeter group,
the system of stability inequalities describing $D_n(X)$
depends on the root system. In fact in \cite{KapovichLeebMillson1} it is shown
that the systems obtained for the root systems $B_3$ and $C_3$
are different (even though they have the same number of inequalities).

\begin{example}
Suppose that $G$ is a reductive algebraic over $\Q$ and $G^{\vee}$ is its
Langlands' dual.
Let $\g$ and $\g^{\vee}$ be the Lie algebras of $G$ and $G^{\vee}$ and
let $X = G/K, X^{\vee}$ be the corresponding symmetric spaces. Let
$\mathfrak{g}= \mathfrak{k} \oplus \mathfrak{p}$ be the Cartan decomposition.
Take a Cartan subalgebra $\a\subset \mathfrak{p}$. The pairing $\<\ ,\ \>$
induces an isomorphism $\a^*\cong \a^{\vee}$, where $\a^{\vee}$ is a Cartan subalgera of $\g^{\vee}$.
We then identify $\a$ with $\a^*$ using the invariant metric. This gives us
an isometry $f: \a\to \a^{\vee}$ which conjugates the spherical Coxeter group
$W_{sph}$ of $\g$ to the spherical Coxeter group
$W_{sph}^{\vee}$ of $\g^{\vee}$. This isometry also carries a Euclidean Weyl chamber
$\De$ of $W_{sph}\acts \a$ onto a  Euclidean Weyl chamber for the action
$W_{sph}^{\vee}\acts \a^{\vee}$. Thus $f$ induces a bijection
$$
f^n: D_n(X)\to D_n(X^{\vee}).
$$
For instance, if $G$ and $G^{\vee}$ are simple complex Lie groups of  type
$B_\ell$ and $C_\ell$ respectively, then the above construction provides
an isometry $D_n(X)\cong D_n(X^{\vee})$.
\end{example}

In the next subsection we will write down the inequalities
\ref{thm:stabineqC} for the root system $B_2$ in the case $n=3$.

\subsection{The polyhedron for the root system $B_2$}\label{inequalities}

In \cite{KapovichLeebMillson1} and \cite{KumarLeebMillson} the stability inequalities for the
polyhedra $D_3(X)$ were computed for all
$X$ of rank $2$ or $3$ (the polyhedron for $G_2$ was computed in
\cite{BerensteinSjamaar}). We now give the example of the polyhedron for the root
system $B_2$. This example will be useful to us later.

Since all symmetric spaces with the same root system give rise to the
same polyhedron we may take $X = SO(5,\C)/SO(5)$. The Weyl chamber
$\Delta$ is given by
$$
\Delta = \{(x,y): x>y>0 \}.$$
Thus the $\De$-lengths $\si_i$ are vectors $(x,y)$ in $\De$.
Since the root system has rank 2,
$\De_{sph}$ has exactly two vertices $\zeta_1, \zeta_2$. Thus we get two
generalized Grassmannians,
$Grass_{\zeta_i}, i=1, 2$. These Grassmannians are
the spaces of isotropic lines and isotropic planes in $\C^4$.
Thus the set of stability inequalities breaks into two subsystems, one for
each $\zeta_i$.
This gives a system of $19$ inequalities in addition to the
inequalities defining the chamber $\De$. Below we have dropped one
of these $19$ inequalities, which was implied by the inequalities defining
the chamber.

The first subsystem of stability inequalities
(corresponding to the Grassmannian of isotropic lines) is given by
\begin{align*}
x_i \leq x_j + x_k , \quad \{i,j,k\}= \{1,2,3\} \\
y_i \leq y_j + x_k ,  \quad \{i,j,k\}= \{1,2,3\}.
\end{align*}
In order to describe the second system we let $S$ be the sum of all the
coordinates of the side-lengths, so
$$
S = x_1 + y_1 + x_2 + y_2 + x_3 + y_3.$$
The second subsystem (corresponding to the Grassmannian of isotropic
planes)
is then given by
$$
x_i  + y_j \leq S/2, \quad 1\leq i,j \leq 3. $$
To the above system of 18 inequalities we also have to add the ``chamber
inequalities'':
$$
x_i \ge y_i, i=1,2,3 \ \hbox{~~and~~} \  y_i \ge 0 , i=1, 2, 3.
$$
In total we get 24 inequalities. Thus we have:

\begin{cor}
Suppose that $X$ is a symmetric space with the (finite) Weyl group of type
$B_2$. Then
there exists a triangle in $X$ whose $\De$-side lengths are vectors
$\si_i=(x_i, y_i)\in \De$, $i=1, 2, 3$
if and only if $\si_i$'s satisfy the above system of 18 inequalities.
\end{cor}

\begin{rem}
Using a computer we have verified that the above system is minimal,
the polyhedron $D_3$ is a cone over a compact
polytope with 15 vertices and 24 top-dimensional faces.
\end{rem}

Note that the Weyl chamber $\Delta\subset \a$ determines a partial order
on $\a$:
$$
\al\le_{\Delta} \be \iff \be- \al \in \Delta.
$$
One may ask if the ``naive'' triangle inequalities
$$
\al \le_\De \be +\ga
$$
are satisfied by the $\De$-side lengths of triangles in a symmetric space
$X$ with the Weyl chamber $\De$. Below is a counter-example:

\medskip
Consider the root system $B_2=C_2$ and the vectors
$$
\al =
\left( \begin{array}{c}
x_1\\
y_1
\end{array}\right),
\quad \be=
\left( \begin{array}{c}
x_2\\
y_2
\end{array}\right)= \ga=
\left( \begin{array}{c}
x_3\\
y_3
\end{array}\right)=
\left( \begin{array}{c}
t\\
t
\end{array}\right),
$$
where $0< t< y_1<x_1< 2t$. The reader will verify that the stability
inequalities above are satisfied by the vectors $\al, \be, \ga$, however
the inequality
$$
\al \le_\De \be +\ga
$$
fails since it would imply that $0<x_1-y_1< (x_2- y_2) + (x_3-y_3)=0$.

\medskip
Nevertheless, one can rewrite the system of stability inequalities
for the root system $B_2=C_2$ in terms of the {\em root cone}
$$
\De^*= \{(x,y): x\ge 0, x+y\ge 0\}= \{v\in \R^2: \forall u\in \De,
u\cdot v\ge 0\}
$$
as follows:
$$
(\al_1, \al_2,\al_3)\in D_3(X) \iff
$$
$$
\al_1, \al_2, \al_3 \in \De,
$$
and
$$
\al_i \le_{\De^*} \al_j + \al_k,
$$
$$
\tau (\al_i)  \le_{\De^*} \tau (\al_j) + \al_k
$$
for all $i, j, k$ such that $\{i, j, k\}=\{1, 2, 3\}$. Here
$\tau(x,y)=(y,x)\in W_{sph}$.

\section{Polygons in Euclidean buildings and the generalized invariant factor problem}
\label{building}

Let $X$ be a thick Euclidean building modeled on a discrete Co\-xeter com\-plex $(E, W_{aff})$.
As in section \ref{Geodesic polygons}, let $L$ be a lattice
in $E$ which contains the translation subgroup
$L_{trans}$ of $W_{aff}$ and which normalizes $W_{aff}$. Note that $L$ acts by automorphisms
on the Coxeter complex $(E_,W_{aff})$.
For the algebraic applications we would like to
determine for which
$\oa{\tau}\in D_n(X) \cap L^n$,
there exists a polygon in $X$ with $L$-integral
side-lengths $\oa\tau$ and the first vertex at a distinguished special vertex  $o\in X$.
In other words, we are interested in the image of the map
\begin{equation}
\label{ref-to-delta}
\iota: D_n^{ref,L}(X)\to D_n(X)\cap L^n.
\end{equation}

In this chapter we will show that in general the map (\ref{ref-to-delta}) is not onto (section
\ref{firstcounterexample}): The counterexamples are based on the idea
of {\em folding} triangles
into apartments (see section \ref{foldtriangle}). In the subsequent chapter \ref{k-computation}
we will prove ``positive results'': Some of them guarantee that
the map (\ref{ref-to-delta}) is onto for certain pairs $(W_{aff},L)$, the other results
establish sufficient conditions for elements of $D_n(X)\cap L^n$ to belong to the
image of $\iota$.

\subsection{Folding polygons into apartments}
\label{foldtriangle}

In this section we describe a construction
which produces billiard triangles in an apartment from triangles in a building .

Suppose that $\De(x,y,z)$ is a triangle
in a (Euclidean or spherical) building $B$. In general it is not contained in an apartment.
However (see \cite[\S 3.2]{KapovichLeebMillson2}) there exists a
finite subdivision of the edge $\ol{xy}$
by points $x_0=x, x_1,...,x_{k-1}, x_k=y$
such that each geodesic triangle $\Delta(z,x_{i},x_{i+1})$ is contained in an apartment in $B$.

\begin{rem}
It is easy to see that there is a uniform upper bound on the number $k$ which depends only
on the Coxeter group $W$.
\end{rem}

For each $i$,
let $a_i\subset B$ be an apartment containing $\De(z,x_i,x_{i+1})$.
We will identify $a_0$ with the model apartment $A$.
We will produce points $x'_i$ in the first apartment $a_0$
such that the triangles $\De(z,x'_i,x'_{i+1})$ are congruent
to the triangles $\De(z,x_i,x_{i+1})$ via apartment isomorphisms
$\al_i:a_i\ra A$.
This is done inductively as follows.
We start with $x'_0=x_0$ and $x'_1=x_1$.
Suppose that $x'_i$ has been constructed.
To find $x'_{i+1}$,
choose $\al_i:a_i\ra a$ so that it carries $\ol{zx_i}$
to $\ol{zx'_i}$.
Put $x'_{i+1}=\al_i(x_{i+1})$.
The procedure yields a {\em billiard triangle}
in the apartment $A$
consisting of two geodesic sides $\ol{zx'_0}$ and $\ol{zx'_k}$
with the same refined lengths as the corresponding sides of
the original triangle $\De(x,y,z)$,
and one piecewise geodesic path $x'_0x'_1\dots x'_k$. The points $x_1',...,x'_{k-1}$
are the {\em break points} of the broken side of this billiard triangle.

\begin{rem}
\label{special0}
If $z$ is a special vertex of $X$,
then, by projecting the apartment $A$ to a Weyl chamber $\Delta$ (with the tip at $z$),
we can assume that the folded triangle is contained in $\Delta$.
\end{rem}

Consider the spherical building $Y_i:=\Si_{x_i}(B)$
whose Coxeter complex is $(\Si_{x'_i}A, W_i)$, where
the group $W_i$ is the stabilizer of $x_i$ in $W$.
Let $\theta_i$ be the canonical projection of $Y_i$ to the Weyl chamber in
this building.

\begin{rem}
Notice that $W_i$ may be smaller than $W_{sph}$.
\end{rem}

Since the refined lengths of $\ol{x_ix_{i+1}}$ and $\ol{x'_ix'_{i+1}}$ are equal,
it follows that $\theta_i(\ola{x_ix_{i-1}})= \theta_i(\ola{x'_ix'_{i-1}})$,
$\theta_i(\ola{x_ix_{i+1}})=\theta_i (\ola{x'_ix'_{i+1}})$.
Because the directions $\ola{x_ix_{i-1}}$ and $\ola{x_ix_{i+1}}$
in the spherical building $Y_i=\Si_{x_i}(B)$ are antipodal, the properties of the
canonical projection $\theta_i$ (see section \ref{geometries}) imply that
the directions $\ola{x'_ix'_{i-1}}$ and $\ola{x'_ix'_{i+1}}$
in the spherical Coxeter complex $(\Si_{x'_i}A, W_i)$
are antipodal modulo the action of $W_i$.

\begin{dfn}
A broken triangle $T\subset A$ with two geodesic sides $\ol{zx'_0}$ and $\ol{zx'_k}$
and one piecewise geodesic path $x'_0x'_1\dots x'_k$ is a {\em billiard triangle} if
at every break point $x'_i$,
the directions $\ola{x'_ix'_{i-1}}$ and $\ola{x'_ix'_{i+1}}$
in the spherical Coxeter complex $\Si_{x'_i}A$
are antipodal modulo the action of the stabilizer $W_i$ of
$x'_i$ in the Coxeter group of $A$.
\end{dfn}

We note that each broken side $x'_0x'_1\dots x'_k$ of a billiard
triangle can be {\em straightened}
in the model apartment $A$, i.e. there exists a geodesic segment
$\ol{x''_0x''_k}\subset A$ (a {\em straightening}
of $x'_0x'_1\dots x'_k$)
such that $x''_0=x'_0$, the metric length of $\ol{x''_0x''_k}$
is the same as of $x'_0x'_1\dots x'_k$,
and the direction of $\ol{x''_0x''_k}$ at $x'_0$ is the same as
the direction of $\ol{x'_0x'_1}$.

\begin{figure}[tbh]
\centerline{\epsfxsize=3.5in \epsfbox{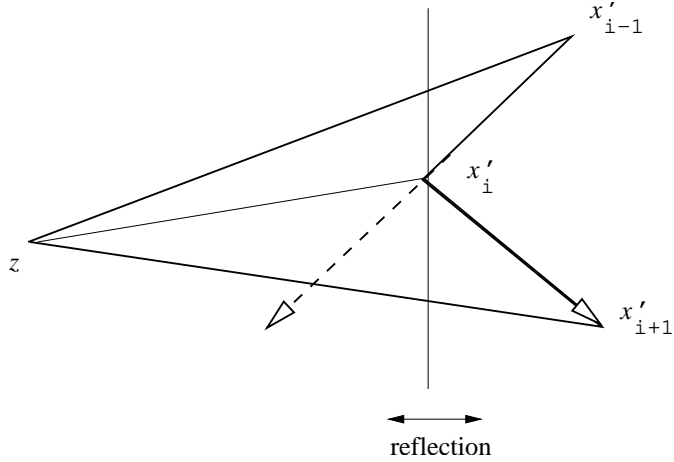}}
\caption{\sl A billiard triangle.}
\label{F8}
\end{figure}

A billiard triangle $T$ in the model apartment
can be unfolded to a geodesic triangle in the building
if and only if, for each break point $x'_i$, $0<i<k$,
the following holds.
Let $\xi'_i,\eta'_i,\zeta'_i\in\Si_{x'_i}A$
be the directions towards $x'_{i-1},x'_{i+1},z$.
The necessary and sufficient condition is:

\begin{condition}
\label{con}
For each $i$ there exists a triangle $\De(\xi_i,\eta_i,\zeta_i)$
in the spherical building $Y_i=\Si_{x'_i}B$
so that $\angle(\xi_i,\eta_i)=\pi$
and the refined lengths of $\ol{\xi_i\zeta_i}$ and $\ol{\eta_i\zeta_i}$
are the same as for $\ol{\xi'_i\zeta'_i}$ and $\ol{\eta'_i\zeta'_i}$ respectively.
\end{condition}

\begin{rem}
\label{rem:con}
A necessary condition for existence of a triangle $\De(\xi_i, \eta_i, \zeta_i)$
is that $\angle(\xi'_i,\zeta'_i)+\angle(\eta'_i,\zeta'_i)\ge \pi$,
which is just the usual (metric) triangle inequality in the spherical building $Y_i$.
\end{rem}

\begin{lem}
\label{unfold}
Suppose that we have an apartment $a'$ in a Euclidean building $X$ and
a billiard triangle $T'\subset a'$,
$T'$ has geodesic sides $\ol{z'x_0'}, \ol{z'x_k'}$ and the broken side
$x_0'x_1'\dots x_k'$. Suppose that
for each break point $x_i', i=1,...,k-1$, there is a wall $H_i\subset a'$
through $x_i'$ which weakly separates\footnote{I.e.,
$H_i$ separates $a'$ into half-apartments $a_+', a-'$ such
that $z$ is in the closure of $a_+'$ and
$x_{i-1}', x_{i+1}'$ are in the closure of $a_-'$.}
$\{x_{i-1}',  x_{i+1}'\}$ from the vertex $z'$.
Assume that the reflection $w$ in the wall
$H_i$ carries the direction $\ov{x_i' x_{i+1}'}$ to the
direction antipodal to $\ov{x_i' x_{i-1}'}$.
Then there exists a geodesic triangle $T=\Delta(x,y,z) \subset X$
so that
$$
\si_{ref}(\ol{z'x'_{0}})= \si_{ref}(\ol{zx}), \si_{ref}(\ol{z'x'_k})= \si_{ref}(\ol{zy})
$$
and the refined side-length of $\ol{xy}$ is the same as
for the straightening of the broken side
${x_0'x_1'\dots x_k'}$.
\end{lem}

\begin{figure}[tbh]
\centerline{\epsfxsize=3.5in \epsfbox{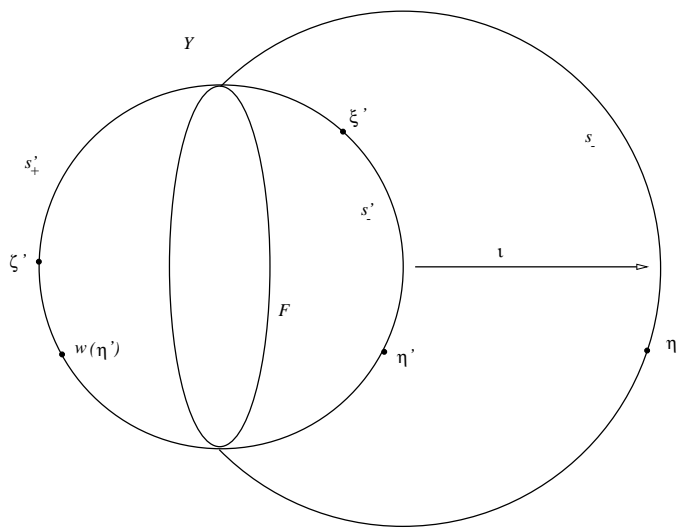}}
\caption{\sl Unfolding.}
\label{unfolding}
\end{figure}

\proof We prove the lemma by verifying Condition \ref{con} for each $i=1,..., k-1$.
The spherical building $Y:= \Si_{x_i'} X$ contains the apartment $s'= \Si_{x_i'}(a')$
and the directions $\xi', \zeta', \eta'$ of the geodesic segments
$$
\ol{x_i' x_{i-1}'}, \ol{x_i' x_{i+1}'}, \ol{x_i' z'}.
$$
The wall $F:=\Si_{x_i'}(H_i)$ in the spherical apartment $s'$
separates $s'$ into half-apartments $s_+', s_-'$, so that $cl(s_+')$ contains $\zeta'$ and $cl(s_-')$
contains the directions  $\xi', \eta'$. Because the building $Y$ is thick, there
exists an apartment $s\subset Y$ which intersects $s'$ along the half-apartment $s_+'$.
Let $s_-:= cl(s\setminus s_+')$.
There exists an isomorphism of Coxeter complexes
$\iota: s'\to s$ which restricts to the identity on $s_+$. Set $\eta:= \iota(\eta'), \zeta:= \zeta', \xi:=\xi'$.
Clearly,
$$
\si_{ref}(\ol{\zeta \xi})= \si_{ref}(\ol{\zeta' \xi'}), \quad \si_{ref}(\ol{\zeta \eta})= \si_{ref}(\ol{\zeta' \eta'}).
$$
It remains to verify that $\angle(\xi, \eta)=\pi$, i.e. that the points $\xi$ and $\eta$ are antipodal.
Note that we also have the third apartment $s''= s'_-\cup s_-$ and the isomorphism of Coxeter complexes
$$
j: s'\to s''
$$
which is the identity on $s_-'$, see Figure \ref{unfolding}.
The restriction of $j$ to $s'_+$ equals $\iota\circ w$.
Thus $j$ carries $w(\eta')$ (which is antipodal to $\xi'$
by assumption) to the point $\eta$, and $j(\xi')=\xi$. Thus the points $\eta, \xi$ are antipodal. \qed


We now state a conjecture describing unfoldable billiard triangles. Let $\Delta\subset E$ be a
Weyl chamber with the tip at a special vertex $o$.
We first define a {\em weak LS path} between two special vertices $x=x_0, y=x_k\in \Delta$
to be a broken geodesic $x_0 x_1 \ldots x_k\subset \Del$ which satisfies Littelmann's
axioms for an LS path in \cite[\S 4]{Littelmann2} except we do not require $dist(\la_{i-1}, \la_i)=1$
in the definition of an $a$-chain. A billiard triangle
is a {\em Littelmann triangle} if the broken side $x_0 x_1 \ldots x_k$ is an LS path.
We say that a billiard triangle in $\Del$ with the geodesic sides
$\ol{ox_0}, \ol{ox_k}$ and a broken side $x_0 x_1 \ldots x_k$ is a {\em generalized
Littelmann triangle} if $x_0 x_1 \ldots x_k$  is a weak LS path.
We recall (see \cite{Littelmann2}) that there exists a Littelmann triangle with geodesic sides
$\ol{ox_0}, \ol{ox_k}$ and the broken side
$x_0 x_1 \ldots x_{k-1} x_k$ iff $V_\ga\subset V_\al \otimes V_\be$ where
$\al:= \si(\ol{o x_0}), \ga=\si(\ol{o x_k})$
and $\be$ is the $\De$-length of the straightening of the broken geodesic $x_0 x_1 \ldots x_{k-1} x_k$,
provided that $\al, \be, \ga$ are characters of the split torus $\ul{T}^\vee \subset \ul{G}^\vee(\C)$.
Here $V_\al, V_\be, V_\ga$ are irreducible representations of the group
$\ul{G}^\vee(\C)$, cf. \S \ref{introduction}.

\begin{conj}
A billiard triangle in $\Del$ with the sides $\ol{ox_0}, \ol{ox_k}$ and
$x_0 x_1 \ldots x_{k-1} x_k$ and with the special vertices $0, x_0, x_k$,
is unfoldable iff it is a generalized Littelmann triangle.
\end{conj}

\subsection{A Solution of Problem {\bf Q2} is not necessarily a solution
of Problem {\bf Q3}}
\label{firstcounterexample}

In this section we first construct an example of a discrete thick Euclidean
building $X$ (modeled on discrete Euclidean Coxeter complex with the root system $R$ of
type $B_2$) and a triangle $P\subset X$ with the $\De$-side
lengths $\tau_1, \tau_2, \tau_3\in L= Q(R^\vee)$ so that:

1. The vertices of $P$ are at vertices of $X$.

2. There is no triangle $P'\subset X$ with vertices at special vertices of $X$ and the
 $\De$-lengths $\tau_1, \tau_2, \tau_3$.

In terms of our basic algebra problems, $(\al, \be, \ga)$ is a solution of
Problem {\bf Q2} but not of Problem {\bf Q3}
for the simply-connected group $\underline{G} = Spin(5)$,

\begin{figure}[tbh]
\centerline{\epsfxsize=3.5in \epsfbox{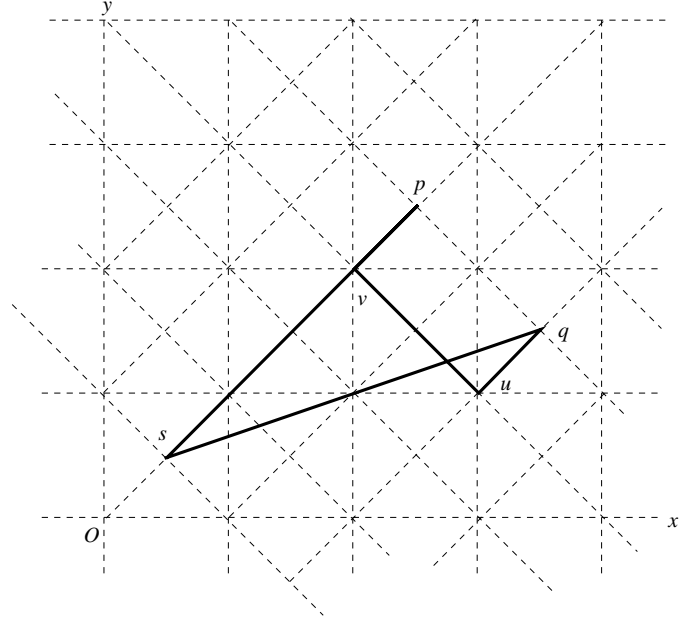}}
\caption{\sl A folded triangle.}
\label{F13}
\end{figure}

We next describe (without proof) analogous counterexamples for the root systems $C_\ell, \ell\ge 2$ and
$L=Q(R^\vee)$.  In terms of our basic algebra problems, this shows that
there are solutions of Problem {\bf Q2} that are not solutions to Problem
{\bf Q3} for the simply-connected groups $\underline{G} = Sp(2\ell),
\ell\ge 2$.

\begin{ex}
\label{notriangle2}
Let $X$ be a thick Euclidean building
with associated discrete Euclidean Coxeter complex $(E,W_{aff})$
of type $B_2$.
Then the map (\ref{ref-to-delta}) is not surjective for $n=3$ and $L= Q(R^\vee)$.
\end{ex}

We recall that the simple roots for the root system $B_2$ are $(1, -1)$ and
$(0, 1)$ and the simple coroots are $(1, -1)$ and $(0, 2)$.
The folded triangle $T$ is represented in Figure \ref{F13}. It has three
vertices: $s, p$ and $q$ and two geodesic sides: $\ol{sp}$ and $\ol{sq}$.
The third side is the broken geodesic segment which consists of three pieces:
$\ol{pv}, \ol{vu}$ and $\ol{uq}$. The vectors in the Weyl chamber
$\De=\{(x,y)\in \R^2: 0\le y\le x\}$ which represent the
corresponding $\De$-side lengths are $(2,2)$, $(3,1)$ and $(2,2)$,
where $(2,2)$ represents the broken side. To see that $T$ can be unfolded
into a geodesic triangle in the corresponding building $X$ one can either
use Lemma \ref{unfold} or simply verify that
the vectors $(2,2)$, $(3,1)$ and $(2,2)$ satisfy the
inequalities of Section \ref{inequalities}.
Note that the vectors $(3,1)$ and $(2,2)$ belong to the coroot lattice.

\begin{figure}[tbh]
\centerline{\epsfxsize=3.5in \epsfbox{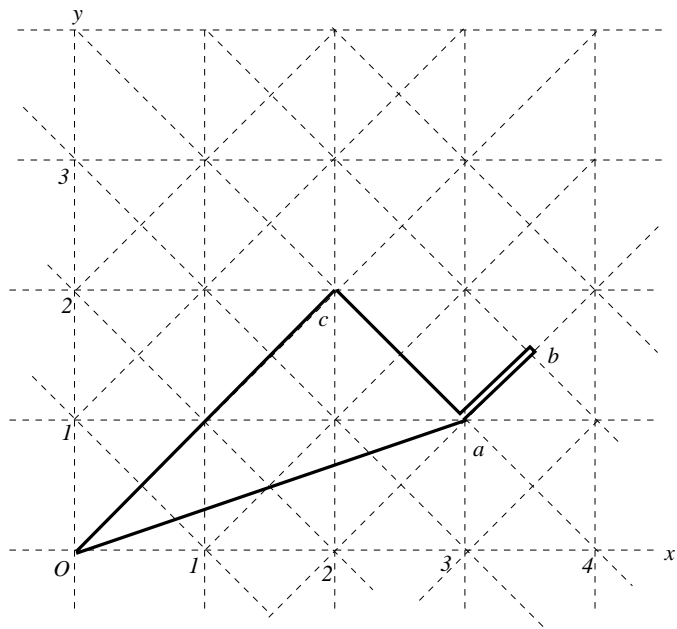}}
\caption{\sl A billiard triangle.}
\label{F15}
\end{figure}

We now prove that there is no triangle $\De(s'',p'',q'')$ in the building $X$
whose $\De$-side lengths are the vectors $(2,2)$, $(3,1)$ and $(2,2)$
and whose vertices are special vertices of $X$. Suppose that such triangle
exists.
As described in section \ref{foldtriangle},
we subdivide the side $\ol{p''q''}$ and fold $\De(s'',p'',q'')$
to a billiard triangle $T'$ in the model apartment $A$. Without loss of generality
we can assume that the side $\ol{s''p''}$ folds onto the geodesic segment
$\ol{oc}$, where $o$ is the origin in $A$ and $c=(2,2)$.
By Remark \ref{special0} we can assume that this triangle is contained in
the Weyl chamber $\Del=\{(x,y): x\ge y\ge 0\}$. Hence the side $\ol{s''q''}$
folds onto the geodesic segment $\ol{oa}$ where $a=(3,1)$.
Note that $d(a, c)$ is strictly
less than $2\sqrt{2}$, which is the magnitude of the vector $(2,2)$. Thus
the broken geodesic segment $f(\ol{p''q''})$ (which
is the image of   $\ol{p''q''}$ under folding) has to have at least one break
point. Next, observe that the only break points in $f(\ol{p''q''})$
can occur at the vertices of the affine Coxeter complex.
Any vertex other than one in the set
\begin{align*}
S:= \{(2, 2), (3, 1), (2.5, 0.5),  (1.5, 1.5),\\
 (2.5, 1.5),  (2.5, 2.5),
(3.5, 0.5), (3.5, 1.5)\}
\end{align*}
would be too far from $a, c$ for a billiard triangle to exist. We note
that all points in $S\setminus \{a, c\}$ are nonspecial vertices
of the Coxeter complex. Thus the only break in $f(\ol{p''q''}$)
which can occur at such a point is ``backtracking''.


Below we exclude breaks at various points of $S$ and will leave the rest
of the possibilities to the reader, since the arguments are similar.
Suppose that the broken geodesic path
$f(\ol{p''q''})$ is the concatenation of the geodesic segments
$$
\ol{ca}, \ol{ab}, \ol{ba},
$$
where $b=(3.5, 1.5)$, see Figure \ref{F15}.
The point $b$ is one of two break points
of $f(\ol{p''q''})$. At the link $\Si_b(A)$ consider the directions:
$$
\xi', \eta', \zeta'
$$
towards the points $a, a$ and $o$. Then
$$
\angle(\zeta', \xi')+\angle(\zeta', \eta')= 2\angle(\zeta', \xi')< \pi.
$$
Hence, according to Remark \ref{rem:con}, the billiard triangle $T'$
cannot be unfolded into a geodesic triangle $\De(s'',p'',q'')\subset X$. Contradiction.

\begin{figure}[tbh]
\centerline{\epsfxsize=3.5in \epsfbox{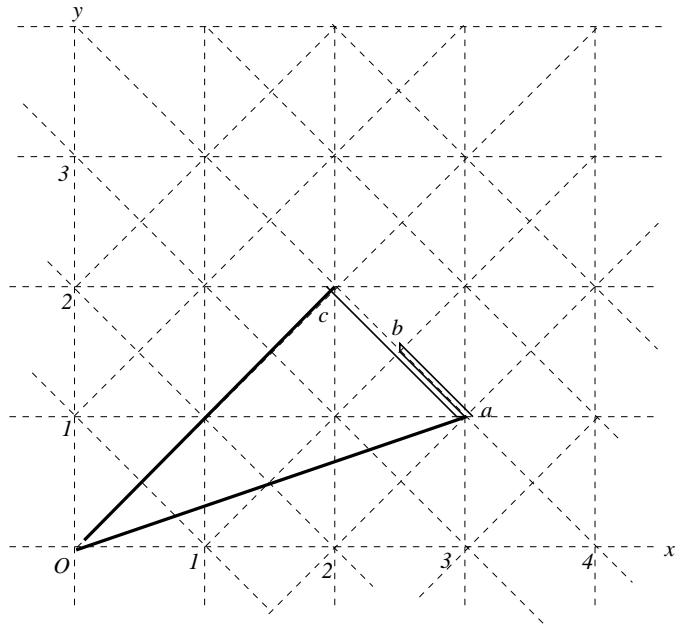}}
\caption{\sl A billiard triangle.}
\label{F16}
\end{figure}

Lastly, consider a break point at the vertex $b=(2.5, 2.5)$;
the broken geodesic path
$f(\ol{p''q''})$ is the concatenation of the geodesic segments
$$
\ol{ca}, \ol{ab}, \ol{ba},
$$
see Figure \ref{F16}. At the link $\Si_a(A)$ of the break point $a$ consider the directions:
$$
\xi', \eta', \zeta'
$$
towards the points $c, b$ and $o$. Again,
$$
\angle(\zeta', \xi')+\angle(\zeta', \eta')= 2\angle(\zeta', \xi')< \pi,
$$
contradiction. (Note that in this example we do not get a contradiction
by considering the link of the break point $b$.) \qed

\begin{rem}
Consider the usual embedding $R_2 \to R_\ell$ of the root system $R_2=B_2$ to $R_\ell=B_\ell$.
The $\De$-side length vectors $(2,2)$, $(3,1)$ and $(2,2)$ correspond to the vectors
$$
\al=\be=(2, 2, 0,..., 0), \ga=(3, 1, 0,...,0)\in L=Q(R^\vee_\ell).
$$
One can show that for each $\ell\ge 3$ and the appropriate discrete building $X$,
the vectors $\al, \be, \ga$ {\bf belong to the image}
of $\iota: D^{ref,L}_3(X) \to D_3(X)$.
\end{rem}

\no Consider now a thick building $X_\ell'$ modeled on the discrete Euclidean Coxeter
complex of type $R_\ell':=C_\ell$, $\ell\ge 3$. Set $L=Q((R_\ell')^\vee)$.
To obtain examples of triples of vectors
which do not belong to the image of
$D^{ref,L}_3(X_\ell') \stackrel{\iota}{\embed} D_3(X_\ell')$ we do the following:

Choose the vectors $\al'=\be'=(2, 0)$ and $\ga'=(2, 1)$.
Then $\al', \be', \ga'$ are in the coroot lattice $L$ of $C_2$.
We claim this choice of side-lengths gives a solution
of Problem {\bf Q2} that is not a solution of Problem {\bf Q3} for the
root system $C_2$.

\begin{lem}
There exists an isomorphism of algebraic groups $\phi:Spin(5) \to Sp(4)$
carrying a split torus of $Spin(5)$ to a split torus of
$Sp(4)$ such that the induced map on Cartan subalgebras
(relative to the coordinates of \cite[pg. 252-255]{Bourbaki}),
is given by the matrix
$$
\frac{1}{2}\left[\begin{array}{cc}
1 & 1\\
1 & -1
\end{array}\right],
$$
\end{lem}

\proof
Let $\phi:Spin(5) \to SL(4)$ be the spin representation. By Theorem G (b)
of \cite{Samelson}, the spin representation of $Spin(2n+1)$ is symplectic
if and only if either $n \equiv 0$ mod $4$ or $n \equiv 1$ mod $4$. Hence,
the image of $\phi$ lies in $Sp(4)$. If follows easily that $\phi$ is
an isomorphism. Let $\epsilon_i, i =1,2$ be the coordinate functionals
as in \cite{Bourbaki}, loc. cit. Since the weights of the spin representation
relative to the basis of the dual of the Cartan subalgebra are $(1/2,1/2), (1/2,-1/2),
(-1/2,1/2)$ and $(-1/2,-1/2)$ we find that the pull-back
by the map, induced on the duals of the Cartan subalgebras, of
the coordinate functional $\epsilon_1 = (1,0)$ is $(1/2,1/2)$ and
the pull-back of the coordinate functional $\epsilon_2 = (0,1)$ is $(1/2,-1.2)$.
Thus the above matrix is the matrix of the map on the duals. Since it is
symmetric is also the matrix of the map on the Cartan subalgebras.
\qed

Since the vectors $\al, \be, \ga$ (from the previous example)
map to the vectors $\al'=\be'$ and $\ga'=(2, 1)$ respectively,
the claim follows.

\medskip
Now consider the natural embedding of
root systems $C_2\embed C_\ell=R_\ell$, $\ell \ge 2$.
One can verify (similarly to the arguments presented in Example \ref{notriangle2})
that the vectors $\al'=\be', \ga'\in L=Q(R^\vee_\ell)$ satisfy the property that
$$
(\al', \be', \ga')\in  D^{L}_3(X_\ell) \setminus \iota( D^{ref,L}_3(X_\ell) ),
$$
where $X_\ell$ is a discrete Euclidean building modeled on
the Euclidean Coxeter complex associated with the root system
$C_\ell$. Thus the triple $(\al', \be', \ga')$ is a solution of
Problem {\bf Q2} for $Sp(2\ell,\mathbb{C})$ but not
a solution of Problem {\bf Q3} for $Sp(2,\mathbb{K})$ where $\mathbb{K}$
is an nonarchimedean local field with value group $\mathbb{Z}$.

\section{The existence of fixed vertices in buildings
and computation of the saturation factors for reductive groups}
\label{k-computation}

As we have seen in the previous chapter, the map
\begin{equation}
\label{ref-to-delta2}
\iota: D_n^{ref,L}(X)\to D_n(X)\cap L^n.
\end{equation}
in general is not surjective.
The goal of this chapter is to find conditions on the root systems, etc., which
would guarantee existence of polygons $P$ in Euclidean buildings $X$ with the
prescribed $L$-integral $\De$-side lengths and vertices at the vertices of $X$.
Moreover, we will find conditions under which the vertices of $P$ are in $G\cdot o$,
where $o\in X$ is a certain special vertex. We will also see that the image of $\iota$
is always contained in the set
$$
D_n^{L,0}(X)=\{(\tau_1,...,\tau_n)\in  D_n(X)\cap L^n : \sum_{i=1}^n \tau_i \in Q(R^\vee)=L_{trans}\}.
$$
We will show that for each data $(W_{aff}, L)$ (where $L_{trans}\subset L\subset N_{aff}$), there exists
a natural number $k$ such that for each $\tau\in D_n(X)\cap L^n$, the vector
$k\tau$ belongs to the image of the map (\ref{ref-to-delta}). We will compute the
{\em saturation factors} $k$ for various classes of $(W_{aff}, L)$. In few cases
we are fortunate and $k=1$, i.e. the map (\ref{ref-to-delta2}) is onto. We then apply these results
to the algebra Problem {\bf Q3}.

\subsection{The saturation factors associated to a root system}
\label{sec:satfactor}

In this section we define and compute {\em saturation factors}
associated with root systems.

\begin{dfn}
Let $(E, W_{aff})$ be a Euclidean Coxeter complex,
$W_{aff}=W_{R,\Z}$. We define the {\em saturation factor} $k_R$
for the root system $R$ to be the least natural number $k$ such
that $k\cdot E^{(0)}\subset E^{(0),sp}= N_{aff}\cdot o$. The
numbers $k_R$ for the irreducible root systems are listed in the
table (\ref{tab}).
\end{dfn}

Below we explain how to compute the saturation factors $k_R$.
First of all, it is clear that if the root system $R$ is reducible
and $R_1,...,R_s$ are its irreducible components, then $k_R=
LCM(k_{R_1},...,k_{R_s})$, where  $LCM$ stands for the {\em least
common multiple}. Henceforth we can assume that the system $R$ is
reduced, irreducible and $n=\dim(V)$. Then the affine Coxeter
group $W_{aff}$ is discrete, acts cocompactly on $E$ and its
fundamental domain (a {\em Weyl alcove}) is a simplex.

Let $\{\al_1,...,\al_n\}$ be the collection of simple roots in $R$
(corresponding to the positive Weyl chamber $\De$) and
$\al_0:=\theta$ be the highest root. Then
\begin{equation}
\label{highestroot} \theta= \sum_{i=1}^n m_i \al_i.
\end{equation}

We can choose as a Weyl alcove $C$ for $W_{aff}$ the simplex
bounded by the hyperplanes $H_{\al_j,0}, H_{\theta,1}$,
$j=1,...,n$. The vertices of $C$ are: $o=x_0$ (the origin) and the
points $x_1,...,x_n$. Each $x_i, i\ne 0$, belongs to the
intersection of the hyperplanes $H_{\al_0,1}, H_{\al_j,0}$, $1\le
j\ne i\le n$. The set of values (mod $\Z$) of the linear
functionals $\al$ ($\al\in R$)  on the vertex set $E^{(0)}$ of the
Coxeter complex, equals $\{\al_i(x_i): i=1,...,n\}$. Note that
$1=\theta (x_i)= m_i\al_i(x_i)$ where the numbers $m_i$ are the
ones which appear in the equation (\ref{highestroot}). Thus
$\al_i(x_i)= \frac{1}{m_i}$.

\begin{lem}
\label{kexists} $k_R= LCM(m_1,...,m_n)$.
\end{lem}
\proof We have: $\al_i(k x_i)\in \Z$ for each $i$, which in turn
implies that $\al(k x_i)\in \Z$ for all $\al\in R, i=1,...,n$.
Hence $\al( k E^{(0)})\subset \Z$ for each $\al\in R$. Since
$N_{aff}= P(R^{\vee})$, this proves that $k E^{(0)} \subset
N_{aff}\cdot o$. If $k\in \N$ is such that $k\cdot E^{(0)}\subset
P(R^{\vee})$, then $m_i$ divides $k$ for each $i=1,...,n$. \qed

In our paper we will also need a generalization of the numbers
$k_R$, which we discuss for the rest of this section. (The reader
who is interested only in simply-connected groups can ignore this
material.) We again consider a general reducible root system $R$.
Suppose that $L'$ is a subgroup of $N_{aff}$ containing the
lattice $L_{trans}=Q(R^\vee)$; we will assume that $L'$ acts as a
{\em lattice} on $E$ (i.e. a discrete cocompact group). Set $L:=
L'\cap V_1\oplus ... \oplus V_s$,  where $V_i$ are the vector
spaces underlying $E_i$. We note that since $L_{trans}\subset L$
and $L_{trans}$ acts as a lattice on $E_1\oplus ... \oplus E_s$,
the discrete group $L$ also acts as a lattice on  $E_1\oplus ...
\oplus E_s$.

Let $p_i$ denote the orthogonal projections $E\to E_i$. Consider
the images $L_i$ of $L$ under the projections $p_i$ ($i=0,...,s$);
since $L\subset N_{aff}$ and $p_i(N_{aff})=N_{aff}^i$, we have the
inclusions
$$
L^i_{trans}\subset L_i \subset N_{aff}^i, i=1,...,s,
$$
where $L^i_{trans}$ is the translation subgroup of $W_{aff}^i$.

\begin{example}
\label{glq} Suppose that $(E, W_{aff})$ is the Coxeter complex
associated with the root system of the group $GL(n)$. Then
$E=\R^n$, $L'=L_{GL(n)}=\Z^n$ is the cocharacter group of the
maximal torus $\ul{T}$ (represented by diagonal matrices) in
$GL(n)$. The group $L'$ is generated by the cocharacters
$e_i=(0,...,0, 1,0...0)$ ($1$ is on the $i$-th place). The coroot
lattice $Q(R^{\vee})$ is generated by the simple coroots
$\al^{\vee}_i=e_{i}-e_{i+1}, i=1,...,n-1$. The metric on $E$ is
given by the trace of the product of matrices. We have the
decomposition $E=E_0\times E_1$ where $E_0$ is 1-dimensional and
is spanned by the vector $e=e_1+...+e_n$, and the space $E_1$ is
the kernel of the map
$$
tr: (x_1,...,x_n)\mapsto \sum_{i=1}^n x_i.
$$
Thus $E_1$ is the (real) Cartan subalgebra of the Lie algebra
$\mathfrak{sl}(n)$ of $SL(n)$, the derived subgroup of $GL(n)$.

The projection $p_1: E\to E_1$ is given by $p_1(u) = u -
\frac{1}{m} tr(u)e$. The group $W_{aff}$ equals $W^1_{aff}$, which
acts on $E_1$ as the Euclidean Coxeter group with the extended
Dynkin diagram of type $\tilde{A}_{n-1}$.

The intersection $L=L'\cap E_1=Q(R^\vee)$, where $R^{\vee}$ is a
root system contained in $E_1$. It is the coroot system of the Lie
algebra $\mathfrak{sl}(n)$ The projection $L_1 = p_1(L')$ is
$P(R^{\vee})$,  the coweight lattice of the Lie algebra
$\mathfrak{sl}(n)$.

\end{example}

Consider the group of isometries $\tilde{W}$ generated by elements
of $W_{aff}$ and $L$. Then $\tilde{W}$ is a Euclidean Coxeter
group with the linear part $W_{sph}$ and translation part $L$,
$\tilde{W}=W_{sph}\ltimes L$. Since $\tilde{W}^i= W^i_{sph}\ltimes
L_i$ normalizes $W_{aff}^i$, for each $i$ we get the induced
action of the finite abelian group
$F_i:=\tilde{W}^i/{W}^i_{aff}\cong L_i/L^i_{trans}$ on the Weyl
alcove $a_i$ of $W_{aff}^i$.

\begin{dfn}
A face $\delta_i\subset a_i$ of a Weyl alcove $a_i$ of
$W^i_{aff}$, will be called {\em $L_i$-admissible} if there exists
an element $g\in F_i$ which preserves $\de_i$ and $\< g\>$ acts
transitively on its vertices ($i=1,...,s$).
\end{dfn}

Note that in the case $L_i=L_{trans}^i$, the only $L_i$-admissible
simplices are the vertices of $\de_i$.

\begin{dfn}
\label{satfactor} For each pair of groups $(W_{aff}^i, L_i)$,
define the {\em saturation factor} $k_i=k(W_{aff}^i,L_i)\in \N$,
to be the smallest natural number $k_i$ such that for each
$L_i$-admissible face $c\subset a_i$, the multiple of its
barycenter $k_i b_c$, belongs to $E^{(0),sp}_i=N_{aff}^i\cdot o$.
We let $k(W_{aff},L)=k(W_{aff}, L')$ denote $LCM(k_1,...,k_s)$. In
the case $L=N_{aff}$ we will use the notation $k_w$ for
$k(W_{aff},L)$.
\end{dfn}

We note that in the case $L=L_{trans}$ we get $k(W_{aff},L)=k_R$
and if $L$ is the weight lattice, $L = P(R^{\vee})$ we have
$k(W_{aff},L)=k_w$.

Our next goal is to compute the saturation factors for various
irreducible root systems and various lattices $L$.

We again assume that the root system $R$ is irreducible and that
its rank $n$ equals $\dim(V)$. Note that the group $F=L/L_{trans}$
acts by automorphisms on the extended Dynkin diagram $\tilde\Ga$
of the root system $R$ (since $F$ acts on the Weyl alcove $a$
which is uniquely determined by the labelled graph $\tilde\Ga$,
whose nodes correspond to the faces of $a$). For $i=1,...,n$ we
mark the $i$-th node (corresponding to $\al_i$) of $\tilde\Ga$
with the natural number $m_i$ which appears the formula for the
highest root (\ref{highestroot}). We mark the $0$-th node of
$\tilde\Ga$ (corresponding to $\theta$) with $1$. Then the
automorphisms of $\t\Ga$ preserve this labelling; the action of
the full group $N_{aff}$ on $\t\Ga$ is transitive on the set of
all the nodes labelled by $1$. Not all automorphisms of
$\tilde\Ga$ can be induced by $F$ even if one takes $L$ as large
as possible, i.e. $L=N_{aff}$. Recall that the action on $a$ comes
from the action of $\tilde{W}$ by conjugation on $W_{aff}$; this
action induces inner automorphisms of the spherical Weyl group
$W_{sph}$. Thus, if $g$ is an automorphism of $\t\Ga$ induced by
an element of $F$ and $g$ fixes a vertex with the label $1$, then
$g$ acts trivially on $\t\Ga$. This does not completely determine
the image of $N_{aff}$ in $Aut(\t\Ga)$ but it will suffice for the
computation of the saturation factors.

Here is the procedure for computing the saturation factor
$k=k(W_{aff},L)$. Given $g\in F$ (including the identity) consider
the orbits of $\<g\>$ in the vertex set of the graph $\tilde\Ga$.
Here and in what follows $\<g\>$ denotes the cyclic subgroup of
$\Isom(E)$ generated by $g$. Let
$\ol{\O}=\{\ol{x}_{i_1},...,\ol{x}_{i_{t}}\}$ be such an orbit.
This orbit corresponds to the orbit $\O=\{x_{i_1},...,x_{i_{t}}\}$
of $\<g\>$ on the vertex set of the Weyl alcove. Take the
barycenter
$$
b(\O)= \frac{1}{t} \sum_{j=1}^t x_{i_j}$$ of the corresponding
vertex set (also denoted $\O$) of the Euclidean Coxeter complex.
For the point $b=b(\O)$ compute the rational numbers $\al_i(b),
i=0,...,n$.

Then find the GCD (the greatest common denominator) of the
rational numbers $\al_i(b), i=0,...,n$, call it $k_{\O}$. Finally,
let
$$
k:= LCM(\{k_{\O}, \hbox{~where $\O$ runs through all orbits of
all~} \<g\>\subset F\}).
$$

\begin{rem}
Instead of taking all $g\in F$ it is enough to consider
representatives of their conjugacy classes in $\tilde{W}/W_{aff}$
(under the conjugation by the full automorphism group of $\t\Ga$).
\end{rem}

It is clear that the  number $k$ computed this way satisfies the
required property:

1. For the barycenter $b$ of each $L$-admissible face of $a$, the
multiple $kb$ belongs to the coweight lattice (which equals
$N_{aff}\cdot o$)

2.  The number $k$ is the least natural number with this property.

\noindent The numbers $k_R$ and $k_w$ are listed in the table
(\ref{tab}) below (the number $i$ in the table is the index of
connection). We will verify the computation in the most
interesting case, namely for the root system of the $A$-type.

\begin{lem}
Suppose that the Dynkin diagram $\Ga$ has type $A_n$ and that
$F\cong \Z/m$. Then $k=k(W_{aff},L)$ equals $m$. In particular, if
$L=N_{aff}$ then we get $k_w=n+1$.
\end{lem}
\proof The group $F\cong \Z/m$ acts on the graph $\t\Ga=\t{A}_{n}$
by cyclic permutations. Let $g\in F$ be a permutation of order
$t$; note that $t$ divides $m$. Then for each orbit $\O$ of $g$
(in the vertex set of the Weyl alcove $a$) we get:
$$
b(\O)= \frac{1}{t} \sum_{x_j\in \O} x_j.
$$
For each $i\ne 0$, $\al_i(b(\O))=0$ if $x_i\notin \O$, and
$\al_i(b(\O))=1/t$ if $x_i\in \O$. For the highest root we get:
$\theta(b(\O))=\frac{t-1}{t}$ if $x_0\in \O$ and $\theta(b(\O))=1$
if $x_0\notin \O$ . In any case, $k_{\O}=t$.

Since all $t$'s divide the order $m$ of the group $F$ (and for the
generator of $F$, $t=m$), the LCM of $k_{\O}$'s taken over all
orbits and all elements of $F$, equals $m$. \qed

Similarly we have

\begin{lem}
Suppose that the Dynkin diagram $\Ga$ has type $D_\ell$ and
$F\cong \Z/2$. Then $k=k(W_{aff},L)$ equals $4$ if $F$ permutes at
least two roots labelled by $2$ and $k=2$ if it does not (the
latter holds for the orthogonal groups).
\end{lem}

We note that for all classical root systems except $D_4$ where
$k_w = k_R$, $k_w$ equals the index of connection $i$ and for all
exceptional root systems, $k_R=k_w$ (so for the computation of
saturation constants for Problem {\bf Q3} $D_4$ behaves like an
exceptional root system).

If $\ul{G}$ is a reductive algebraic group with the cocharacter
lattice $L$ and the associated affine Weyl group $W_{aff}$, we let
$K(G)=k_{inv\; fact}(\ul{G}):= k(W_{aff},L)$. The numbers
$k_{inv\; fact}$ listed below appear in the discussion of the
saturation for the invariant factor problem for the group
$\ul{G}$, see section \ref{k-computation}.

We then obtain:

\begin{equation}
\label{tab} {\mini
\begin{array}{|c|c|c|c|c|c|c|}
\hline
~ & ~ & ~ & ~& ~ & ~& ~ \\
\hbox{Root system} & \ul{G} & \theta & i & k_R & k_w & k(\ul{G})\\
~ & ~ & ~ & ~& ~ & ~& ~ \\
\hline
~ & ~ & ~ & ~& ~ & ~& ~ \\
A_\ell & SL(\ell+1) & \al_1+...+\al_\ell & \ell+1& 1& \ell+1& 1\\
~ & ~ & ~ & ~& ~ & ~& ~ \\
\hline
~ & ~ & ~ & ~& ~ & ~& ~ \\
A_\ell & GL(\ell+1) &~ & ~ & ~ & ~ & 1 \\
~ & ~ & ~ & ~& ~ & ~& ~ \\
\hline
~ & ~ & ~ & ~& ~ & ~& ~ \\
A_\ell & PSL(\ell+1) & ~ & ~ & ~ & ~ & \ell+1 \\
~ & ~ & ~ & ~& ~ & ~& ~ \\
\hline
~ & ~ & ~ & ~& ~ & ~& ~ \\
B_\ell & Spin(2\ell+1), SO(2\ell+1)  & \al_1+ 2\al_2+...+2\al_\ell & 2& 2 & 2& 2\\
~ & ~ & ~ & ~& ~ & ~& ~ \\
\hline
~ & ~ & ~ & ~& ~ & ~& ~ \\
C_\ell & Sp(2\ell), PSp(2\ell) & 2\al_1+ 2\al_2+...+2\al_{\ell-1}+\al_\ell & 2& 2& 2& 2\\
~ & ~ & ~ & ~& ~ & ~& ~ \\
\hline
~ & ~ & ~ & ~& ~ & ~& ~ \\
D_\ell& Spin(2\ell), SO(2\ell)& \al_1+ \al_2+\al_3+ 2\al_4+...+2\al_\ell & 4& 2 & 4 & 2\\
~ & ~ & ~ & ~& ~ & ~& ~ \\
\hline
~ & ~ & ~ & ~& ~ & ~& ~ \\
D_\ell, \ell>4 & PSO(2\ell)& ~& ~ & ~ & ~ &  4\\
~ & ~ & ~ & ~& ~ & ~& ~ \\
\hline
~ & ~ & ~ & ~& ~ & ~& ~ \\
D_4 & PSO(8)& \al_1+ \al_2+\al_3+ 2\al_4 & 4& 2 & 2& 2\\
~ & ~ & ~ & ~& ~ & ~& ~ \\
\hline
~ & ~ & ~ & ~& ~ & ~& ~ \\
G_2 & G & 3\al_1+2\al_2& 1& 6 & 6& 6\\
~ & ~ & ~ & ~& ~ & ~& ~ \\
\hline
~ & ~ & ~ & ~& ~ & ~& ~ \\
F_4 & G & 2\al_1+ 3\al_2+4\al_3 +2\al_4 &  1 & 12 & 12& 12\\
~ & ~ & ~ & ~& ~ & ~& ~ \\
\hline
~ & ~ & ~ & ~& ~ & ~ & ~\\
E_6 & \t{G}, Ad(G) &\al_1+\al_2+ 2\al_3+2\al_4+2\al_5+3\al_6 & 3& 6 & 6 & 6\\
~ & ~ & ~ & ~& ~ & ~& ~ \\
\hline
~ & ~ & ~ & ~& ~ & ~ & ~\\
E_7 & \t{G}, Ad(G) & \al_1+ 2\al_2+2\al_3+2\al_4+3\al_5+& 2 & 12 & 12& 12\\
~ & ~& +3\al_6+ 4\al_7 & ~ & ~& ~ & ~\\
\hline
~ & ~ & ~ & ~& ~ & ~& ~ \\
E_8 & G & 2\al_1+2\al_2+3\al_3+3\al_4+ 4\al_5+ & 1 & 60& 60& 60\\
~&~ & +4\al_6+ 5\al_7 +6\al_8 & ~ &  ~ & ~& ~\\
\hline
\end{array}
}
\end{equation}

Here $\t{G}$ denotes the simply-connected algebraic group, the
symbol $Ad(G)$ denotes the algebraic group of adjoint type, i.e.
the quotient of $\t{G}$ by its center. In the case of root systems
with the index of connection equal to $1$, $Ad(G)=\t{G}$, so we
have used the symbol $G$ to denote the unique connected algebraic
group with the given root system. Note that for the
non-simply-connected classical groups we always get the order of
the fundamental group as the saturation factor (except for the
group $PSO(8)$).

\subsection{The existence of fixed vertices}

We begin this section with few simple remarks about existence of polygons in $X$ with the vertices
of a given type. As before, $G\subset Aut(X)$ is such that
the $G$-stabilizer of the model apartment $A$ acts on $A$ through the group $\tilde{W}$,
whose translational part is $L=L_{\ul{G}}$.
Thus, if $\ol{xy}\subset X$ is a geodesic segment
with $x\in G\cdot o$, then $y\in G\cdot o$ iff $\si(x,y)\in L$.
Therefore the following are equivalent for a polygon $P\subset X$,
whose $\De$-lengths are in $L=L\cdot o$:

1. All vertices of $P$ are in $G\cdot o$.

2. One (say, the first) of the vertices of $P$ is in $G\cdot o$.

We also have

\begin{lem}
\label{special}
For a vector $\oa\tau\in D_n(X)\cap L^n$ the following are equivalent:

1. There exists a polygon $P\subset X$ with the $\De$-lengths $\oa\tau$, such that
all (equivalently, one of) the vertices of $P$ are in $G\cdot o$.

2. There exists a polygon $P\subset X$ with the $\De$-lengths $\oa\tau$, such that
all (equivalently, one of) the vertices of $P$ are special.
\end{lem}
\proof It is clear that (1)$\Ra$(2), let's prove the converse. Let $P=x_1\cdots x_n$ be a
polygon in $X$ with the $\De$-lengths $\oa\tau$, such that
the first vertex $x_1$ of $P$ is special. Let $y_1,...,y_{n+1}\in E$ be such that
$\ol{y_i y_{i+1}}\subset E$ have the same refined lengths
as the segments $\ol{x_i x_{i+1}}$, $i=1,...,n$.
By assumption, the vertex $y_1$ is special.
Recall that the normalizer $N_{aff}$ of $W_{aff}$ in the group of translations of $E$, acts
transitively on the set of special vertices. We identify the model apartment $E$ with an apartment in $X$
containing the origin $o$. Let $T\in N_{aff}$ be such that $T(y_1)=o$. Set
$z_i:= T(y_i)$, $i=1,...,n$. Then $z_i$ are in $L$-orbit of $o$ for each $i$ and according to
the Transfer Theorem \ref{transfer}, there exists a geodesic polygon $Q=u_1\cdots u_n\subset X$
with the refined side-lengths same as for $\ol{z_i z_{i+1}}, i=1,...,n$. In particular,
$\si(u_i, u_{i+1})= \tau_i$ for each $i$ and all vertices of $Q$ are in $G\cdot o$. \qed

\medskip
Thus a vector $\oa\tau$ belongs to the image of the map (\ref{ref-to-delta2})
iff there exists a polygon $P\subset X$ with special vertices and the $\De$-side
lengths $\oa{\tau}$.

As we have seeing in chapter \ref{gau},
the existence problem for polygons with the given $\De$-side lengths in $X$ is equivalent to
the existence of a fixed point for a certain map $\Phi_\psi: X\to X$. In this section we will
try to find conditions under which $\Phi_\psi$ fixes a {\em special vertex} in the building $X$.
We first analyze the case when the Euclidean Coxeter complex is irreducible.

We therefore assume that the affine Coxeter group $W_{aff}$
corresponds to a reduced irreducible root system $R$ of rank
$\ell$ in an $\ell$-dimensional vector space $V$, with the
underlying affine space $E$, see chapter \ref{datum}. Recall that
$L$ is a lattice such that
\begin{equation}
\label{ink}
L_{trans}=Q(R^{\vee}) \subset L\subset N_{aff}=P(R^{\vee}).
\end{equation}
As before, $\tilde{W}$ denotes the subgroup of $\Isom(E)$
generated by elements of $W_{aff}$ and $L$. Then $\tilde{W}$ is a
Euclidean Coxeter group with the linear part $W_{sph}$ and
translation part $L$. Since $\tilde{W}$ normalizes $W_{aff}$, we
get the induced action of the finite abelian group
$\tilde{W}/{W}_{aff}\cong L/L_{trans}=F$ on the Weyl alcove $\al$
of $W_{aff}$. The notion of $L$-admissible faces of Weyl alcoves
of $(E,W_{aff})$ introduced in section \ref{sec:satfactor},
carries over to the building $X$. Note that in the case
$L=L_{trans}$, the only $L$-admissible faces are the vertices of
$X$.

Suppose that
$\tau=(\tau_1,...,\tau_n)$, where $\tau_i\in L$ for each $i$.
Let
$$
\psi=((m_1,\xi_1),...,(m_n,\xi_n))$$
be a weighted configuration on $\tits X$ of the type $\tau$, where $m_i=
\tau_i/|\tau_i|$. Consider the map $\Phi=\Phi_\psi: X\to X$. Then
$\Phi= \Phi_n\circ ...\circ \Phi_1: X\to X$ where each $\Phi_i$ acts on every
apartment $a_i$ asymptotic to $\xi_i$ by the translation $T_i$
(by the vector which has the same $\De$-length
as $\tau_i$). Thus the assumption $\tau_i\in L\subset N_{aff}$ implies that each $T_i$
is an automorphism of the Euclidean Coxeter complex (which we identify with $a_i$).
In particular, the map $\Phi$ preserves the simplicial structure of $X$.

\begin{thm}
\label{fixedbarycenter}
Assume in addition that $\Phi$ has a fixed point in $X$. Then:

1. $\Phi$ fixes the barycenter of an $L$-admissible face in $X$.

2. If $\sum_{i=1}^n \tau_i\in L_{trans}$ then $\Phi$ fixes a vertex of $X$.

3. More generally, if $\hat{L}$ is a lattice normalized by $W_{sph}$,
such that $L_{trans}\subset \hat{L}\subset L$,
and  $\sum_{i=1}^n \tau_i\in \hat{L}$, then $\Phi$ fixes the barycenter of an $\hat{L}$-admissible
face in $X$.

4. If $\Phi$ fixes a special vertex of $X$ then $\sum_{i=1}^n \tau_i\in L_{trans}$.
In particular, if the Coxeter complex associated with $X$ is the irreducible
Coxeter complex $(E,W_{aff})$ of type $A_{\ell}$,
we have:
$$
\sum_{i=1}^n \tau_i\in L_{trans} \iff \Phi \hbox{~~fixes a vertex of~~} X.
$$
\end{thm}
\proof 1. For a fixed point $x=x_1$ of $\Phi$, let $\de$ denote the smallest face  (a simplex)
in $X$ containing $x$. Consider the apartments $a_1\subset X$ (containing $x_1, x_2, \xi_1$),
$a_2\subset X$ (containing $x_2, x_3, \xi_2$),..., $a_n\subset X$ (containing $x_{n}, x_1, \xi_n$).
Then $\Phi_1(\de)\subset a_1\cap a_2$, $\Phi_2\circ \Phi_1(\de)\subset a_2\cap a_3$,...,
$\Phi_n\circ...\circ \Phi_1(\de)\subset a_n\cap a_1$. Since $\Phi$ preserves the simplicial
structure of $X$ and $\de$ is the smallest simplex containing $x$, we get: $\Phi(\de)=\de$.

The map $\Phi|_\de : \de\to \de$ is the composition of maps
$$
(T_n \circ ... \circ T_2\circ T_1)|_\de,
$$
where each $T_i:a_i\to a_i$ is a translation.
Let $\varphi_1: E\to a_1$ be an (isometric) parameterization as in the definition of a space modeled
on a Coxeter complex.
We can assume that $\de':= \varphi_1^{-1}(\de)$ is a face of the Weyl alcove $\alpha\subset E$.
There exists a parameterization $\varphi_2: E\to a_2$ so that
$\varphi_2^{-1}\circ \varphi_1=id$ on the domain of this composition. Similarly (like in the definition
of a developing map of a geometric structure) we choose parameterizations
$\varphi_i, i=1,...,n$, so that for each $i$, $\varphi_{i+1}^{-1}\circ \varphi_i=id$ on the domain of this
composition. Then we get:
$$
(T_n \circ ... \circ T_2\circ T_1)|_\de= (T_n \circ \varphi_n \circ \varphi_{n-1}^{-1}
\circ T_{n-1} \circ ... \circ \varphi_2^{-1} \circ T_2\circ \varphi_2\circ \varphi_1^{-1}\circ
T_1)|_{\de}.
$$
Let $T_i':= \varphi_i^{-1}\circ T_i \circ \varphi_i: E\to E$; these maps
are translations by the vectors $\tau_i'$, which are in the $W_{aff}$-orbits
of the vectors $\tau_i$, $i=1,...,n$. Set $T':= T_n'\circ ... \circ T_1'$.
Then the ``holonomy map''
$$
T'=\varphi_n^{-1}\circ T_n \circ \varphi_n \circ \varphi_{n-1}^{-1}
 \circ ... \circ \varphi_2^{-1} \circ T_2\circ \varphi_2\circ \varphi_1^{-1}\circ
T_1\circ \varphi_1: E\to E$$
sends the simplex $\de'$ to itself. Moreover,
$$
T'|_{\de'}= \varphi_n^{-1}\circ \Phi\circ \varphi_1|_{\de'}=
\varphi_n^{-1}\circ \varphi_1 \circ \varphi_1^{-1} \circ \Phi\circ \varphi_1|_{\de'}=
w\circ \varphi_1^{-1}\circ  \Phi\circ \varphi_1|_{\de'},
$$
where $w\in W_{aff}$.
Thus, by letting $\Phi'|_{\de'}:=\varphi_1^{-1}\circ \Phi\circ \varphi_1|_{\de'}$,
we conclude that $\Phi'|_{\de'}$ admits the extension by $w^{-1}T'$ to the entire
$E$.
Since $w^{-1}T'\in \tilde{W}$, $\Phi$ induces an automorphism
of the simplex $\de$ by an element of the group $\tilde{W}$. We note that
$\Phi$ does not necessarily permute all the vertices of $\de$. Consider however an orbit
of $\<\Phi\>$ on the vertex set of $\de$. It is clear that this orbit spans an $L$-admissible
simplex $c$ in $\de\subset X$; since $\Phi|_\de$ is a 1-Lipschitz automorphism,
it preserves the barycenter of the simplex $c$ and hence the first assertion of the
theorem follows.

2. Let us prove the second assertion. Set
$t=\sum_{i=1}^n \tau_i$. Note that the map $T_i'$ is the translation in $E$
 by the vector $\tau_i'$ where $\tau_i'=w_i(\tau_i)$ for some $w_i\in W_{sph}$.
The composition $T'= T_n'\circ ... \circ T_1'$ is the translation by the vector
$$
v=\sum_{i=1}^n \tau_i'=  \sum_{i=1}^n w_i(\tau_i).
$$
We claim that $v\in L_{trans}$ iff $t\in L_{trans}$. We leave the proof of the following elementary lemma
to the reader:

\begin{lem}
If $\tau\in P(R^\vee)$ then for each $w\in W_{sph}$ we have:
$$
w(\tau)-\tau\in Q(R^\vee).
$$
\end{lem}

\no Therefore
$$
v= \sum_{i=1}^n (w_i(\tau_i) - \tau_i) + \sum_{i=1}^n \tau_i= s+ t,
$$
where $s\in Q(R^\vee)=L_{trans}$. This proves the claim.

\medskip
Recall that $\Phi'|_{\de'}$ admits the extension by $w^{-1}T'$ (where $w\in W_{aff}$)
to the entire model apartment $E$.
Thus, if $t\in L_{trans}$ then  $\Phi'|_{\de'}$  is the restriction
of an element $g\in W_{aff}$. Since the alcove $\al$ is a fundamental domain for the action of $W_{aff}$ on
$E$ and $\de'$ is a face of $\al$, we conclude that  $\Phi'|_{\de'}=id$.
Hence $\Phi$ fixes a vertex of $X$.

3. The proof of 3 is analogous to the proof of 2 and we leave it to the reader.

4. Lastly, suppose that $\Phi$ fixes a special vertex $x$ of $X$, then
$\de=\{x\}$, $\de'=\{x'\}$ is a special vertex in $E$,
and $\Phi'(x')=x'$ implies that $T'(x')=w(x')$, where $w\in W_{aff}$.
Since $x'$ is a special vertex, $w= T\circ w'$, where $w'\in W_{aff}$ fixes
$x'$ and $T\in L_{trans}$. Thus $T'(x')=T(x')$, which implies that $T'\in L_{trans}$,
and hence $t\in L_{trans}$. \qed

\noindent As an immediate corollary of the above theorem we get:

\begin{cor}
\label{fixedvertex}
1. For each $\oa{\tau}=(\tau_1,...,\tau_n)\in  D_n(X)\cap L^n$
such that $\sum_i \tau_i\in L_{trans}$,
there exists a polygon $P\subset X$ such that
$\si(P)=\oa{\tau}$ and the vertices of $P$ are at the vertices of $X$.

2. If there exists a polygon $P\subset X$ with special vertices and $\si(P)=\oa{\tau}$,
then  $\sum_i \tau_i\in L_{trans}$.
\end{cor}

\begin{figure}[tbh]
\centerline{\epsfxsize=3.5in \epsfbox{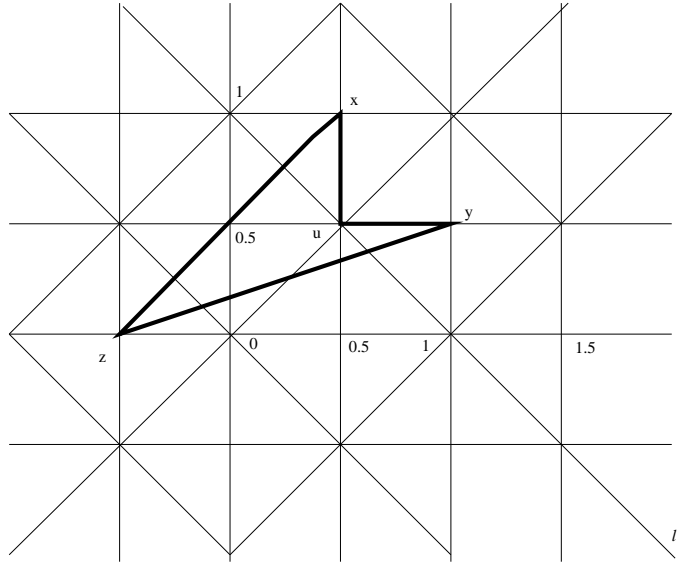}}
\caption{\sl A billiard triangle.}
\label{F4}
\end{figure}

\begin{example}
For the irreducible Coxeter complex $(E,W_{aff})$ of the type $C_{2}$,
there are vectors $\tau_1, \tau_2, \tau_3\in L=P(R^\vee)$ such that
$\Phi$ fixes a vertex of $X$ but $\sum_{i=1}^3 \tau_i\notin Q(R^\vee)=L_{trans}$.
\end{example}
\proof Consider the billiard triangle $T$ shown in Figure \ref{F4}.
It has geodesic sides $\ol{zx}$, $\ol{zy}$
and the broken side $\ol{xuy}$.
Its $\De$-side lengths are $L_{trans}$-integral
and it has non-special vertices.

Let us check that this triangle
can be unfolded to a geodesic triangle $\De(z,x,y')$ in $X$ (with vertices at the vertices of $X$),
in the sense of section \ref{foldtriangle}.
This can be done as follows. We note that the geodesic $l$ (weakly) separates
$z$ from $x$ and $y$. Hence we can unfold the billiard triangle $T$ to a geodesic triangle
in $X$ using Lemma \ref{unfold}.

Then the vectors $\tau_1, \tau_2$ representing
the sides $\ol{xz}$ and $\ol{xy'}$ of $\De(x, y', z)$ are in the coroot lattice $Q(R^\vee)$,
but the vector $\tau_3$ representing the side $\ol{zy}$ is not in the coroot lattice.
Hence $\sum_{i=1}^3 \tau_i\notin Q(R^\vee)$. \qed

Thus the fixed point of $\Phi$ may not be special.
However, if $N_{aff}$ acts transitively on the vertices of $E$,
every vertex is special and we get:

\begin{cor}
\label{transit}
Suppose that
the normalizer $N_{aff}$ of $W_{aff}$ in $V$
acts transitively on the vertices of $E$. Then the image of
$D_n^{ref,L}(X)$ in $D_n(X)$ equals
$$
D_n^{L,0}(X):=\{(\tau_1,...,\tau_n)\in D_n(X)\cap L^n, \sum_{i=1}^n \tau_i\in L_{trans}\}.
$$
\end{cor}

The example of a Coxeter complex for which the hypothesis of Corollary \ref{transit} holds
is given by $(E,W_{aff})$ for which the finite Weyl group $W_{sph}$ is of type $A_\ell$.
Then each vertex of the Coxeter complex $(E, W_{aff})$ is special.

\bigskip
Next we consider the general case when the Coxeter complex $(E,W_{aff})$ is reducible.
As in section \ref{geometries}, we have the corresponding
deRham decomposition of the building $X$:
$$
X= X_0 \times X_1 \times ...\times X_s,
$$
where $X_0$ is the flat deRham factor, the group $W^i_{aff}$ ($i>0$)
acts as an irreducible affine Coxeter group on the Euclidean space $E_i$. The Weyl chamber $\De$
is the product $\prod_{i=0}^s \De_i$ of Weyl chambers, where of course $\De_0=E_0$.
Let $p_i$ denote the orthogonal projection $X\to X_i$, $i=0, 1,...,s$; by abusing
notation we will also use $p_i$ to denote the orthogonal projections
$E\to E_i$. Given a lattice $L'$ in $E$ we get the inclusions
$$
L^i_{trans}\subset p_i(L)=L_i \subset N_{aff}^i, i=1,...,s.
$$
Note that
$$
D_n(X_0)= \{ (\tau_1,...,\tau_n): \sum_{i=1}^n \tau_i=0\}.
$$

We have
\begin{prop}
\label{prop:reduction}
1. $D_n(X)= \prod_{i=0}^s D_n(X_i)$.

2. For each $\oa\tau\in \De^n$, there is a polygon
$P\subset X$ with the $\De$-side lengths
$\oa\tau$ and the vertices at the (special) vertices of $X$ iff
for each $i$ there exists a
polygon $P_i\subset X_i$ with the $\De_i$-side lengths
$p_i(\oa\tau)$ and the vertices at the (special) vertices of $X_i$.

3. $D_n^{ref,L'}(X)=\prod_{i=0}^s D_n^{ref,L_i'}(X_i)$.
\end{prop}
\proof We will prove the third assertion since the proofs of (1) and (2) are similar.
First of all, if $P\subset X$ is a polygon with $L'$-integral refined side-lengths then
its projections $P_i:= p_i(P)\subset X_i$ are also
polygons with $L_i$-integral refined side-lengths. Conversely, if
$P_i\subset X_i$ are polygons with $L_i'$-integral refined side-lengths
and $x_{ij}$ is the $j$-th vertex of the polygon $P_i$, we set
$$
x_j:= (x_{0j},..., x_{sj})\in X,
$$
this point is a vertex of $X$ since all $x_{ij}$'s
are vertices of the corresponding buildings $X_i$.
This defines the polygon $P=x_1 \cdots x_n$ in $X$.
It is clear that the refined side-lengths of
$P$ are $L'$-integral. \qed

The above proposition implies that Corollaries
\ref{fixedvertex} and \ref{transit} remain valid
for Euclidean buildings with reducible Coxeter complexes.

As in section \ref{sec:satfactor}, given a lattice $L'$ in $E$ we
define the lattice $L$ as the intersection of $L'$ and the
translation group of $E_1\oplus ...\oplus E_s$. We now can compute
the image of the map (\ref{ref-to-delta2}) ``up to saturation'':

\begin{thm}
\label{sat} 1. Let $k=k(W_{aff},L')=k(W_{aff},L)$ be the
saturation factor defined in section \ref{sec:satfactor}. Then for
each $\oa\tau\in D_n(X)\cap (kL')^n$, there exists an $n$-gon
$P\subset X$ with special vertices and  $\De$-side lengths
$\tau_1,...,\tau_n$.

2. Let $R$ be the root system corresponding to the Coxeter complex
$(E,W_{aff})$ and let $k:=k_R=k(W_{aff},L_{trans})$ be the
saturation factor defined in section \ref{sec:satfactor}. Then for
each $\oa\tau\in D_n(X)\cap (kL')^n$ satisfying
$$
\sum_i \tau_i \in L_{trans}=Q(R^\vee)
$$
there exists an $n$-gon $P\subset X$ with
special vertices and  $\De$-side lengths $\tau_1,...,\tau_n$.
\end{thm}
\proof 1. We first consider the case when the
root system $R$ is reduced, irreducible and has rank equal the dimension of the space $V$; then $L=L'$.
Let $P'\subset X$ be a polygon with the $\De$-lengths $\frac{1}{k}\tau_1,...,\frac{1}{k}\tau_n$.
Since $\frac{1}{k}\tau_i\in L$ for each $i$,
Theorem \ref{fixedbarycenter} implies that the polygon $P'$ can be chosen so that its first
vertex is the  barycenter of an $L$-admissible simplex in $X$.
Let $\ol{x_i x_{i+1}}$, $i=1,2,...,n$, be oriented segments in the model apartment $E$,
which represent the refined side-lengths of the polygon $P'$; the point $x_1$
(and hence $x_{i}$ for each $i=2,...,n+1$) is the barycenter of an $L$-admissible simplex in $E$.
We regard $E$ as the vector space $V$ by identifying the origin with the special vertex $o$.
Then, according to the definition of $k=k(W_{aff},L)$, the segments
$\ol{(k x_i) (k x_{i+1})}=\ol{y_i y_{i+1}}$ have end-points at the {\em special vertices} of $E$.
Consider the affine transformation $\iota: E\to E$ which is the homothety
$x\mapsto kx$. Then $\iota: (E,W_{aff}) \to (E,W_{aff})$ is
an injective endomorphism of the affine Coxeter systems
and the Transfer Theorem \ref{transfer} implies that
there exists a polygon $P$ in $X$ whose refined side-lengths are  represented by
$\ol{y_i y_{i+1}}$, $i=1,...,n$. The vertices of $P$ are special and
its $\De$-side-lengths are $\tau_1,...,\tau_n$.

We now consider the general case when the Coxeter complex $(E,W_{aff})$ is reducible.
Let $P'\subset X$ be a polygon with the $\De$-lengths
$\frac{1}{k}\tau_1,...,\frac{1}{k}\tau_n$.
Since $p_0( \sum_{i=1}^n \tau_i )=0$, we have
$$
\frac{1}{k}\sum_{i=1}^n \tau_i \in L'\cap V_1\oplus ... \oplus V_s= L.
$$
Hence for each $j=1,...,s$ we get
$$
\frac{1}{k}\sum_{i=1}^n p_j(\tau_i) \in L_j=p_j(L).
$$
Therefore, by Theorem \ref{fixedbarycenter} (part 3, where we take $\hat{L}=L_j$, $j=1,...,s$),
the polygon $P'$ can be chosen so that for each $j=1,...,s$,
the vertices of $p_j(P')$ are barycenters of $L_j$-admissible faces of $X_j$.
Recall that $k=LCM(k_1,...,k_s)$.

By the irreducible case, for each $j=1,...,s$ there is a polygon $P_j'\subset X_j$
with special vertices and side-lengths $\frac{k_j}{k}p_j(\tau_1),...,p_j(\tau_n)$.
Let $\ol{x_{i,j} x_{i+1,j}}$, $i=1,2,...,n, j=1,...,s$, be oriented segments in the model apartment $E_j$,
which represent the refined side-lengths of the polygon $P'_j$.
It follows from the Transfer Theorem \ref{transfer} that
for each $j=1,...,s$ there exists a polygon $P_j\subset X_j$ whose refined
side-lengths are represented by the segments $\ol{(m_j x_{i,j}) (m_j x_{i+1,j})}$, where $m_j=k/k_j$.
Note that  $m_j \frac{k_j}{k}p_j(\tau_i)= p_j(\tau_i)$ for each $i$ and $j$.
By the definition of $k_j$, for each $j$ the vertices of $P_j$ are special vertices of $X_j$.
Now Proposition \ref{prop:reduction} implies that
there exists a polygon $P\subset X$  with the special vertices,
whose projections to $X_j$ are $P_j$ for each $j$ and the $\De$-side lengths of $P$ are $\tau_1,...,\tau_n$.

2. The proof of this assertion is analogous to the proof of 1. First, suppose that the Coxeter complex
is irreducible. Then, according to Theorem \ref{fixedbarycenter} (part 2),
a polygon $P'\subset X$ with the $\De$-lengths $\frac{1}{k}\tau_1,...,\frac{1}{k}\tau_n$
can be chosen so that its vertices are vertices of $X$.
Let $\ol{x_i x_{i+1}}$, $i=1,2,...,n$, be oriented segments in the model apartment $E$,
which represent the refined side-lengths of the polygon $P'$.
Then, by the definition of $k=k_R$, the end-points of the segments $\ol{(k x_i) (k x_{i+1})}$
are special vertices of $X$. The rest of the argument is the same as for (1). \qed

\subsection{Saturation factors for reductive groups}

We now apply the results from the previous section
to the generalized Invariant Factor Problem for nonarchimedean reductive
Lie groups $G$.
Suppose that $\ul{G}$ is a connected reductive algebraic Lie group over $\K$, where $\K$ is a field
with discrete valuation $v$. Let $G:= \ul{G}(\K)$.
The group $\ul{G}$ determines a Bruhat-Tits building $X$,
Bruhat-Tits root system $R$, the corresponding affine Coxeter group $W_{aff}$ and
the extended cocharacter lattice $L:=L_{\ul{G}}$ satisfying the double inclusion (\ref{ink}).
Define $k_{inv\; fact}(G):=k(W_{aff},L)$.

As an immediate corollary of Theorem \ref{sat}, we obtain:

\begin{cor}
\label{simplesaturation}
1. For $k=k_{inv\; fact}(G)$ and any
$\tau \in L^n \cap D_n(X)$,
there exists an $n$-gon $P$ in $X$ with
$\De$-side-lengths equal to $k\tau$ and vertices in $G \cdot o$.

2. For $k=k_R$ and any
$$
\oa\tau \in D^{L,0}_n(X)= \{ \oa\tau \in D_n(X)\cap L^n : \sum_{i=1}^n \tau_i \in Q(R^\vee)\},
$$
there exist an $n$-gon $P$ with
$\De$-side lengths equal to $k\tau$ and vertices in $G \cdot o$.

3. If $\ul{G}$ is semisimple and simply-connected
then $D^{L,0}_n(X)=D^L_n(X)$ and\newline $k_{inv\; fact}(G)=k_R$.
\end{cor}
\proof Parts (1) and (2) follow from
Theorem \ref{sat} (parts (1) and (2) respectively) and Lemma \ref{special}.

To prove (3) note that in this case $Q(R^\vee)=L$ and thus
$$
k_{inv\; fact}(G)=k_R=k(W_{aff}, Q(R^\vee)). \qed
$$

\begin{rem}
We emphasize that $k_{inv\; fact}(G)$ and $k_R$ depend only on $G$
and not the choice of $\tau$.
\end{rem}

By specializing to the case $n=3$ and using the equivalence (see \S
\ref{Geodesic polygons}) of Problem RGTI for Euclidean buildings
with the Problem {\bf Q3}, we get:

\begin{cor}
\label{simplesaturation2}
For each $k\in \N$ divisible by $k_{inv\; fact}(G)$ and each
triple $(\al, \be, \ga)\in k(L\cap \De)^3$
which satisfies the {\em generalized triangle inequalities} for the space $X$, there
exist representatives $A, B, C$  of the double coset classes
$$
\ul{G}(\O) \al(\pi) \ul{G}(\O), \quad \ul{G}(\O) \be (\pi) \ul{G}(\O), \quad
\ul{G}(\O) \ga (\pi) \ul{G}(\O),
$$
such that $ABC=1$.
\end{cor}

\begin{example}
Let $\ul{G}= Spin(5)$. Then  Example \ref{notriangle2} shows that there
is a triple
$(\al, \be, \ga)\in (L\cap \De)^3, L=Q(R^\vee)$,
such that $\al, \be, \ga$ satisfy the generalized triangle inequalities and such that one cannot find $A, B, C$ in
the corresponding double cosets of $Spin(5, \K)$ such that $ABC=1$.

Thus one cannot take $k=1$ for the case $\underline{G} = Spin(5)$.
\end{example}

\begin{example}
$k_{inv\; fact}(GL(m))=1$. In other words, let $X$ be the
Bruhat-Tits building associated with the group $GL(m, \K)$,
where $\K$ is a field with nonarchimedean discrete
valuation and the value group $\Z$.
Let $\oa{\tau}\in \Z^m$ have entries arranged in decreasing order.
Then $\oa{\tau}$ satisfies the generalized triangle inequalities
(i.e. $\oa{\tau}\in D_n(X)$) if and only if there exists an
$n$-gon $P$ in $X$
with   the first vertex at the origin (the point stabilized by $GL(m, \O)$) and
$\De$-side lengths $\oa{\tau}$.
\end{example}
Although this example is a special case of Theorem \ref{sat} (part 2),
we present a complete proof for the benefit of the reader.

\proof Recall (see  Example \ref{glq}) that the Euclidean Coxeter complex $(E,W_{aff})$ of $X$
is reducible, $E=E_0\oplus E_1$, where $E_1$ is the kernel of the map
$$
tr: (x_1,...,x_m)\mapsto x_1+...+x_m,
$$
and $E_0$ is the span of $(1,...,1)$.
Let $V_1$ be the vector space underlying $E_1$.
We may identify $V_1$ with the Lie algebra of traceless real diagonal matrices,
the real Cartan subalgebra of the Lie algebra of $SL(m)$.
Moreover (see  Example \ref{glq}), $L=L_{GL(m)}=\Z^m$, is the cocharacter lattice;
the orthogonal projection of $L$ to
$V_1$ is the coweight lattice of $SL(m)$, and the root system is of type $A_{m-1}$.
On the other hand, suppose that
$\oa\tau=(\tau_1,...,\tau_n)$. Then, since $\oa\tau\in D_n(X)$, there exists a
polygon $\hat{P}$ in $X$ with the $\De$-side lengths
$\tau_i$. Let $X=X_0\times X_1$ be the deRham decomposition
of $X$, where $X_0=E_0\cong \R$ is the flat factor.
Since the projection of
$\hat{P}$ to the flat deRham factor $E_0$ of $X$ is also a polygon, we get:
$$
p_0(\tau_1+...+\tau_n)=tr(\tau_1+...+\tau_n)=0.$$
Hence $\tau_1+...+\tau_n\in L\cap E_1=Q(R^\vee)$.
Therefore we can apply Theorem \ref{fixedbarycenter} (part 2 or 4)
 to conclude that there exists a polygon $P'\subset X_1$ whose vertices are
(special) vertices of $X_1$ and whose $\De$-side lengths are $p_1(\tau_1),..., p_1(\tau_n)$.
By combining this with the polygon $p_0(\hat{P})$
via Proposition \ref{prop:reduction} we conclude
that there exists a polygon $P\subset X$ whose vertices are (special) vertices of $X$
and whose $\De$-side lengths are $\tau_1,..., \tau_n$.
Thus $k_{inv\; fact}(GL(m))=1$. \qed

\begin{thm}
1. Let $\ul{DG}$ be the derived (i.e. commutator) subgroup of
$\ul{G}$. Then
\begin{equation}
\label{deriv}
k_{inv\; fact}(\ul{G})= k_{inv\; fact}(\ul{DG}).
\end{equation}
2. Moreover, if $\ul{G_1},...,\ul{G_s}$ are the simple factors of $\ul{DG}$ then
$$
k_{inv\; fact}(\ul{G})= LCM( k_{inv\; fact}(\ul{G_1}),..., k_{inv\; fact}(\ul{G_s}) ).
$$
\end{thm}

\proof
The equation (\ref{deriv}) follows from three observations:

We first claim that the lattice $L'=L_{\ul{G}}\cap V_1\oplus ... \oplus V_s$ equals
$L_{\ul{DG}}$. This is because we can assume that the maximal split tori
$\ul{T}_{\ul{DG}}$ and $\ul{T}_{\ul{G}}$ of $\ul{DG}$ and $\ul{G}$ are related by
$$
\ul{T}_{\ul{DG}}=\ul{T}_{\ul{G}} \cap \ul{DG}.$$
Hence the cocharacter lattice of $\ul{DG}$ is the sublattice of $L_{\ul{G}}$
consisting of those cocharacters whose image is contained in $\ul{DG}$. Passing to tangent
vectors at the identity yields the claim.

Second, by definition, the coroot system $R_{\ul{DG}}^{\vee}$ of $\ul{DG}$ is the same as the coroot system $R_{\ul{G}}^{\vee}$
of $\ul{G}$ regarded as a root system in $E_1\oplus ...\oplus E_s$. Thus $\ul{G}$ and $\ul{DG}$
have the same affine Weyl group $W_{aff}$.

Lastly, the projection of the lattice $L_{\ul{G}}$ into $V_1\oplus ... \oplus V_s$ is
contained in $P(R^\vee)\cap V_1\oplus ... \oplus V_s$.

With these observations, the same argument as in the case of $GL(m)$ goes through:
$$
k(W_{aff}, L_{\ul{G}})= k(W_{aff}, L_{\ul{DG}})$$
(the saturation factor $k_{inv\; fact}$ depends only on the group
$W_{aff}$ and the cocharacter lattice), which implies (\ref{deriv}).

To prove the second assertion of the Theorem note that the Euclidean Coxeter complex for
$\ul{DG}$ is the direct product of the Euclidean Coxeter complexes for its factors $\ul{G_i}, i=1,...,s$.
Let $\ul{T}_i$ be maximal split tori of $\ul{G}_i$, $i=1,...,s$ and $L_i=X_*(\ul{T}_i)$
be their cocharacter lattices. Let $\ul{T}_{\ul{G}}:= \prod_i \ul{T}_i$, and $L':= X_*(\ul{T}_{\ul{G}})$;
then each $L_i$ is the projection of $L'$ to $V_i$, where $(V_i, W_{aff}^i)$ are the Euclidean
Coxeter complexes associated with $\ul{G_i}$.

Then Definition \ref{satfactor} implies that
$$
k(W_{aff}, L')= LCM( k(W_{aff}^1, L_1), ...., k(W_{aff}^s, L_s)).  \qed
$$

\bigskip
We can also give a complete solution (in terms of the stability inequalities)
for the Invariant Factor Problem {\bf Q3} in the case of
nonarchimedean reductive Lie groups $G$ with root system of type $A_{m-1}$.
For instance, $\ul{G}$ could be the quotient of  $SL(m)$ by a subgroup of its center.
As before, let $L$ denote the lattice $L_{\ul{G}}$ and
let $X:= X_{\ul{G}}$ be the Bruhat-Tits building associated with $\ul{G}$.

\begin{thm}
\label{atype}
Let $\ul{G}$ be a reductive algebraic group over $\K$,
so that the group $\ul{G}$ has the associated root system $R$ of type $A_{m-1}$.
Then the natural embedding
\begin{equation}
\label{embinteger1}
D_3^{ref,L}(X)\embed D_3^L(X)= \{
(\al, \be, \ga)\in D_3(X)\cap {L}^3 : \al+\be+\ga\in Q(R^{\vee})\}
\end{equation}
is onto.
\end{thm}
\proof The assertion follows from Corollary \ref{simplesaturation}
and the fact that $k_R=1$ for the root system of type $A$. \qed

Specializing to the case of $GL(m)$ and $SL(m)$ we get the following corollary:

\begin{cor}
Let $X={\mathbb P}_{m}$ be the symmetric space of symmetric positive-definite
$m\times m$-matrices. Then
there exists a solution to the Invariant Factor Problem {\bf Q3} for the case
$G=SL(m,\K)$ (or $GL(m,\K)$) and $K=SL(m,\O)$ (resp. $GL(m,\O)$)
if and only if $\alpha$, $\beta$ and $\gamma$
are integer diagonal matrices so that $tr(\al+\be+\ga)=0$ and the
projections
$$\al- \frac{1}{m} tr(\al)I, \be- \frac{1}{m} tr(\be)I,
\ga- \frac{1}{m} tr(\ga)I$$
satisfy the stability inequalities \ref{thm:stabineqC} for the space $X$. Here $I$ is the
identity matrix.
\end{cor}

\section{The comparison of Problems Q3 and Q4}
\label{3&4}

\subsection{The Hecke ring}

In this section we will (for the most part) follow the notation of
\cite{Gross}. We urge the
reader to consult this paper for more details. However we will let
$\K$ and $\mathcal{O}$ be as in section \ref{tough}.
We will assume that the valuation $v$ is discrete, the field $\K$ is locally compact
and the residue field is finite of order $q$ and uniformizer $\pi$.

\begin{rem}
The assumption that $\K$ is locally compact is equivalent to $\K$ being a totally-disconnected local field,
see \cite[pg. 5]{Taibleson}. By applying Theorem \ref{transfer},
it follows that the main result of this section, Theorem \ref{q4->q3},
also holds in the case when $\K$ is not locally compact, for example the case in which
 $\K = \C ((t))$ the field of fractions of the ring of formal
power series $\C[[t]]$.
\end{rem}

We let $\ul{G}$ be a
connected reductive algebraic group over $\K$. We will assume that
$\ul{G}$ is split over $\K$. Then $\underline{G}$ is the general fiber of
a group scheme (also denoted $\underline{G}$) over $\mathcal{O}$ with reductive
special fiber. We fix a maximal split torus $\underline{T}\subset
\underline{G}$ defined over
$\mathcal{O}$.
We put $G = \ul{G}(\K), K= \ul{G}(\mathcal{O})$ and
$T =\ul{T}(\K)$.
We let $\ul{B} \subset \underline{G}$ be a Borel subgroup normalized
by $\underline{T}$ and put $B= \underline{B}(\K)$. We let $U$ be the unipotent
radical of $B$ whence $B = TU$. Let $X=X_{\ul{G}}$ denote the Bruhat-Tits building
associated with the group $\ul{G}$; $o\in X$ is a distinguished special vertex stabilized by
the compact subgroup $K$.

We have already defined in section \ref{tough} the free abelian groups
(of rank $l= \dim(\underline{T})$)
$X_*(\underline{T})$ and $X^*(\underline{T})$ and a perfect $\Z$-valued pairing
$\<\ ,\ \>$ between them. The first contains the coroots $R^{\vee}$, the second
contains the roots $R$ of $\underline{G}$.

The roots are the characters of $\underline{T}$
that occur in the adjoint representation on the Lie algebra of $\underline{G}$.
The subset $R^+$ of the roots that occur in representation on the Lie algebra
of $\underline{B}$ forms a positive system and the indecomposable elements
of that positive system form a system of simple roots $\Pi$. We let $W$ denote the
corresponding (finite) Weyl group.

The root basis $\Pi$ determines a positive Weyl chamber $P^+$ in
$X_*(\underline{T})$, by
$$
P^+ = \{\lambda \in X_*(\underline{T}): \ \<\lambda,\alpha\> \geq 0, \alpha \in \Pi\}.
$$
We define a partial ordering on $P^+$ by $\la > \mu$ iff the difference $\la-\mu$
is a sum of positive coroots.

\medskip
We will use the notation $0$ for the trivial cocharacter.

\medskip
We define the element $\rho \in X^*(\underline{T})\otimes \Z[1/2]$ by
$$
2\rho = \sum\limits_{\al\in R^+} \alpha.$$
We recall that $\rho$ is the sum of the fundamental weights of $R$.
Then, for all nontrivial $\lambda \in P^+$,  the half-integer
$\<\lambda,\rho\>$ is positive.

\begin{dfn}
The Hecke ring $\mathcal{H} = \mathcal{H}(G,K)$ is
the ring of all locally constant, compactly supported functions $f: G\lra \Z$
which are $K$-biinvariant. The multiplication in $\mathcal{H}$ is by the
convolution
$$
f \cdot g (z) = \int_G f(x) \cdot g(x^{-1}z)dx$$
where $dx$ is the Haar measure on $G$ giving $K$ volume $1$.
\end{dfn}

We claim that the function $f\cdot g$ also takes values in $\Z$.
Indeed, $f$ and $g$ are finite sums of characteristic functions of
$K$-double cosets. Thus it suffices to prove the claim in the
case that $f$ and $g$ are both characteristic functions of
$K$-double cosets. In this case it is immediate that their
convolution product is the characteristic function of the set
of products of the elements in the two double cosets. This product
set is itself a finite union of $K$-double cosets. This implies
that the structure constants $m_{\la,\mu}(\nu)$
of the Hecke ring, defined below, are nonnegative integers.
The characteristic function of $K$ is the unit element of $\mathcal{H}$.
For the proof of the next lemma see \cite[\S 2]{Gross}.

\begin{lem}
$\mathcal{H}$ is commutative and associative.
\end{lem}

In fact much more is true. For $\lambda \in X_*(\underline{T})$ let
$c_{\lambda}$ be the characteristic function of the corresponding
$K$-double coset $\la(\pi)\in K\backslash G/K$.

\begin{lem}
\begin{enumerate}
\item The assignment $\lambda \lra c_{\lambda}$ induces an isomorphism
of free abelian groups $\Z[P^+] \lra \mathcal{H}$.
\item The Hecke ring is a filtered ring with the filtration levels indexed by
the ordered abelian semigroup $P^+$. In particular, we have
\begin{equation}
\label{heckesum}
c_{\lambda} \cdot c_{\mu} =
c_{\lambda + \mu} + \sum\limits_{\nu<\lambda + \mu}
m_{\lambda,\mu}(\nu) c_{\nu}.
\end{equation}
\end{enumerate}
\end{lem}

Here and in what follows $\sum\limits_{\nu<\mu}$ denotes the sum over elements
$\nu$ from $P^+$. We note that this sum is finite.

\medskip
We will prove in Lemma \ref{polynomialstructureconstants}, that the structure
constants $m_{\lambda,\mu}(\nu)$ are polynomials in $q$ with integer
coefficients.
Here and below we will keep track of the dependence of certain
quantities on the cardinality $q$ of the residue field. Thus $q$ will play
the role of a variable in what follows. One of the main points here is that
the structure constants $n_{\lambda,\mu}(\nu)$ of the representation ring
$R(G^{\vee})$ do not depend on $q$.  They will be encoded in the
coefficients of the polynomials $m_{\lambda,\mu}(\nu)$ (along with
the coefficients of the Kazhdan-Lusztig polynomials $b_{\lambda}(\mu)$).

We recall the definition of the structure constants $m_{\alpha,\beta,\gamma}(\delta)$:
Given double cosets
of $\al(\pi), \be(\pi), \ga(\pi)$ in $K\backslash G/K$, consider the
characteristic functions $c_{\al}$ , $c_{\be}$ and $c_{\ga}$ of these double cosets;
then decompose the triple product
$$
c_{\alpha} \cdot c_{\beta} \cdot c_{\gamma} =
\sum_{\delta} m_{\alpha,\beta,\gamma}(\delta) \  c_{\delta}
$$
in the Hecke algebra $\mathcal{H}$. This defines  $m_{\alpha,\beta,\gamma}(\delta)$.
Our primary interest is $m_{\alpha,\beta,\gamma}(0)$
(where $0$ is the trivial character); it will be viewed as a function of the variable
$q$.

In our proof of the saturation conjecture for $GL(m)$ we will need the following
lemma, where  $\gamma^{*}$ denotes the weight contragredient to $\gamma$. Then $\gamma^{*}(\pi)$
is a representative for the double coset obtained by inverting the elements in
the one represented by $\gamma(\pi)$.

\begin{lem}
\label{contragredientstructureconstant}
$$m_{\alpha,\beta,\gamma}(0) = vol(K \gamma(\pi) K) \
 m_{\alpha,\beta}(\gamma^{*}).$$
 \end{lem}
\proof
We define an inner product $((\ , \ ))$ on $C_0(G)$, the space of
compactly supported complex-valued functions on $G$, by
$$
((f,g)):=\int_G f(x) \overline{g(x)} dx.$$

We first claim that
\begin{equation*}
((c_{\alpha},c_{\beta})) = \begin{cases} vol\ K \alpha(\pi) K, & \text{if $\alpha
=\beta$}\\
0, & \text{otherwise}.
\end{cases}
\end{equation*}
Indeed
$$((c_{\alpha},c_{\beta})) = \int_G c_{\alpha}(x) c_{\beta}(x) dx.$$
The claim follows by noting that the product function is identically
one if the cosets coincide and otherwise it is identically zero.

Next we observe that
$$
((c_{\alpha},c_{\beta})) = c_{\alpha}\cdot c_{\beta^{*}}(1).
$$
Indeed, the right-hand side is equal to
$$\int_G c_{\alpha}(x) c_{\beta^{*}}(x^{-1}) dx =
\int_G c_{\alpha}(x) c_{\beta}(x) dx.$$
Hence
$$ c_{\alpha} \cdot c_{\beta^{*}}(1) = \begin{cases} vol\  K \alpha(\pi) K,
& \text{if $\alpha=\beta$}\\
0, &\text{otherwise.}
\end{cases}
$$
\no Finally
$$
m_{\alpha,\beta,\gamma^{*}}(0) =
\sum\limits_{\delta}m_{\alpha,\beta}(\delta) (c_{\delta}\cdot c_{\gamma^{*}})(1)
=m_{\alpha,\beta}(\gamma) vol\ K \gamma(\pi) K. \qed$$

\subsection{A geometric interpretation of $m_{\alpha,\beta,\gamma}(0)$}
\label{geo=Hecke}

In this section we will prove the following

\begin{thm}
\label{trcount}
$m_{\alpha,\beta,\gamma}(0)$ is the number of triangles in the building $X$
with first vertex $o$ and side-lengths $\alpha,\beta,\gamma$.
\end{thm}

We recall that given $\tau\in \Delta$ we define $S(o,\tau)=\{x\in X: \si(o,x)=\tau\}$, the
``$\Delta$-sphere of radius $\tau$ and center at $o$''.

\begin{lem}
For each $\tau \in \Delta$ the group $K$ acts transitively on
$S(o,\tau)$.
\end{lem}
\proof
Let $x\in S(o,\tau)$. By the properties of the action of $G$ on the building
$X$, presented in the List \ref{list}, there exists $g\in G$ such that
$g\cdot o = o$ and $g\cdot x \in \Delta\subset A$
where $A$ is the model apartment.
Since $g\cdot o = o$ we have $g=k \in K$, whence
$$
\overrightarrow{o~ k\cdot x}=
\overrightarrow{o~ \tau}, \text{(equality of vectors)}.$$
If follows that $k\cdot x=\tau$ and $x=k^{-1}\tau$.
\qed

As a consequence we can identify the right $K$--quotients of the $K$--double cosets in $G$ to the
$\Delta$--spheres $S(o,\tau)$ in the building. The proof of the following lemma is then clear.

\begin{lem}
Let $\alpha \in X_*(\underline{T})$.
Then the map $k\lra k\cdot \alpha(\pi)$ induces a bijection between the quotient
$K\alpha(\pi)K/K$ and $S(o,\alpha(\pi))$.
\end{lem}

We define
$$
f:K\alpha(\pi)K/K \times K\gamma^{*}(\pi)K/K \lra K\backslash G /K$$
by sending $(g_1, g_2)$ to the double coset represented by
$f(g_1,g_2) := g_1^{-1} \cdot g_2$.
The reader will observe that $f$ is well-defined. Thus $f$ induces
a map (again denoted $f$) from $S(o,\alpha(\pi)) \times S(o,\gamma^{*}(\pi))$
to $K\backslash G /K$.

Now let $\alpha,\beta,\gamma \in X_*(\underline{T})$. We define the set
$$
\mathcal{T}_{\alpha,\beta,\gamma} := \{ (x,y)\in K\alpha(\pi)K/K \times
K\gamma^{*}(\pi)K/K : f(x,y) \in K\beta(\pi)K \}.
$$
Note that this set is finite since $q<\infty$.
Let $\Delta(o,x_1,x_2)$ be a triangle in the building $X$ with the $\De$-side lengths
$\alpha, \beta, \gamma$.  Then $f(x_1,x_2) \in K\beta(\pi)K$ and we obtain
a map $F:\Delta(o,x_1,x_2)\mapsto (x_1,x_2)$ from the space of triangles in $X$ with side-lengths
$\alpha, \beta, \gamma$ and first vertex at $o$ into the set
$\mathcal{T}_{\alpha,\beta,\gamma}$. We leave the proof of the following
lemma to the reader.

\begin{lem}
$F$ is a bijection.
\end{lem}

\begin{rem}
The appearance of the contragredient coweight $\gamma^{*}$ comes about
because we require $\si(x_2, o)=\ga$, but  $\si(x_2, o)$ is the contragredient
of the $\si(o, x_2)$.
\end{rem}

Write
$$
K\alpha(\pi)K = \cup_{i=1}^I x_i K  \text{~~~and~~~}  K\gamma^*(\pi)K =
\cup_{j=1}^J y_j K,$$
where both $I$ and $J$ are finite (since $q<\infty$).
The theorem will be a consequence of the following lemma.

\begin{lem}
$m_{\alpha,\beta,\gamma}(0) = \#(\mathcal{T}_{\alpha,\beta,\gamma})$.
\end{lem}
\proof
We have
\begin{align*}
m_{\alpha,\beta,\gamma}(0) =  &\  c_{\alpha}\cdot c_{\beta} \cdot c_{\gamma}(1)
                           =  \int_G(\int_G c_{\alpha}(x) c_{\beta}(x^{-1}y)
                             c_{\gamma}(y^{-1}) dx) dy \\
                           =  & \int_G \left(  \sum_{i=1}^I \int_K
                    c_{\alpha}(x_ik)c_{\beta}(k^{-1}x_i^{-1}y) dk \right)
                    c_{\gamma}(y^{-1}) dy\\
                    =  & \sum_{i=1}^I c_{\alpha}(x_i) \int_G
                    c_{\beta}(x_i^{-1}y) c_{\gamma}(y^{-1}) dy \\
                            = & \sum_{i=1}^I \sum_{j=1}^J
                            c_{\alpha}(x_i)\int_K  c_{\beta}(x_i^{-1}y_jk)
                            c_{\gamma}(k^{-1}y_j^{-1}) dk\\
                            = & \sum_{i=1}^I \sum_{j=1}^J c_{\alpha}(x_i)
                            c_{\beta}(x_i^{-1}y_j)c_{\gamma}(y_j^{-1})\\
                            =  &\sum_{i=1}^I \sum_{j=1}^J c_{\beta}(x_i^{-1}y_j)
                            = \#(\mathcal{T}_{\alpha,\beta,\gamma})   \qed
\end{align*}

As a consequence of Theorem \ref{trcount} and
Lemma \ref{contragredientstructureconstant}, we find that the structure constants
for the Hecke algebra are determined by the geometry of the building.

\begin{thm}\label{geometrydeterminesstructureconstants}
Let $\alpha,\beta,\gamma \in X_*(\underline{T})$.
Then we have
\begin{equation}
\label{struc}
m_{\alpha,\beta}(\gamma)=
\frac{\#(\mathcal{T}_{\alpha,\beta,\gamma^{*}})}{\#(S(o,\gamma(\pi))) }.
\end{equation}
\end{thm}

\proof
We have only to apply Lemma \ref{contragredientstructureconstant} and
observe that since we are assuming that $vol(K) = 1$ we
have $vol(K \gamma(\pi) K) = vol (K \gamma(\pi) K/K) = \#(S(o,\gamma(\pi)))$.
\qed

\subsection{The Satake transform}

In this section we define an integral transform $S$, the Satake transform,
from compactly supported, $K$-biinvariant functions on $G$ to
left $K$-invariant, right $U$-invariant
functions on $G$.

Let $\delta:B \lra \mathbb{R}^*_+$ be the modular function of $B$,
\cite[\S 3]{Gross}. We may regard $\delta$ as a left $K$-invariant, right
$U$-invariant function on $G$. By the Iwasawa decomposition for $G$,
\cite[pg. 51]{Tits}, any such function is determined by its restriction to $T$.
We normalize the Haar measure $du$ on $U$ so that the open subgroup $K \cap U$
has measure $1$.
For a compactly supported $K$-biinvariant function $f$ on $G$ we define its
{\em Satake transform} as a function $Sf(g)$ on $G$ given by
$$
Sf(g)= \delta(g)^{1/2} \cdot \int_U f(gu)du.$$

 Then $Sf$ is a left $K$-invariant, right $U$-invariant
function on $G$ with values in $\Z[q^{1/2}, q^{-1/2}]$; this function
is determined by its restriction to $T/T\cap K \cong X_*(\underline{T})$.
 The main result of \cite{Satake} (see also \cite[pg. 147]{Cartier}), is that the image of
$S$ lies in the subring
$$
(\Z[X_*(\ul{T})])^W \otimes \Z[q^{1/2},q^{-1/2}] \cong
 R(G^{\vee})\otimes \Z[q^{1/2},q^{-1/2}],$$
where $\cong$ is a ring isomorphism.
 Here and below, $G^{\vee}= \ul{G}^{\vee}(\C)$.
Furthermore we have (see \cite{Satake}, \cite[pg. 147]{Cartier}):

\begin{thm}
The Satake transform gives a ring isomorphism
$$
S: \mathcal{H}\otimes \Z[q^{1/2},q^{-1/2}] \cong
R(G^{\vee})\otimes \Z[q^{1/2},q^{-1/2}].$$
 \end{thm}
For $\lambda \in P^+$  let $ch V_{\lambda}$ be the
character of the irreducible representation $V_{\lambda}$ of $G^{\vee}$
(see e.g. \cite[pg. 375]{FultonHarris}).
We may identify $ch V_{\lambda}$ with a
$W$-invariant element of $\Z[X^*(\underline{T}^{\vee})] =
\Z[X_*(\underline{T})]$.

\medskip
In what follows we will need two bases for the free
$\Z[q^{1/2},q^{-1/2}]$-module $R(G^{\vee})\otimes \Z[q^{1/2},q^{-1/2}]$.
The first basis is
$$
\mathcal{S}:=\{S(c_{\lambda}): \lambda \in P^+\},$$
the second basis is
$$
\mathcal{R}:= \{ch V_{\lambda}: \lambda \in P^+\}.$$

The change of basis matrices relating these two bases are both upper triangular
(see Lemma \ref{expansionlemma} below) with entries in the ring
$\Z[q^{1/2},q^{-1/2}]$. We define the
order (at $\infty$) $ord(f)$ of the element $f=\sum _{k=-M}^{k=N}a_k q^{k/2}
\in \Z[q^{1/2},q^{-1/2}]$ by $ord(f) =N$, provided that $a_N \neq 0$.
Note that if $f$ is a polynomial
in $q$ then $ord(f) = 2 deg(f)$ where $deg(f)$ is the degree of $f$ in $q$.
We will accordingly extend the degree to $\Z[q^{1/2},q^{-1/2}]$ by defining
$$
deg(f) = 1/2 \ ord(f).$$
Thus the extended degree takes values in the half integers.
We define $deg_{\mathcal{S}}(F)$ for $F\in R(G^{\vee})\otimes \Z[q^{1/2},q^{-1/2}]$
to be the maximum of the degrees of the components of $F$ {\em when $F$
is expressed in the basis $\mathcal{S}$}. We will use a similar convention
(i.e., expanding it terms of the basis $\mathcal{S}$)
when we speak of the ``leading term" of $F$.
We retain the notation $deg$ (without subscript) for the ordinary notion of degree
$\in \frac{1}{2}\Z$ for the Laurent polynomials in $\Z[q^{1/2}, q^{-1/2}]$.
Analogously, we define $deg_{\mathcal{R}}(F)$ when $F$ is expanded in terms of the basis
$\mathcal{R}$.

In what follows we will need two formal properties of $deg_{\mathcal{R}}$
and $deg_{\mathcal{S}}$. We leave their proofs to the reader.
The first is the {\em ultrametric} inequality
which we state for $deg_{\mathcal{S}}$
$$
deg_{\mathcal{S}}(\sum_{i=1}^n F_i) \leq \max_{1 \leq i \leq n}(
deg_{\mathcal{S}}(F_1),\cdots, deg_{\mathcal{S}}(F_n)).$$

The second is the obvious identity (again stated for $deg_{\mathcal{S}}$)
$$deg_{\mathcal{S}}( p \cdot F) = deg(p) + deg_{\mathcal{S}}(F), \ p \in
\Z[q^{1/2},q^{-1/2}].$$

For the following lemma, due to Lusztig, we refer the reader to \cite[(3.11) and (3.12)]{Gross},
\cite[(6.10)]{Lusztig2}. See \cite[\S 2]{Haines},
for the statement of the lemma in the generality we require here.

\begin{lem}
\label{expansionlemma}
1. There exist polynomials $a_{\lambda}(\mu)$ in $q$ such that
$$
S(c_{\lambda})= q^{\<\lambda,\rho\>} ch V_{\lambda} + \sum\limits_{\mu < \lambda}
a_{\lambda}(\mu) q^{\<\mu,\rho\>} ch V_{\mu}.$$

2. Conversely, there exist polynomials $b_{\lambda}(\mu)$ in $q$ such that
$$q^{\<\lambda,\rho\>} ch V_{\lambda} = S(c_{\lambda}) +
\sum\limits_{\mu < \lambda} b_{\lambda}(\mu) S(c_{\mu}).$$
\end{lem}

\medskip
The degree estimate below follows from \cite{Lusztig2},
who proved that the polynomials $b_{\lambda}(\mu)$
are {\em Kazhdan-Lusztig polynomials}.
These polynomials have many remarkable properties but the only property
we need here is the degree estimate. The inequality  of the next lemma will
play a critical role in our proofs.

\begin{lem}
\label{KazhdanLusztigestimate}
For all $\mu<\la, \mu\in P^+$ we have:
$$
deg(b_{\lambda}(\mu)) < \ \<\lambda - \mu,\rho\> \leq \  \<\lambda,\rho\>.$$
\end{lem}

This lemma when combined with Lemma \ref{expansionlemma} has the following
consequences.

\begin{lem}\label{characterestimate}

\hfill

\begin{enumerate}
\item $deg_{{\mathcal S}}( ch V_{\lambda}) <  0, \ \text{for all} \
\lambda \in P^+\setminus \{0\}$ .
\item $deg_{\mathcal{S}} \leq deg_{\mathcal{R}}.$
\end{enumerate}
\end{lem}

\proof
By Lemma \ref{expansionlemma} we get
$$
 ch V_{\lambda} = q^{-\<\lambda,\rho\>} S(c_{\lambda}) +
\sum\limits_{\mu < \lambda} q^{-\<\lambda,\rho\>} b_{\lambda}(\mu) S(c_{\mu}).$$
\noindent Hence Lemma \ref{KazhdanLusztigestimate} implies
$$
deg_{{\mathcal S}}(ch V_{\lambda}) < 0.$$

Since $deg_{\mathcal{R}}(ch V_{\lambda}) = 0$ we find that the inequality
of the second statement of the lemma holds on a basis. Consequently it
holds for all $F \in \mathcal{H}\otimes \Z[q^{1/2},q^{-1/2}] \cong
R(G^{\vee})\otimes \Z[q^{1/2},q^{-1/2}]$.
Indeed, writing $F = \sum_{i=1}^n p_i \ ch V_i$ we obtain, using
the ultrametric inequality and the first statement of the lemma,
$$ deg_{\mathcal{S}}( F) \leq \max_{1 \leq i \leq n}
(deg_{\mathcal{S}}(p_1 ch V_1), \cdots, deg_{\mathcal{S}}(p_n ch V_n))$$
$$
\leq (deg(p_1), \cdots, deg(p_n)) = deg_{\mathcal{R}}( F).$$
\qed

Recall that $0$ denote the trivial cocharacter.

\begin{lem}
\label{polynomialstructureconstants}
There exists a polynomial $M_{\alpha,\beta,\gamma}(q)$ in the variable $q$,
such that we have the equality of functions in $q$:
$$
m_{\alpha,\beta,\gamma}(0) = M_{\alpha,\beta,\gamma}(q).$$
Furthermore, all structure constants $m_{\alpha,\beta}(\gamma)$ are polynomials
in $q$.
\end{lem}
\proof
We will prove that all the structure constants
$m_{\alpha,\beta}(\gamma)$ are polynomials in $q$. The first statement will follow from this.

Following \cite{Gross} we define $\phi_{\lambda} \in
R(G^{\vee})\otimes \Z[q^{1/2},q^{-1/2}]$ by
$$
\phi_{\lambda} = q^{\<\lambda,\rho\>} ch V_{\lambda}.
$$
The elements $\{\phi_{\lambda}: \lambda \in X^*(\ul{T}^{\vee})\}$
give a new  basis for the $\Z[q^{1/2},q^{-1/2}]$-module
$R(G^{\vee})\otimes \Z[q^{1/2},q^{-1/2}]$.

We claim that the structure
constants for the ring relative to this basis are polynomials in $q$.
Indeed, we have
\begin{align*}
\phi_{\alpha} \cdot \phi_{\beta} = & q^{\<\alpha + \beta, \rho \>}
ch V_{\alpha}\cdot ch V_{\beta} = q^{\<\alpha + \beta, \rho \>}
\sum\limits_{\gamma \leq \alpha + \beta} n_{\alpha + \beta}(\gamma)\
ch V_{\gamma} \\
= & \sum\limits_{\gamma \leq \alpha + \beta} n_{\alpha + \beta}(\gamma)\
q^{\<\alpha + \beta - \gamma,\rho\>} \phi_{\gamma}.
\end{align*}
Now, since $\gamma \leq \alpha + \beta$, the coweight $\alpha+\beta -\gamma$
is a sum of positive coroots and consequently
$\<\alpha + \beta - \gamma,\rho\>$ is an nonnegative integer.

We can now prove the second statement. Let $\alpha$ and $\beta$ be given.
We expand $S(c_{\alpha})$ and $S(c_{\beta})$ in terms of the
basis of $\phi_{\lambda}$'s thereby introducing
the polynomials $a_{\lambda}(\mu)$. We then multiply the resulting
expressions. According to the paragraph above, the result is an expression
in the $\phi_{\lambda}$'s with polynomial coefficients in $q$. We then
substitute for the $\phi_{\lambda}$'s using Lemma \ref{expansionlemma}
introducing the polynomials $b_{\lambda}(\mu)$. \qed

\subsubsection{A remarkable nonvanishing property of the polynomials
$M_{\alpha,\beta,\gamma}$}
It follows from Theorem \ref{trcount} that if $q$ is a prime power then
the values $M_{\alpha,\beta,\gamma}(q)$ are nonnegative. We now show that if
for some prime power $q$ the value $M_{\alpha,\beta,\gamma}(q)$ is nonzero
then all values at prime powers are nonzero. We note that the examples
we computed below show that the polynomial $M_{\alpha,\beta,\gamma}$
can have coefficients of both signs.

\begin{thm}\label{nonvanishing}
If a polynomial $M_{\alpha,\beta,\gamma}$ is nonzero for some prime
power $q = p^e$ then it is nonzero at all prime powers.
\end{thm}
\proof
We assume that $q$ is given such that
$$
M_{\al,\be,\ga}(q)\ne 0.
$$
This means that
$$
(\al,\be,\ga)\in Sol({\bf Q3}, \ul{G}(\K_q)),
$$
where $\K_q$ is a local field with residue field of order $q$. In other words,
in the Euclidean building $X_q$ corresponding to the group $\ul{G}(\K_q)$
there exists a triangle $\tau$ with $\De$-side lengths
$\al, \be, \ga$, whose vertices are special vertices of $X$.  It then follows
from Part 5 of Theorem \ref{mainthm} or, equivalently,
from the Transfer Theorem \ref{transfer}, that
$$
(\al,\be,\ga)\in Sol({\bf Q3}, \ul{G}(\K)),
$$
for an {\em arbitrary} local field $\K$ (with any order $q'\ge 2$ of the residue field).
Therefore  $M_{\al,\be,\ga}(q')\ne 0$ for all prime powers $q'$.  \qed

\subsection{A solution of Problem {\bf Q4} is a solution of Problem {\bf Q3}}

We recall the definition of the structure constants $n_{\alpha,\beta,\gamma}(\delta)$:
$$
ch(V_{\alpha})\cdot ch(V_{\beta})\cdot ch(V_{\gamma}) =
\sum_{\delta\in P^+} n_{\alpha,\beta,\gamma}(\delta) \  ch(V_{\delta}).$$
We are interested in comparing the coefficient
$n_{\alpha,\beta,\gamma}(0)$ (corresponding to the trivial character $0$
of $G^{\vee}$) with $m_{\alpha,\beta,\gamma}(0)$.

The main goal of this section is to prove the following

\begin{thm}
\label{q4->q3}
If the triple $\alpha,\beta,\gamma$ is a solution of Problem {\bf
Q4} then it is also a solution of Problem {\bf Q3}. More
precisely,
$$
n_{\alpha,\beta,\gamma}(0) \neq 0 \Lra m_{\alpha,\beta,\gamma}(0) \neq 0.$$
\end{thm}

We first prove the following:

\begin{thm}
\label{4->3}
Suppose that $\underline{G}$ is of a reductive algebraic group over $\K$
which is split over $\K$. Then:

(a) The degree of the polynomial $M_{\alpha,\beta,\gamma}(q)$ is at
most $\<\alpha + \beta + \gamma, \rho \>$.

(b) The leading coefficient of $M_{\alpha,\beta,\gamma}(q)$, i.e., the coefficient
at $q^{\<\alpha + \beta + \gamma ,\rho \>}$, is equal to $n_{\alpha,\beta,\gamma}(0)$.
\end{thm}

{\em Proof of Theorem \ref{4->3}}. 
For the benefit of the reader, below we explain the idea behind the proof of
Theorem \ref{4->3}.
We have to compute the coefficient of $0$ in the Hecke triple product
$c_{\alpha}\cdot c_{\beta} \cdot c_{\gamma}$. Instead, we will compute
the coefficient of $0$ in the triple product
$S(c_{\alpha})\cdot S(c_{\beta}) \cdot S(c_{\gamma})$ in the representation
ring $R(G^{\vee})\otimes \Z[q^{1/2},q^{-1/2}]$ where the triple product is
expanded {\em relative to the basis $\mathcal{S}$}.
In order to do this we will use the formula from
Lemma \ref{expansionlemma} and compute (initially) in the basis $\mathcal{R}$.
However, to prove the theorem we must again apply
Lemma \ref{expansionlemma} to rewrite the result in terms of the
basis $\mathcal{S}$.

The theorem will follow from the next three lemmas.
We owe the first of these to Jiu-Kang Yu.

\begin{lem}
\label{basic}
For all $\mu\in P^+$ satisfying $\mu<\la$ we have:
$$
deg (a_{\lambda}(\mu) q^{\<\mu,\rho\>}) < \<\lambda,\rho\>.$$
\end{lem}
\proof
The proof is by induction on $\lambda$ with respect to the partial order $<$.
We remind the reader that  each  set $\{\mu\in P^+ : \mu < \lambda\}$ is finite.
If the set  $\{\mu\in P^+ : \mu < \lambda\}$ is empty then there is nothing
to prove.

Now assume that we are given $\lambda \in P^+$ and that we have
proved the above estimate for all predecessors of $\lambda$ in the partial
order $<$. We have, according to Lemma \ref{expansionlemma},
$$
S(c_{\lambda}) = q^{\<\lambda,\rho\>} ch V_{\lambda} -
\sum\limits_{\mu < \lambda} b_{\lambda}(\mu) S(c_{\mu})
$$
and
$$
S(c_{\mu})= q^{\<\mu,\rho\>} ch V_{\mu} +
\sum\limits_{\eta < \mu} a_{\mu}(\eta) q^{\<\eta,\rho\>} chV_{\eta}.$$
Thus
$$
S(c_{\lambda}) = q^{\<\lambda,\rho\>} ch V_{\lambda} -
\sum\limits_{\mu<\lambda} b_{\lambda}(\mu) q^{\<\mu,\rho\>} ch V_{\mu} -
\sum\limits_{\mu < \lambda}\sum\limits_{\eta <
\mu} b_{\lambda}(\mu) a_{\mu}(\eta) q^{\<\eta,\rho\>}ch V_{\eta}.$$
We make the change of variable $\mu \to \tau, \eta \to \mu$ in the last sum
to obtain
$$
S(c_{\lambda}) = q^{\<\lambda,\rho\>} ch V_{\lambda} -
\sum\limits_{\mu<\lambda} b_{\lambda}(\mu) q^{\<\mu,\rho\>} ch V_{\mu} -
\sum\limits_{\mu < \lambda}\sum\limits_{\{\tau : \mu <
\tau < \lambda\}} b_{\lambda}(\tau) a_{\tau}(\mu) q^{\<\mu,\rho\>}ch V_{\mu}.$$

Since $ch V_{\mu},\mu \in P^+,$ is a basis of the $\Z[q^{-1/2},q^{1/2}]$-module
$R(G^{\vee})\otimes \Z[q^{1/2},q^{-1/2}]$,
by combining the previous formula with the part 1 of Lemma
\ref{expansionlemma}, we obtain:
$$
a_{\lambda}(\mu) q^{\<\mu,\rho\>} = - b_{\lambda}(\mu) q^{\<\mu,\rho\>}
-\sum\limits_{\{\tau : \mu <
\tau < \lambda\}} b_{\lambda}(\tau) a_{\tau}(\mu) q^{\<\mu,\rho\>}.$$

The degree of the first term satisfies the required estimate by Lemma
\ref{KazhdanLusztigestimate}. We estimate the degrees of the terms in the
second sum.
By the induction hypothesis, $deg(a_{\tau}(\mu) q^{\<\mu,\rho\>}) < \<\tau,\rho\>$
and by Lemma \ref{KazhdanLusztigestimate}  we have
$$
deg(b_{\lambda}(\tau)) < \<\lambda - \tau, \rho\>.$$
Consequently, for all $\mu$, we get
$$
 deg(b_{\lambda}(\tau) a_{\tau}(\mu) q^{\<\mu,\rho\>}) <
\<\lambda - \tau,\rho\> + \<\tau,\rho\> =\<\lambda,\rho\>. \qed
$$

\begin{cor}
$deg(a_{\lambda}(\mu)) < \<\lambda - \mu,\rho\>$.
\end{cor}

%

\begin{lem}
\label{triple}
For $\al, \be, \ga\in P^+$ such that $\al+\be+\ga$ is a nontrivial character,
we have:

\begin{enumerate}
\item $
deg_{{\mathcal R}}
(S(c_{\alpha})\cdot S(c_{\beta}) \cdot S(c_{\gamma})
- q^{\<\alpha + \beta + \gamma,\rho\>} ch V_{\alpha} \cdot ch V_{\beta}
\cdot ch V_{\gamma}) < \<\alpha + \beta + \gamma, \rho\> $.
\item $ deg_{{\mathcal R}}
(S(c_{\alpha})\cdot S(c_{\beta}) \cdot S(c_{\gamma}) =
\<\alpha + \beta + \gamma, \rho\>. $
\end{enumerate}
\end{lem}
\proof By expanding the triple product
$$
S(c_{\alpha})\cdot S(c_{\beta}) \cdot S(c_{\gamma})=
$$
\begin{align*}
(q^{\<\ga,\rho\>}ch V_\ga + \sum\limits_{\lambda<\gamma} a_{\gamma}(\lambda) q^{\<\lambda,\rho\>}
ch V_\la) \\
\cdot
(q^{\<\be,\rho\>}ch V_\be + \sum\limits_{\mu<\be} a_{\be}(\mu) q^{\<\mu,\rho\>}
ch V_\mu) \\
\cdot
(q^{\<\al,\rho\>}ch V_\al + \sum\limits_{\nu<\al} a_{\al}(\nu) q^{\<\nu,\rho\>}
ch V_\nu)
\end{align*}
in terms of the  basis ${\mathcal R}:=\{ch V_{\tau}\}$ we get the following types of summands:
$$
q^{\<\alpha + \beta + \gamma,\rho\>} ch V_{\alpha} \cdot ch V_{\beta}
\cdot ch V_{\gamma},
$$
$$
q^{\<\alpha + \beta , \rho\>}ch V_{\alpha}
\cdot ch V_{\beta}
\sum\limits_{\lambda<\gamma} a_{\gamma}(\lambda) q^{\<\lambda,\rho\>}
ch V_{\lambda},
$$
$$
q^{\<\ga,\rho\>} ch V_\ga
\cdot
(\sum\limits_{\mu<\be} a_{\be}(\mu) q^{\<\mu,\rho\>}
ch V_\mu) \\
\cdot
( \sum\limits_{\nu<\al} a_{\al}(\nu) q^{\<\nu,\rho\>}
ch V_\nu),
$$
and similar ones obtained by permuting $\al, \be, \ga$, and finally:
$$
a_\al(\nu) a_\be(\mu) a_\ga(\la) q^{\<\nu+\mu+\la,\rho\>}
ch V_{\alpha} \cdot ch V_{\beta} \cdot ch V_{\gamma},
$$
where in the latter case $\nu<\al, \mu<\be, \la<\ga$, $\nu, \mu, \la\in P^+$.
The first term ($q^{\<\alpha + \beta + \gamma,\rho\>} ch V_{\alpha} \cdot ch V_{\beta}
\cdot ch V_{\gamma}$)
cancels out, we estimate the degree of each of the remaining terms separately. We will
do it in the case of the 2-nd and 4-th term and leave the 3-rd term to the reader:

First, since $\la<\ga$,
\begin{align*}
deg_{{\mathcal R}} [q^{ \<\alpha + \beta , \rho\> }ch V_{\alpha}
\cdot ch V_{\beta} \cdot a_{\gamma}(\lambda) q^{\<\lambda,\rho\>}
ch V_{\lambda}] \le
\\
 \<\alpha + \beta, \rho\> + deg( a_{\gamma}(\lambda)q^{\<\la,\rho\>} )<\\
\hbox{~(by Lemma \ref{basic})~~}\\
< \<\alpha + \beta, \rho\> +  \<\ga, \rho\> =\<\al+\be+\ga, \rho\>.
\end{align*}
Similarly,
\begin{align*}
deg_{{\mathcal R}}[a_\al(\nu) a_\be(\mu) a_\ga(\la) q^{\<\nu+\mu+\la,\rho\>}
ch V_{\alpha} \cdot ch V_{\beta} \cdot ch V_{\gamma}] = \\
deg [a_\al(\nu) q^{\<\nu,\rho\> } ]+ deg [a_\be(\mu)  q^{\<\mu,\rho\> } ] +
deg [a_\ga(\la)  q^{\<\la,\rho\> } ] <
\\
\hbox{~(by Lemma \ref{basic})~~}\\
 < \<\al+\be+\ga, \rho\>.
\end{align*}
The first statement immediately implies the second and the lemma follows. \qed

\begin{cor}
\label{Striple}
Under the assumptions of the above lemma, we have:
$$
deg_{{\mathcal S}}
(S(c_{\alpha})\cdot S(c_{\beta}) \cdot S(c_{\gamma})
- q^{\<\alpha + \beta + \gamma,\rho\>} ch V_{\alpha} \cdot ch V_{\beta}
\cdot ch V_{\gamma}) < \<\alpha + \beta + \gamma, \rho\> .$$
\end{cor}
\proof The assertion follows from Lemma \ref{triple} since
by Lemma \ref{characterestimate} we have $deg_{\mathcal{S}}
\le deg_{\mathcal{R}}$. \qed

\medskip

Now we can prove the degree estimate in Theorem \ref{4->3}. Observe that since
$M_{\al,\be,\ga}(q)$ is one of the coefficients (the coefficient of the
trivial double coset) of the triple product
$c_\al\cdot c_\be \cdot c_\ga$ when expressed in terms of the double
coset basis for the Hecke ring, we have
$$deg(M_{\al,\be,\ga}(q)) \leq deg_{\mathcal{S}}(c_\al\cdot c_\be \cdot c_\ga)$$
and hence
$$
deg(M_{\al,\be,\ga}(q)) \leq deg_{\mathcal{S}} (S(c_\al)\cdot S(c_\be)\cdot S(c_\ga))\le
deg_{\mathcal{R}} (S(c_\al)\cdot S(c_\be)\cdot S(c_\ga)) =
\<\alpha + \beta + \gamma,\rho\>.
$$
Here the last equality follows from Lemma \ref{triple}. This proves
the first assertion of Theorem \ref{4->3}.

We now identify the leading term of $M_{\al,\be,\ga}(q)$.

Let $N_{\alpha,\beta,\gamma}$ be the element of the Hecke ring defined
by
$$
N_{\alpha,\beta,\gamma}:=q^{\<\alpha + \beta + \gamma,\rho\>} ch V_{\alpha} \cdot ch V_{\beta}
\cdot ch V_{\gamma},$$

Let $\mathbb{I}$ denote the identity element of the representation ring
(i.e. the character of the trivial representation).
We identify this element (via the Satake isomorphism) with the trivial
double coset $c_0$.
We have
\begin{lem}
$$
deg_{{\mathcal S}}(N_{\alpha,\beta,\gamma}
- n_{\alpha,\beta,\gamma}(0) q^{\<\alpha + \beta + \gamma,\rho\>}\mathbb{I}) <
\<\alpha + \beta + \gamma, \rho\> .$$
\end{lem}
\proof
By definition of the structure constants $n_{\alpha,\beta,\gamma}(\delta)$
we have
$$
q^{\<\alpha + \beta + \gamma,\rho\>}ch V_{\alpha} \cdot
ch V_{\beta}\cdot ch V_{\gamma}
=\sum\limits_{\delta \leq \alpha + \beta + \gamma} n_{\alpha,\beta,\gamma}(\delta)
\  q^{\<\alpha + \beta + \gamma - \delta,\rho\>} q^{\<\delta,\rho\>}ch V_{\delta}.$$
But by Lemma \ref{characterestimate}, $\delta\ne 0 \Lra
deg_{{\mathcal S}}(q^{\<\delta,\rho\>} ch V_{\delta}) < \<\delta,\rho\>$.
Noting that $\mathbb{I}$ is the character of the trivial representation
the lemma follows.\qed

We now combine the previous lemma, the ultrametric inequality and
Corollary \ref{Striple}
to obtain the estimate

\begin{align*}
&deg_{\mathcal{S}}(S(c_{\alpha})\cdot S(c_{\beta})
\cdot S(c_{\gamma}))
- n_{\alpha,\beta,\gamma}(0) q^{\<\alpha + \beta + \gamma,\rho\>}\mathbb{I})
\\
&\leq
\max( deg_{\mathcal{S}}(S(c_{\alpha})\cdot S(c_{\beta}) \cdot S(c_{\gamma})
-N_{\alpha,\beta,\gamma}), \  deg_{\mathcal{S}}(N_{\alpha,\beta,\gamma} -
n_{\alpha,\beta,\gamma}(0)
q^{\<\alpha + \beta + \gamma,\rho\>}\mathbb{I})) \\
&< \<\alpha + \beta + \gamma, \rho\> .
\end{align*}

Now the previous estimate on the difference of elements in the Hecke
ring implies that degrees of {\em all} coefficient polynomials of the difference
when expanded in terms of the double coset basis satisfy the
same estimate. Now the coefficient of the identity element $c_0 = \mathbb{I}$
of the triple product is $m_{\alpha,\beta,\gamma}(0) = M_{\alpha,\beta,\gamma}(q)$.
Thus applying the above observation to the identity coefficient of the
difference we obtain

$$deg(M_{\alpha,\beta,\gamma}(q) - n_{\alpha,\beta,\gamma}(0)
q^{\<\alpha + \beta + \gamma,\rho\>}) < \<\alpha + \beta + \gamma, \rho\>.$$

This proves Theorem \ref{4->3} and in particular
implies that the polynomial $M_{\alpha,\beta,\gamma}$ is nonzero.
Hence for nonarchimedean local fields with sufficiently large residue fields
(i.e. sufficiently large $q$)
the structure constant $m_{\alpha,\beta,\gamma}(0)$ is nonzero.
But by Theorem \ref{nonvanishing} if the structure constant
$m_{\alpha,\beta,\gamma}(0)$ is nonzero for some nonarchimedean local field
$\mathbb{K}$ it is nonzero for all nonarchimedean local fields.
Hence Theorem \ref{q4->q3} is proved.

The above theorem also immediately implies  the following

\begin{cor}
If $n_{\alpha,\beta,\gamma}(0)\ne 0$ then
$\<\alpha + \beta + \gamma ,\rho \>\in \Z$.
\end{cor}

\medskip
We now give another explanation why
$$
n_{\alpha,\beta,\gamma}(0) \neq 0 \Lra
\<\alpha + \beta + \gamma,\rho\>\in \Z.$$

\begin{thm}\label{sum}
\begin{enumerate}
\item $m_{\alpha,\beta,\gamma}(0)\neq 0 \Lra \alpha+\beta+\gamma \in
Q(R^{\vee})\Lra \<\alpha + \beta + \gamma,\rho\>\in \Z$.

\item $n_{\alpha,\beta,\gamma}(0)\neq 0 \Lra \alpha+\beta+\gamma \in
Q(R^{\vee})$.
\end{enumerate}
\end{thm}
\proof (1) follows from Theorem \ref{fixedbarycenter}, part 4. Hence (2)
follows from (1) via Theorem \ref{4->3}, however we give a direct proof below.

The center of $G^{\vee}$ acts on each of the factors of the
triple tensor product by a scalar, hence it acts by a scalar on
the triple tensor product itself. But if $n_{\alpha,\beta,\gamma}(0) \neq 0$
this scalar is necessarily $1$. Hence the center fixes the (highest) weight
vector with weight $\alpha + \beta + \gamma$, whence $\alpha + \beta + \gamma$
annihilates the center. Accordingly it is in the character lattice of the
adjoint group $Ad\  G^{\vee}$. Thus $\alpha + \beta + \gamma \in Q(R^{\vee})$. \qed

\subsection{A solution of Problem {\bf Q3} is not necessarily a solution of
Problem {\bf Q4}}
\label{3not4}

In this section we consider two examples: $\underline{G} = SO(5)$,
whence $G^{\vee} = Sp(4,\C)$ and the group $\underline{G}$ of type $G_2$
and $G^{\vee} = G_2(\C)$.

We begin with $\underline{G} = SO(5)$; we assume that we have chosen Witt bases for $\C^5$
and $\C^4$. We let $T^{\vee}$ be the  split torus consisting of
diagonal matrices in $Sp(4, \C)$ and $\ul{T}$ be the split torus of $SO(5)$ consisting
of diagonal matrices. We use the
rectangular coordinates $x_1, x_2$ in
the Cartan subalgebra $\mathfrak{a}^{\vee} \subset \mathfrak{sp}(4,\C)$ such that the
simple roots are $x_1 - x_2, 2x_2$. We  set
$$
\alpha = \beta = \gamma = (1,1) \in X^*(\ul{T}^{\vee}) =
X_*(\underline{T}).$$
Note that $\al+\be+\ga$ belongs to the coroot lattice of $SO(5)$,
i.e., the condition stated in Theorem \ref{sum} is satisfied.

We give two proofs of the next lemma, the first computational, the second
conceptual.

\begin{lem}
\begin{align*}
& n_{\alpha,\beta,\gamma}(0) = 0.\\
& m_{\alpha,\beta,\gamma}(0) = q^5 -q.
\end{align*}
Thus this choice of side-lengths is a solution Problem {\bf Q3}
but not of Problem {\bf Q4}.
\end{lem}

\begin{rem}
In fact the second equality implies the first because, in this case,
$\<\alpha + \beta + \gamma,\rho\> = 6$. Hence, (according to Theorem \ref{4->3})
$M_{\alpha,\beta,\gamma}(q) =
n_{\alpha,\beta,\gamma}(0) q^6 + \text{lower terms}$. But the coefficient
of $q^6$ in the above formula for $m_{\alpha,\beta,\gamma}(0)$ is $0$.
\end{rem}
\proof
The first equality is obvious if we observe that under the isomorphism
between $Sp(4)$ and $Spin(5)$ the representation of $Sp(4, \C)$ with the highest
weight $(1,1)$ corresponds to the standard (vector) representation of $SO(5, \C)$
on $\C^5$ (considered as a nonfaithful representation of $Spin(5, \C)$).
Now it is standard that the tensor square of this representation consists
of three irreducible summands, the trivial representation (corresponding to
the invariant quadratic form), the exterior square of the standard
representation and the harmonic quadratic polynomials on $\C^5$. Since none
of these are equivalent to the standard representation and all representations
of $Sp(4, \C)$ are self-dual the first statement follows.

The proof of the second statement requires more work. From \cite[pg. 231]{Gross},
we have:
$$
S(c_{(1,1)}) = q^2 ch V_{(1,1)} - 1.$$
Upon taking the cube and calculating in the representation ring of $Sp(4, \C)$
one obtains
\begin{align*}
S(c_{(1,1)})^3 = q^6 [ch V_{(3,3)} + 2 ch V_{(3,1)} + ch V_{(2,0)}+ 3 ch V_{(1,1)}] \\
- 3 q^4 [ch V_{(2,2)} + ch V_{(2,0)} + 1]
+ 3 q^2 ch V_{(1,1)} - 1.
\end{align*}
Upon substituting for the $ch V_{\lambda}$ using Lemma \ref{expansionlemma}
(and computing the appropriate Kazhdan-Lusztig polynomials $b_{\lambda}(\mu)$)
the lemma follows.
\qed

As for the conceptual proof of (the second part of) the lemma, note that
the side-lengths belong to the root-subgroup $SL(2)$ corresponding to the positive
root $x_1 + x_2$. But it is evident that one can construct an equilateral
triangle with side-lengths equal to the positive root
and vertices of the correct type in the tree for $SL(2, \K)$.
Since this tree is convex in the building for $SO(5, \K)$
and the fixed vertex for $SL(2,\mathcal{O})$ in this embedded tree is $o$
(the fixed vertex for $SO(5,\mathcal{O})$), the lemma
follows. Equivalently if we can realize the trivial double coset as
a product of the three $SL(2,\mathcal{O})$--double cosets in $SL(2, \K)$
belonging to $(1,1),(1,1),(1,1)$ then we can certainly do it in $SO(5,\K)$.

\begin{rem}
This last sentence is the essence of the counterexample. The solutions to
Problems {\bf Q3} and {\bf Q4} behave differently with respect to inclusions of subgroups.
The solutions of Problem {\bf Q3} ``push forward", the solutions of Problem {\bf Q4} do not.
\end{rem}

We now give another example of a triple of coroots $\alpha, \beta,
\gamma$ such that the triple is a solution to Problem {\bf Q3} but
not of Problem {\bf Q4}. This example was motivated by unpublished observations of
S.\ Kumar and J.\ Stembridge.

We take $\ul{G} = G_2$. Hence $G^{\vee}
=G_2$ as well. We take $\alpha=\lambda_1$ and $\beta =
\gamma=\lambda_2$. Here $\lambda_1$ is the first fundamental
weight (the highest weight of the unique irreducible seven
dimensional representation) and $\lambda_2$ is the second
fundamental weight (the highest weight of the adjoint
representation).

\begin{lem}
$m_{\lambda_1,\lambda_2,\lambda_2}(0) = q^5(q+1)(q^6-1).$
\end{lem}

\proof
\ From \cite[pg. 231, (5.7)]{Gross}, we have
\begin{align*}
S(c_{\lambda_1}) = & q^3 chV_{\lambda_1} -1\\
S(c_{\lambda_2}) = & q^5 chV_{\lambda_2} -q^3 chV_{\lambda_2} - q^4.
\end{align*}

Multiplying (using LiE) one obtains
$$
S(c_{\lambda_1}) \cdot S(c_{\lambda_2}) =
q^8 chV_{\lambda_1+\lambda_2} + (q^8 - q^6) (chV_{2\lambda_1}+
chV_{\lambda_1}) -q^6 chV_{\lambda_2} -q^4 S(c_{\lambda_1})
-S(c_{\lambda_2}).$$

Hence, by Lemma \ref{expansionlemma} we have
$$
m_{\lambda_1,\lambda_2}(\lambda_2) = b_{\lambda_2,\lambda_1 + \lambda_2}(q)
+(q^2-1) b_{\lambda_2,2\lambda_1}(q) - q - 1.$$

Using $b_{\lambda_2,\lambda_1 + \lambda_2}(q) = 1+q$ and
$b_{\lambda_2,2\lambda_1}(q)=1$ we obtain
$$
m_{\lambda_1,\lambda_2}(\lambda_2)= q^2 -1.$$

Now Lemma \ref{contragredientstructureconstant} implies that
$$
m_{\lambda_1,\lambda_2,\lambda_2}(0)=m_{\lambda_1,\lambda_2}(\lambda_2)
\cdot vol(K\lambda_2(\pi)K).
$$
\ From \cite[pg. 735]{Gross}, we get
$vol(K\lambda_2(\pi)K) = deg(c_{\lambda_2})
=q^{10} + q^9 + q^8 +q^7 + q^6 + q^5.$
The lemma follows. \qed

\begin{cor}
$$n_{\lambda_1,\lambda_2,\lambda_2}(0)=0.$$
\end{cor}
\proof
\ From \cite[pg. 231, (5.7)]{Gross},  we have
$\<\lambda_1,\rho\> = 3$ and $\<\lambda_2,\rho\>=5$
whence $\<\lambda_1 + \lambda_2 + \lambda_2,\rho\> = 13$. Hence, the
structure constant has degree (one) less than the maximum possible degree.
\qed

\begin{figure}[h]
\centerline{\epsfxsize=3.5in \epsfbox{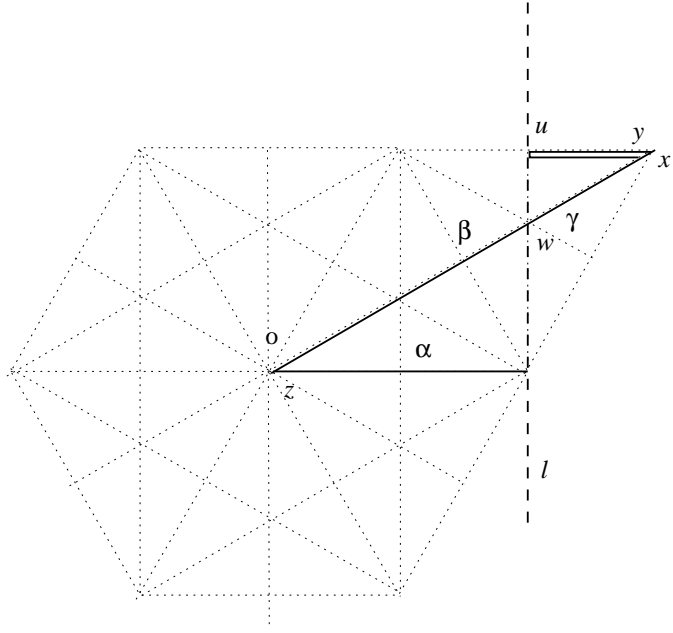}}
\caption{\sl A billiard triangle in the
affine Coxeter complex of type $G_2$.}
\label{g2example}
\end{figure}

We now give a geometric proof (using unfolding) that
$m_{\lambda_1,\lambda_2,\lambda_2}(0)$ is
nonzero.
We show that there exists a triangle $T'$ (with $\De$-side-lengths
$\lambda_1, \lambda_2, \lambda_2$)
in a discrete Euclidean building $X$ modeled on the discrete affine Coxeter complex of type $G_2$,
so that all vertices of $T'$ are special. We start with the billiard triangle
$T$ in Figure \ref{g2example}, contained in a model apartment $A\subset X$;
this triangle  has two geodesic sides $\ol{zx}, \ol{zy}$ and one broken side
$\ol{xuy}$. Note that all three vertices of this billiard triangle are special vertices of the
Coxeter complex. The break point $u$ on the broken side of the billiard
triangle belongs to the wall $l$ of the Coxeter complex, the wall $l$
separates the broken side $\ol{xuy}$ from the vertex $z=o$.
Thus, by Lemma \ref{unfold}, one can unfold the billiard triangle $T$ to a geodesic triangle
$T'\subset X$ preserving the refined side-lengths. One can also easily see that $n_{\lambda_1,\lambda_2,\lambda_2}(0)=0$
using Littelmann's path model.

\subsection{The saturation theorem for $GL(\ell)$}
\label{saturationforgl}

In this subsection we explain Green's idea for proving that for
$GL(\ell)$
$$
Sol({\mathbf Q}{\mathbf 3}) \subset Sol({\mathbf Q}{\mathbf 4}).
$$
\no Let $\K$ be a valued field with a discrete complete valuation $v$
and let $\O\subset \K$ be the subring of elements with nonnegative valuation.
Pick a uniformizer $\pi\in \K$ for the valuation $v$. We will assume that the residue field is finite.

Let $\ul{G}$ be a connected reductive group algebraic group  over $\K$ and set
$G=\ul{G}(\K)$ and let $G^\vee=\ul{G}^\vee(\C)$ be its Langlands dual.
We will assume that $\ul{G}$ is split over $\K$. Recall that
the {\em structure constants}  $m_{\alpha,\beta}(\gamma)$
and $n_{\alpha,\beta}(\gamma)$ for the Hecke ring of
$G$ and the representation ring of $G^{\vee}$ respectively, are given by
$$
c_{\alpha} \cdot c_{\beta} =
\sum_{\ga} m_{\alpha,\beta}(\gamma) c_{\ga},
$$
$$
ch(V_{\alpha})\cdot ch(V_{\beta})=
\sum_{\ga} n_{\alpha,\beta}(\ga) \  ch(V_{\ga}).
$$
We will say that a dominant weight $\mu$ for $GL(\ell)$ is a {\em partition} if
all components of the integer vector $\mu$ are nonnegative.

The starting point is that if
the dominant weights $\alpha,\beta,\gamma$ are partitions then
there is a standard formula \cite[pg. 161]{Macdonald}
that expresses
the structure constant $m_{\alpha,\beta}(\gamma)$ as the number of
finite ${\O}$-module extensions
$$
1\to A \to B \to C\to 1,
$$
where $A\cong \oplus_{i=1}^\ell \mathcal{O}/(\pi^{\alpha_i})$,
$B\cong \oplus_{i=1}^\ell \mathcal{O}/(\pi^{\gamma_i})$ and
$C\cong \oplus_{i=1}^\ell \mathcal{O}/(\pi^{\beta_i})$.

We will say that such an extension is of type $(\alpha, \beta, \gamma)$.
Now for $GL(\ell)$ there is an explicit formula (the Littlewood-
Richardson rule) for the structure constant $n_{\alpha,\beta}(\gamma)$
(in case $\alpha,\beta,\gamma$ are partitions)
as the number of Littlewood--Richardson sequences (of partitions)
of type $(\alpha, \beta, \gamma)$, see \cite[pg. 68, 90]{Macdonald}.
Then the observation of J.~Green is the following

\begin{lem}
Suppose there exists an extension of $\mathcal{O}$-modules of type
$(\alpha,\beta,\gamma)$. Then there exists a Littlewood-Richardson sequence
(of partitions) of  type $(\alpha,\beta,\gamma)$.
\end{lem}

Thus for the group $GL(\ell)$ we find that if $\alpha,\beta,\gamma^*$
are partitions and the structure constant
$m_{\alpha,\beta}(\gamma^*)$ is nonzero then the structure constant
$n_{\alpha,\beta}(\gamma^*)$ is nonzero (the superscript $*$ denotes the
contragredient dominant weight).

We next remove the assumption that $\alpha,\beta,\gamma^*$ are partitions.
To this end let $\lambda$ be the cocharacter of $GL(\ell)$ that
sends $z\in \Q\setminus \{0\}$ to the scalar matrix with diagonal entries equal to $z$.
If $\alpha$ is a cocharacter then there exists $k$ such that
$\lambda^k \cdot \alpha$ is a partition. Now since the matrix
$\lambda(\pi)$ is scalar it is clear that for any $\alpha$
$c_{\alpha} \cdot c_{\lambda} = c_{\alpha \cdot \lambda}$.
It then follows that
$$
m_{\la^k \cdot \alpha, \lambda^k \cdot \beta}
(\lambda^{2k} \cdot \gamma^*) = m_{\alpha,\beta}(\gamma^*).
$$
Now the previous operation of convolving with the function $c_{\lambda_1}$
corresponds, under Langlands' \  duality, to multiplying a character
by $\chi_1$, the character of the determinant representation.
Thus
$$
m_{\alpha,\beta}({\gamma^*})\neq 0 \Leftrightarrow
m_{\lambda^k \cdot \alpha, \la^k \cdot \beta}
(\lambda^{2k} \cdot \gamma^*) \neq 0
\Leftrightarrow
n_{\chi^k \cdot\alpha, \chi^k \cdot \beta}
(\chi^{2k} \cdot \gamma^*) \neq 0 \Leftrightarrow
n_{\alpha,\beta}(\gamma^*) \neq 0.
$$

Finally, by Lemma
\ref{contragredientstructureconstant} we have
$$m_{\alpha,\beta,\gamma}(0) \neq 0 \Leftrightarrow
m_{\alpha,\beta}(\gamma^*) \neq 0.$$ Also it is immediate that
$n_{\alpha,\beta}(\gamma^*) = n_{\alpha,\beta,\gamma}(0)$. Thus,
if  $(\alpha,\beta,\gamma)\in Sol({\mathbf Q}{\mathbf 3})$
then $(\alpha,\beta,\gamma)\in Sol({\mathbf Q}{\mathbf 4})$
as well. By combining this with Theorem
\ref{q4->q3}, we conclude that for $GL(\ell)$ Problems {\bf Q3}
and {\bf Q4} are equivalent. Since we have proved that $1$ is a
saturation factor for the Problem {\bf Q3} it is also a saturation
factor for {\bf Q4}.

\subsection{Computations for the root systems $B_2$ and $G_2$}
\label{semigroup}

Given a root system $R$, the intersection $D_3(X)\cap P(R^\vee)$ is a semigroup whose finite generating set
(a {\em Hilbert basis}) can be computed once we know the explicit {\em stability inequalities} defining the convex cone
$D_3(X)$. We have performed these computations for the root systems $R= B_2, G_2$ using the program
4ti2 (which could be found at http://www.4ti2.de) and the stability inequalities established in \cite{KapovichLeebMillson1}.
Let $\varpi_1, \varpi_2$ be the long and short fundamental coweights for the root system $R$;
below we will use the coordinates $[x,y]$ for the coweight $x\varpi_1+ y\varpi_2$.

\bigskip
{\bf $B_2$ computation.}
The Hilbert basis in the case $R=B_2$ consists of the following 8 triples $(\al,\be,\ga)$
and their permutations under the $S_3$ action:
$$
\left(\left[\begin{array}{c}
0\\
1\end{array}\right],
\left[\begin{array}{c}
0\\
1\end{array}\right],
\left[\begin{array}{c}
0\\
0\end{array}\right]\right), \quad
\left(\left[\begin{array}{c}
1\\
0\end{array}\right],
\left[\begin{array}{c}
1\\
0\end{array}\right],
\left[\begin{array}{c}
0\\
0\end{array}\right]\right),
$$
$$
\left(\left[\begin{array}{c}
1\\
0\end{array}\right],
\left[\begin{array}{c}
0\\
1\end{array}\right],
\left[\begin{array}{c}
0\\
1\end{array}\right]\right),
\quad
\left(\left[\begin{array}{c}
1\\
0\end{array}\right],
\left[\begin{array}{c}
1\\
0\end{array}\right],
\left[\begin{array}{c}
0\\
2\end{array}\right]\right),
$$
$$
\left(\left[\begin{array}{c}
1\\
0\end{array}\right],
\left[\begin{array}{c}
1\\
0\end{array}\right],
\left[\begin{array}{c}
1\\
0\end{array}\right]\right),
\quad
\left(\left[\begin{array}{c}
1\\
0\end{array}\right],
\left[\begin{array}{c}
1\\
0\end{array}\right],
\left[\begin{array}{c}
1\\
1\end{array}\right]\right),
$$
$$
\left(\left[\begin{array}{c}
0\\
1\end{array}\right],
\left[\begin{array}{c}
0\\
1\end{array}\right],
\left[\begin{array}{c}
0\\
1\end{array}\right]\right),
\quad
\left(\left[\begin{array}{c}
1\\
0\end{array}\right],
\left[\begin{array}{c}
1\\
0\end{array}\right],
\left[\begin{array}{c}
0\\
1\end{array}\right]\right).
$$
We note that the first 5 generators are represented by flat triangles contained in an apartment in $X$
and therefore are solutions of {\bf Q4} (and of course of {\bf Q3}).
The last three generators are not solutions of {\bf Q3} since
$\al+\be+\ga\notin Q(R^\vee)$.

However, a direct computation (using Littelmann triangles)
shows that for every generator $(\al,\be,\ga)$ among the last 3 generators
in our list, $(2\al,2\be,2\ga)$ is a solution of {\bf Q4}. Therefore, since the solution set of {\bf Q4}
forms a semigroup, we obtain

\begin{prop}
If $(\al,\be,\ga)$ belongs to $D_3(X)\cap X_*(\ul{T})$ then the triple
$(2\al,2\be,2\ga)$ is a solution of the problem {\bf Q4} for the group $Sp(4)$.
\end{prop}

We next observe that the solution set of the problem {\bf Q3} for the group $Spin(5)$ is not a semigroup. Indeed,
using the coordinates as in the Example \ref{notriangle2},
consider the triples $(\al'=(1,1), \be'=(1,1), \ga'=(1,1))$, $(\al''=(1,1), \be''=(1,1), \ga''=(2,0))$.
As in section \ref{3not4}, the first triple is a solution of the problem {\bf Q3} for $Spin(5)$, the second triple
represents a flat triangle with the vertices $(0,0), (1,1), (2,0)$ in the apartment.
However $(\al'+\al'', \be'+\be'', \ga'+\ga'')=(\al,\be,\ga)$ is the triple from the Example \ref{notriangle2},
and therefore is not a solution of {\bf Q3}.

\bigskip
{\bf $G_2$ computation.} The Hilbert basis in the case $R=G_2$ consists of the following triples
$(\al,\be,\ga)$ and their permutations under the $S_3$ action:
$$
\left(\left[\begin{array}{c}
0\\
1\end{array}\right],
\left[\begin{array}{c}
0\\
1\end{array}\right],
\left[\begin{array}{c}
0\\
0\end{array}\right]\right), \quad
\left(\left[\begin{array}{c}
1\\
0\end{array}\right],
\left[\begin{array}{c}
1\\
0\end{array}\right],
\left[\begin{array}{c}
0\\
0\end{array}\right]\right),
$$
$$
\left(\left[\begin{array}{c}
0\\
1\end{array}\right],
\left[\begin{array}{c}
0\\
1\end{array}\right],
\left[\begin{array}{c}
0\\
1\end{array}\right]\right),  \quad
\left(\left[\begin{array}{c}
1\\
0\end{array}\right],
\left[\begin{array}{c}
1\\
0\end{array}\right],
\left[\begin{array}{c}
1\\
0\end{array}\right]\right),
$$
$$
\left(\left[\begin{array}{c}
1\\
0\end{array}\right],
\left[\begin{array}{c}
1\\
0\end{array}\right],
\left[\begin{array}{c}
0\\
3\end{array}\right]\right),
\quad
\left(\left[\begin{array}{c}
1\\
0\end{array}\right],
\left[\begin{array}{c}
2\\
0\end{array}\right],
\left[\begin{array}{c}
0\\
3\end{array}\right]\right),
$$
$$
\left(\left[\begin{array}{c}
1\\
0\end{array}\right],
\left[\begin{array}{c}
0\\
1\end{array}\right],
\left[\begin{array}{c}
0\\
1\end{array}\right]\right),
\quad
\left(\left[\begin{array}{c}
1\\
0\end{array}\right],
\left[\begin{array}{c}
0\\
1\end{array}\right],
\left[\begin{array}{c}
0\\
2\end{array}\right]\right),
$$
$$
\left(\left[\begin{array}{c}
1\\
0\end{array}\right],
\left[\begin{array}{c}
1\\
0\end{array}\right],
\left[\begin{array}{c}
0\\
2\end{array}\right]\right),
$$
$$
\left(\left[\begin{array}{c}
1\\
0\end{array}\right],
\left[\begin{array}{c}
1\\
0\end{array}\right],
\left[\begin{array}{c}
0\\
1\end{array}\right]\right),
\quad
\left(\left[\begin{array}{c}
1\\
0\end{array}\right],
\left[\begin{array}{c}
1\\
0\end{array}\right],
\left[\begin{array}{c}
1\\
1\end{array}\right]\right).
$$
The first 9 generators are solutions of {\bf Q4} and the last two generators are only solutions of {\bf Q3}.
However one can show (arguing analogously to the Example \ref{notriangle2}) that the sum of the last two generators,
$$
\left(\left[\begin{array}{c}
2\\
0\end{array}\right],
\left[\begin{array}{c}
2\\
0\end{array}\right],
\left[\begin{array}{c}
1\\
2\end{array}\right]\right)
$$
is not a solution of {\bf Q3}. Hence the solution set of {\bf Q3} for the root system $G_2$ is not a semigroup.
On the other hand, a direct computation shows that for both $k=2, k=3$ the triples
$$
k\cdot \left(\left[\begin{array}{c}
1\\
0\end{array}\right],
\left[\begin{array}{c}
1\\
0\end{array}\right],
\left[\begin{array}{c}
0\\
1\end{array}\right]\right),
\quad
k\cdot \left(\left[\begin{array}{c}
1\\
0\end{array}\right],
\left[\begin{array}{c}
1\\
0\end{array}\right],
\left[\begin{array}{c}
1\\
1\end{array}\right]\right)
$$
are solutions of {\bf Q4}. Therefore, since the solution set of {\bf Q4} is a semigroup and
because each natural number $k\ge 2$ has the form $2n+ 3m, n, m\in \N\cup \{0\}$, we get:

\begin{prop}
For each $k\in \N\setminus \{1\}$ the semigroup
$$
k\cdot (D_3(X)\cap P(R^\vee))
$$
is contained in the solution set of {\bf Q4}.
\end{prop}

\section{Appendix: Decomposition of tensor products and Mumford quotients
of products of coadjoint orbits}
 \label{append}

\subsection{The existence of semistable triples and nonzero invariant vectors
in triple tensor products}

In this section we will assume that $G$ is a reductive complex Lie group with
Lie algebra $\mathfrak{g}$ and Weyl group $W$. We let $K$ be a maximal compact
subgroup of $G$ and $T$ be a maximal torus in $G$ and $B$ be a Borel
subgroup containing $T$. We let $\h$ denote the Lie algebra of $T$. We may choose $T$
so that it is preserved by the Cartan involution of $G$ corresponding to $K$.
Let $\g =\k \oplus \p$ denote the Cartan decomposition of the Lie algebra $\g$ of $G$ and
let $\a:= \h\cap \p$.
Then multiplication by $\sqrt{-1}$ interchanges $\p$ and $\k$.
We let $X:= G/K$ be the associated symmetric space.
The choice of $B$ is equivalent to a choice
of positive roots or a positive chamber $\Delta \subset \mathfrak{a}$.
We will be sloppy throughout with the difference between the solution
sets for Problem {\bf Q1} for $G$ and $G^{\vee}$
since these sets are canonically isomorphic (say by using an $Ad\, K$ invariant
metric on $\k$). In the similar fashion we will identify $\k$ with $\k^*$
and $\p$ with $\p^*$.

\medskip
We will use the following notation. Let
$\lambda \in \mathfrak{a}^*$ be a dominant weight for the torus $T$.
Then we let $V_{\lambda}$ be an irreducible $G$--module with the
highest weight $\lambda$ (so $V_{\lambda}$ is unique up to isomorphism).
We let $ \lambda^*$ be the highest weight of the dual representation,
so $ \lambda^*= w_0 ( - \lambda)$ where $w_0$ is the longest element in $W$.

For each fundamental weight $\la$ belonging to $\a^*\subset \p^*$
we define the {\em coadjoint orbit}, $\O_{\la}\subset \k^*$, to be  the
$Ad^*K$--orbit of $\sqrt{-1}\lambda$.

Let $\la \in \mathfrak{a}^*$ be a dominant weight.
Then the coadjoint orbit $\mathcal{O}_{\la}$ as above carries a natural
homogeneous complex structure (see \cite[Chapter 1]{Vogan}) so that
the coadjoint action of $K$ on $\mathcal{O}_{\la}$ extends to
a holomorphic action of $G$. The dominant weight $\la$ defines a
very ample line bundle $\mathbb{L}_{\la}$ over $\mathcal{O}_{\la}$ and the
action of
$G$ on $\mathcal{O}_{\la}$ extends to a holomorphic action
$G \acts \mathbb{L}_{\la}$.

The key point in what follows is the
famous theorem of Borel and Weil below,
Theorem \ref{BW}, see e.g. \cite[Chapter 1] {Vogan}.

\begin{thm}
\label{BW}
 The space of holomorphic sections $\Gamma(\mathcal{O}_{\la}, \mathbb{L}_{\la})$
of $\mathbb{L}_{\la}$ is isomorphic as a $G$--module to the
$G$--module $V_{\la}$.
\end{thm}

Now let $\alpha$ and $\beta$ be dominant weights.
In this section we will discuss the problem of finding the possible
irreducible constituents of tensor products $V_{\alpha} \otimes V_{\beta}$.
Of course this is equivalent to finding for which triples $\alpha$, $\beta$,
$\gamma$, the space of $G$--invariants
$$
(V_{\alpha} \otimes V_{\beta}\otimes
V_{\gamma})^{G} \cong Hom_{G}(V_{\gamma^*}, V_{\alpha} \otimes V_{\beta})$$
is nonzero.

Given dominant weights $\al, \be, \ga$ we define complex $G$--manifold
$$
\mathcal{O}_{\alpha,\beta,\gamma} :=
 \mathcal{O}_{\alpha}\times \mathcal{O}_{\beta} \times \mathcal{O}_{\gamma}.
$$
We also have the outer tensor product
$$\mathbb{L}_{\alpha,\beta,\gamma}:=
\mathbb{L}_{\alpha}
\boxtimes \mathbb{L}_{\beta } \boxtimes \mathbb{L}_{\gamma}$$
which is a very ample $G$--invariant line bundle over
$\mathcal{O}_{\alpha,\beta,\gamma}$.

Since the line bundle $\mathbb{L}_{\alpha,\beta,\gamma}$ is very ample,
it determines a $G$---equivariant holomorphic embedding $F$ of
$\mathcal{O}_{\alpha,\beta,\gamma}$
into $\mathbb{P}((V_{\alpha} \otimes V_{ \beta}\otimes V_{ \gamma})^*)$.
We will use the notation
$$
{\mathcal M}_{\al,\be,\ga}:= F(\mathcal{O}_{\alpha,\beta,\gamma}) //G
$$
for the Mumford quotient associated with this line bundle. See the end of section \ref{Gaussmaps} for further discussion;
note that here we have shortened the notation
${\mathcal M}_{(\al,\be,\ga), sst}(B)$ (in \S \ref{Gaussmaps}) to
${\mathcal M}_{\al,\be,\ga}$.

\begin{lem}
The moduli space of triangles in the infinitesimal symmetric space $\p$ with $\De$--side lengths
$\alpha$, $\beta$, $\gamma$ is canonically homeomorphic to the Mumford quotient
$\mathcal{M}_{\alpha,\beta,\gamma}$ defined above.
\end{lem}
\proof
By a theorem of Kempf and Ness \cite{KempfNess}, the Mumford quotient
$\mathcal{M}_{\alpha,\beta,\gamma}$ is canonically homeomorphic to
the symplectic quotient $\mathcal{O}_{\alpha,\beta,\gamma}//Ad^*(K)$.
It is a standard argument (essentially the formula for the
moment map of the action of the diagonal subgroup of a product) that
$$
\mathcal{O}_{\alpha,\beta,\gamma}//Ad^*(K)=
\{(\lambda,\mu,\nu) \in \mathcal{O}_{\alpha,\beta,\gamma}:
\lambda + \mu + \nu =0\}/Ad^*(K).$$
However the latter is canonically isomorphic to
the moduli space of triangles in $\p$ with the $\De$--side lengths $\alpha,\beta,\gamma\in \De$. \qed

\medskip
Thus
$$
\mathcal{M}_{\alpha,\beta,\gamma}\ne \emptyset \iff (\alpha,\beta,\gamma)
\in D_3(\p) \iff (\alpha,\beta,\gamma)
\in D_3(X).
$$

\bigskip
The main goal of this section is to prove
Theorem \ref{vectorsandngons}, that relates the existence of
semistable triples in $\mathcal{O}_{\alpha,\beta,\gamma}$
and the existence of nonzero invariants in triple tensor products of irreducible
representations of $G$.

\begin{thm}
\label{vectorsandngons}
For each triple of dominant weights $\al, \be, \ga\in \mathfrak{a}^*$ the following are
equivalent:

1. There exists a positive integer $k>0$ such that
$(V_{k\alpha} \otimes V_{k \beta}\otimes V_{k\gamma})^{G}
\neq {0}$.

2. $\mathcal{M}_{\alpha,\beta,\gamma} \neq \emptyset$, i.e. there exists
a semistable point in $\mathcal{O}_{\alpha,\beta,\gamma}$.
\end{thm}

\begin{rem}
Theorem \ref{vectorsandngons} implies that
if $(V_{ \alpha} \otimes V_{ \beta}\otimes V_{\gamma})^{G}\neq {0}$
then there exists a weighted semistable configuration on $\tits X$ of type
$(\alpha, \beta, \gamma)$.
However the converse is false as we have seen in section \ref{3not4}.
\end{rem}

Theorem \ref{vectorsandngons}  will follow from the next two lemmas.

Let $A$ be the graded ring associated to the projective variety
 $F(\mathcal{O}_{\alpha,\beta,\gamma})$. Then by definition
$A^{(k)}$ is the set of restrictions to
$F(\mathcal{O}_{\alpha,\beta,\gamma})$
of the homogeneous polynomials of degree $k$ in the projective
coordinates given by a basis of
$$
V_{ \alpha} \otimes V_{ \beta}\otimes V_\gamma=
((V_{ \alpha} \otimes V_{ \beta}\otimes V_{ \gamma})^*)^*.$$
Thus
$$
A^{(k)}= \Gamma (\mathcal{O}_{\alpha,\beta,\gamma},
\mathbb{L}_{\alpha, \beta, \gamma})^{\otimes k}/I^{(k)}.$$
Here $I^{(k)}$ is the degree $k$--summand of the graded ideal of polynomials in the sections
of $\mathbb{L}_{\alpha,\beta ,\gamma}$ that vanish on $\mathcal{O}_{\alpha,\beta,\gamma}$.
Equivalently, $I^{(k)}$ is the degree $k$ component of the ideal $I$ of
polynomials that vanish on $F(\mathcal{O}_{\alpha,\beta,\gamma})$
in the ring of  polynomials on \newline
$\mathbb{P}((V_{ \alpha} \otimes V_{ \beta}\otimes V_{ \gamma})^*)$.

Now another definition of the Mumford quotient is  $Proj (A^{G})$
where $A^{G}$ is the subring of $G$--invariants  of the graded ring $A$.
We define another graded ring $R$ by
$$
R = \bigoplus_{k=0}^{\infty}\Gamma(\mathcal{O}_{\alpha,\beta,\gamma},
(\mathbb{L}_{\alpha}\boxtimes \mathbb{L}_{\beta } \boxtimes
\mathbb{L}_{\gamma})^{\otimes k}).$$
Below and in what follows, if $U_1, U_2$ are $G_1, G_2$--modules, then
$U_1\boxtimes U_2$ denotes the $G_1\times G_2$--module with the underlying
vector space equal to the tensor product $U_1\otimes U_2$.

We define a graded $G \times G \times G$--module $R^{\prime}$ by
$$
 R^{\prime} =  \bigoplus_{k=0}^{\infty}(V_{k \alpha} \boxtimes
V_{k \beta}\boxtimes V_{k \gamma}).$$
 We will abuse notation and use $R^{\prime}$ to denote the restriction
 of the previous module to the diagonal in $G \times G \times G$.
Then, we have an isomorphism of $G$--modules
$$
R^{\prime} =  \bigoplus_{k=0}^{\infty}(V_{k \alpha} \otimes
V_{k \beta}\otimes V_{k \gamma}).$$

\begin{lem}
There is a canonical isomorphism of graded $(G \times G \times G)$--modules\newline
$R\cong R^{\prime}$.
\end{lem}
\proof
The lemma follows from the following three equations. First, let $M_1$ and $M_2$ be
complex manifolds and $\mathbb{L}_1$, resp. $\mathbb{L}_2$, be a holomorphic line
bundle over $M_1$, resp. $M_2$. First we have

$$
(\mathbb{L}_1 \boxtimes \mathbb{L}_2)^{\otimes k} =
\mathbb{L}_1^{\otimes k} \boxtimes \mathbb{L}_2^{\otimes k}.$$

Next, we have
$$
\Gamma(M_1 \times M_2, \mathbb{L}_1 \boxtimes \mathbb{L}_2) \cong
\Gamma(M_1, \mathbb{L}_1)\otimes \Gamma(M_2, \mathbb{L}_2).
$$
Finally, for a natural number $k$, a dominant weight $\al$ and
the corresponding orbit $\mathcal{O}_{\alpha}$
we have
$$
\Gamma(\mathcal{O}_{\alpha}, \mathbb{L}_{\al}^{\otimes k})\cong
\Gamma(\mathcal{O}_{\alpha}, \mathbb{L}_{k\al}). \qed$$

Next we show that $A$ and $R$ are isomorphic. This is the exceptional
feature of the homogeneous situation.
The key point, the surjectivity of the natural map below,
was pointed out to us by Lawrence Ein.

\begin{lem}
There is a natural $G \times G \times G$--equivariant isomorphism
 from the  graded ring $A$ to  the graded ring $R$.
\end{lem}
\proof
To obtain the desired map from the graded ring $A$ to the graded ring $R$
we observe that for any complex manifold $M$ and holomorphic line bundle
$\mathbb{L}$ over $M$ there is a natural map (usually not onto) from
$\Gamma(M , \mathbb{L})^{\otimes k}$ to
$\Gamma(M , \mathbb{L}^{\otimes k})$. Hence there is a natural map from
$A^{(k)}$ to $R^{(k)}$. The exceptional feature here is that this natural
map is {\em onto}.

To see this we note that by
the theorem of Borel and Weil the action of $G \times G \times
G$ on $R^{(k)}$
is irreducible. But the image of $A^{(k)}$ in $R^{(k)}$ is an invariant
subspace. Hence $A^{(k)}$ maps onto $R^{(k)}$. Clearly the map is injective
(because we have divided by the ideal $I^{(k)})$.
\qed

Theorem \ref{vectorsandngons} follows by taking $G$--invariants from the
$G$--isomorphism $A\cong R^{\prime}$. Indeed,
there exists a triangle with $\De$--side lengths $\alpha$,$\beta$,$\gamma
\Leftrightarrow$ the Mumford quotient $\mathcal{M}_{\alpha,\beta,\gamma}$ is
nonempty $\Leftrightarrow (\bigoplus_{k=1}^{\infty}A^{(k)})^{G} \neq 0 \Leftrightarrow
\bigoplus_{k=1}^{\infty}(V_{k\alpha} \otimes V_{k \beta}\otimes
V_{k \gamma})^{G} \neq 0$. \qed

\begin{rem}
The reason that the existence of a nonzero $G$--invariant does not follow
from the existence of a triangle is that
the lowest degree $G$--invariant might not have  degree $1$, that is, it is
possible that
$$
(R^{(1)})^{G} = \{0\} \ \hbox{but} \ (R^{(k)})^{G} \neq 0
\ \hbox{~ for some $k > 1$}.
$$
\end{rem}

\subsection{The semigroups of solutions to Problems {\bf Q1} and {\bf Q4}}

We will use $S_G$ to denote the
semigroup of dominant characters of $T$ henceforth. We begin this section by
noting that the set of triples of
dominant characters $(\alpha,\beta,\gamma)$ that belong to $D_3(X)$
is a subsemigroup $S_{triangle} \subset S_G^3$ (since it is
determined by a system of homogeneous linear inequalities). Moreover,
since these inequalities have integral coefficients, the subsemigroup
$S_{triangle}$ is finitely generated.   We set
$$
S_{rep}:=Sol({\bf Q4}, G)\subset D_3(X)\cap S_G^3,
$$
i.e. the set of triples $(\al,\be,\ga)$ for which
$$
(V_\al \otimes V_\be \otimes V_\ga)^G\ne 0.
$$

We recall that if $S_1 \subset S_2$ is an inclusion of semigroups then
the {\em saturation} of $S_1$ in $S_2$ is the semigroup of elements of $s \in S_2$ such
that for some positive integer $n$ we have $ns \in S_1$.

\begin{thm}
\label{semig}
\begin{enumerate}
\item The set $S_{rep}$ is a subsemigroup
of the semigroup $S_G^3$ (and of $S_{triangle}$).
\item The saturation of $S_{rep}$ in $S_G^3$ is $S_{triangle}$.
\end{enumerate}
\end{thm}

First we need a general lemma.
Let $G_1$ be a complex reductive group, $B_1$ be a Borel subgroup of $G_1$  and
$\lambda$ and $\mu$ be dominant characters
(as in the beginning of the appendix). We recall that the irreducible representation
$V_{\lambda + \mu}$ is always an irreducible constituent of multiplicity 1 in the tensor product
$V_{\la} \otimes V_{\mu}$. Let $\pi:V_{\la} \otimes V_{\mu}\to V_{\lambda + \mu}$
be the $G_1$--equivariant projection.

\begin{lem}
Let $v_1 \in V_{\la}$ and $v_2 \in V_{\mu}$ be nonzero vectors. Then
$$\pi(v_1 \otimes v_2) \neq 0.$$
\end{lem}
\proof Let $M = G_1/B_1$. We apply the Borel-Weil
Theorem  to obtain
$V_{\lambda} = \Gamma(M, \mathbb{L}_{\la})$, $V_{\mu} = \Gamma(M, \mathbb{L}_{\mu})$
and $V_{\lambda+ \mu} = \Gamma(M, \mathbb{L}_{\la+\mu})$. Then $v_1$ corresponds to
a section $s_1$ and $v_2$ corresponds to a section $s_2$. Hence $\pi(v_1 \otimes v_2)$
corresponds to the product section $s_1\cdot s_2$ of the product line bundle
$\mathbb{L}_{\la+\mu}$. But since $M$ is irreducible the product section cannot be zero.
\qed

 Now we can prove the theorem.

\proof Let $(\alpha_i, \beta_i,\gamma_i) \in S_{rep}$ and let $v_i$
be a nonzero $G$--invariant vector in
$$
V_{\al_i} \otimes V_{\be_i} \otimes V_{\ga_i}, i=1, 2.$$
 We take $G_1:= G\times G\times G$ and $B_1:= B\times B\times B$, where $B\subset G$
is a Borel subgroup as before; we have the $G_1$--modules
$$
V_\la:= V_{\al_1}\boxtimes V_{\be_1}\boxtimes V_{\ga_1},
V_\mu:= V_{\al_2}\boxtimes V_{\be_2}\boxtimes V_{\ga_2}.
$$
Then
$$
V_\la\otimes V_\mu\cong (V_{\al_1}\otimes V_{\al_2})\boxtimes (V_{\be_1}\otimes V_{\be_2})
\boxtimes (V_{\ga_1}\otimes V_{\ga_2}).
$$
By the above lemma,  the vector $\pi(v_1 \otimes v_2)$
is nonzero. Since it is clearly $G$--invariant,
$$
(V_{\al_1+\al_2}\boxtimes V_{\be_1+\be_2}\boxtimes V_{\ga_1+\ga_2})^G\ne 0.
$$
Hence $(\al_1+\al_2, \be_1+\be_2, \ga_1+\ga_2)\in S_G$ and the
first statement of Theorem follows.

 The second statement follows from Theorem \ref{vectorsandngons}.
 \qed

An alternative proof of Part 1 of Theorem \ref{semig} follows from the description of $S_{rep}$ given in  \cite{BZ}
by Berenstein and Zelevinsky.

\begin{thm}
\label{uniformsaturation}
There exists a (nonzero) $k\in \N$ depending only on the group $G$ so that
$k\cdot S_{triangle}\subset S_{rep}$.
\end{thm}
\proof By Gordan's lemma, see \cite[Proposition 1, Page 12]{Fultontoric}, the semigroup $S_{triangle}$
is finitely generated. Choose a finite set of generators $s_1,...,s_m$ of this
 semigroup. For each generator $s_i$ there exists a positive integer $k_i$ such that
$$
k_i s_i\in S_{rep},
$$
see Theorem \ref{vectorsandngons}. Then take $k=LCM(k_1,...,k_m)$. \qed

Therefore, since $S_{rep}$ contains finitely generated semigroup $k\cdot S_{triangle}$
as a subsemigroup of finite index, we get

\begin{cor}
$S_{rep}$ is finitely generated.
\end{cor}

\bigskip
\noindent Michael Kapovich: \newline
Department of Mathematics, \newline
University of California, \newline
Davis, CA 95616, USA
\newline
 kapovich$@$math.ucdavis.edu

\smallskip
\noindent
Bernhard Leeb: \newline
Mathematisches Institut,\newline
Universit\"at M\"unchen, \newline
Theresienstrasse 39, \newline
D-80333 M\"unchen, Germany\newline
leeb$@$mathematik.uni-muenchen.de

\smallskip
\noindent John J. Millson: \newline
Department of Mathematics, \newline
University of Maryland, \newline
College Park, MD 20742, USA
\newline
 jjm$@$math.umd.edu

\end{document}